\documentclass[10pt]{article}
\usepackage{amsmath}
\usepackage{amssymb}
\usepackage{amsthm}
\font\smallit=cmti10
\font\smalltt=cmtt10

\makeatletter

\renewcommand\section{\@startsection {section}{1}{\z@}
{-30pt \@plus -1ex \@minus -.2ex}
{2.3ex \@plus.2ex}
{\normalfont\normalsize\bfseries\boldmath}}

\renewcommand\subsection{\@startsection{subsection}{2}{\z@}
{-3.25ex\@plus -1ex \@minus -.2ex}
{1.5ex \@plus .2ex}
{\normalfont\normalsize\bfseries\boldmath}}

\renewcommand{\@seccntformat}[1]{\csname the#1\endcsname. }

\newcommand\blfootnote[1]{
\begingroup
\renewcommand\thefootnote{}\footnote{#1}
\addtocounter{footnote}{-1}
\endgroup
}

\makeatother

\newtheorem{theorem}{Theorem}
\newtheorem{lemma}[theorem]{Lemma}
\numberwithin{theorem}{section} 
\newtheorem{corollary}[theorem]{Corollary}
\theoremstyle{definition}
\newtheorem{definition}[theorem]{Definition}
\newtheorem*{example}{Example}
\theoremstyle{remark} 
\newtheorem*{remark}{Remark}

\usepackage{graphicx} 
\usepackage{amsfonts}
\usepackage[section]{placeins}
\usepackage{comment}
\usepackage{amsmath}
\usepackage{tabularx}
\usepackage{varwidth}
\usepackage{verbatim}
\usepackage{hyperref}

\def\Q{{\mathbb Q}}

\def\C{{\mathcal C}}

\def\case#1{\noindent {\tt Case~#1:}}

\begin{document}
\begin{figure}
\vspace*{-62pt}
\vspace{.75in}
\end{figure}
\vspace*{-45pt}
\leftline{\smalltt\#A1} \vskip -12.5pt
\centerline{\smalltt  INTEGERS 26 (2026)}
\vskip 3pt \hrule																													
\begin{center}
\uppercase{\bf \boldmath Tilings of an isosceles triangle}\blfootnote{DOI: 
}
\vskip 20pt
{\bf Michael Beeson}\\
{\smallit  San Jos\'e State University, San Jos\'e, California}\\ 

{\tt beesonpublic@gmail.com}\\ 
\end{center}
\vskip 20pt
\vskip 30pt

\centerline{\bf Abstract}
\noindent
An {\em $N$-tiling} of triangle $ABC$ is a way 
to cut $ABC$ into $N$ congruent smaller triangles.  The 
smaller triangle is the {\em tile}.  When $ABC$ is isosceles
with base angles $\alpha$, and not equilateral, 
there are only four possible tiles (aside from a tile similar to 
$ABC$):  a right-angled tile with 
one angle $\alpha$,  a tile with angles $(\alpha, \beta, 2\alpha)$,
 a tile with angles $(\alpha,\beta, 2\pi/3)$, or a tile
with angles satisfying $3\alpha + 2\beta = \pi$ or $2\alpha+3\beta=\pi$ (and
in all but the first case, with $\alpha$ not a 
rational multiple of $\pi$).   We study the first three cases in this paper,
proving several theorems about what values of $N$ can or cannot occur.

\pagestyle{myheadings}
\markright{\smalltt INTEGERS: 26 (2026)\hfill}
\thispagestyle{empty}
\baselineskip=12.875pt
\vskip 30pt

\section{Introduction}
An \emph{$N$-tiling} of triangle $ABC$ by triangle $T$ is a way
of writing $ABC$ as a union of $N$ triangles congruent to $T$,
overlapping only at their boundaries.
The triangle $T$ is the \emph{tile}. 
When $ABC$ is isosceles
with base angles $\alpha$, and not equilateral, 
there are only four possible tiles (aside from a tile similar to 
$ABC$):  a right-angled tile with 
one angle $\alpha$,  a tile with angles $(\alpha, \beta, 2\alpha)$,
 a tile with angles $(\alpha,\beta, 2\pi/3)$, or a tile
with angles satisfying $3\alpha + 2\beta = \pi$ or $2\alpha+3\beta=\pi$ (and
in all but the first case, with $\alpha$ not a 
rational multiple of $\pi$). 

In this article, we studied each of these families of tilings, and obtained
the following results.
For tilings by a right triangle, $N$ has to be a square, or an even sum of squares,
or six times a square; in particular 
it cannot be a prime with $N \equiv 3 \pmod 4$
and all these possibilities actually occur.  We prove
that unless $ABC$ is a right isosceles triangle, $N$
has to be even.

 For tilings by  
$(\alpha, \beta, 2\alpha)$, we show that the tile is necessarily 
rational (the ratios of its sides are rational), and we give
 a necessary condition for the existence of a tiling.  This
 condition implies that when an isosceles and not equilateral
$ABC$ is $N$-tiled by such a tile, 
$N$ cannot be a prime number, or even squarefree.

In the last case, when the tile
has a 120 degree angle, we also prove that the tile must be
rational, and find a necessary condition for the existence
of a tiling.  That condition rules
out $N < 33$.  The smallest known such tiling has $N = 2673$.
Bonfioli \cite{bonfioli2026} has shown that $N$ cannot
be prime in such a tiling.

The proofs in this paper have been checked by computer, in the Lean proof
assistant, except for the constructions of tilings; see
Appendix~\ref{section:lean}.

Our results fit into a larger research program, begun by Laczkovich
\cite{laczkovich1995}.  Laczkovich studied the possible shapes of
tiles and triangles that can possibly be used in tilings, and obtained
results that will be described below.  The reader who is new to 
the subject may want to see examples of $N$-tilings for various 
shapes of $ABC$; such pictures can be found in \cite{beeson-noseven}.
Here we give only examples relevant to the case of $ABC$ isosceles.

\begin{figure}[ht]
\begin{center}
\includegraphics[width=0.7\textwidth]{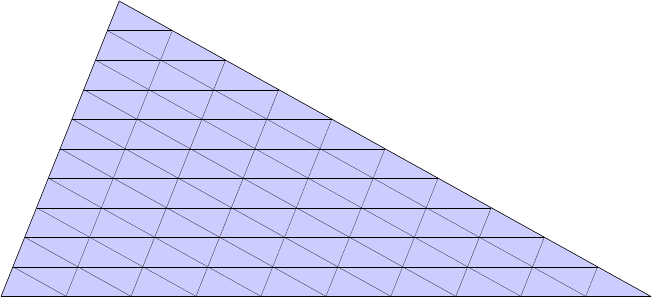}
\end{center}
\caption{A quadratic tiling of an arbitrary triangle.}
\label{figure:bigquadratic}
\end{figure}

First we point out that 
{\em any} triangle can be decomposed into $n^2$ congruent triangles by  drawing $n-1$ equally spaced lines parallel to each of the three sides of 
the triangle, as illustrated in Figure~\ref{figure:bigquadratic}.
Moreover, the large (tiled) triangle is similar to the 
small triangle (the ``tile'').  We call such a tiling a {\em quadratic tiling.} 

 It follows that if we have a tiling of a triangle $ABC$ into $N$ congruent triangles, and $m$ is any 
integer,  we can tile $ABC$ into $Nm^2$ triangles by subdividing the first tiling, replacing each of the $N$ triangles by $m^2$ smaller ones.
Hence, the set of $N$ for which an $N$-tiling of some triangle exists is closed under multiplication by squares.

\begin{figure}[ht]
\includegraphics[width=0.45\textwidth]{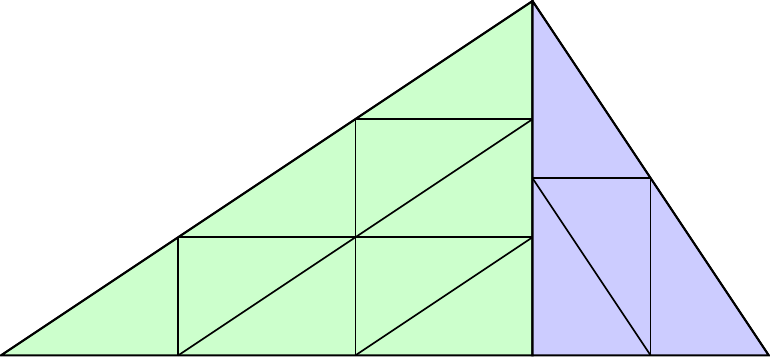} \quad
\includegraphics[width=0.45\textwidth]{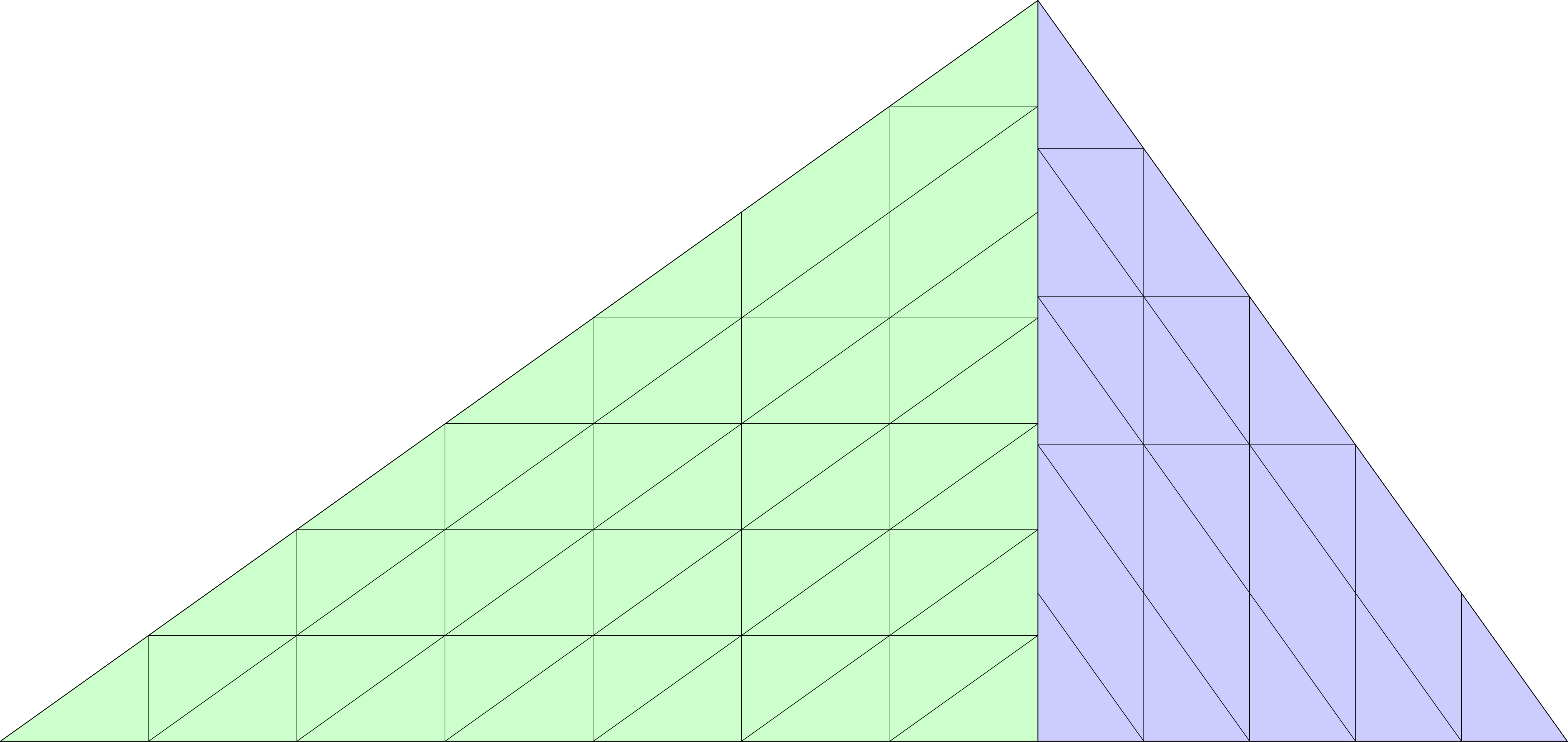}
\caption{Biquadratic tilings with $N = 13 = 3^2 + 2^2$ and $N=74 = 5^2 + 7^2$}.
\label{figure:biquadratic}
\end{figure}

Sometimes it is possible to combine two quadratic tilings (using the same 
tile) into a single tiling, as shown in Figure~\ref{figure:biquadratic}.
We will explain how these tilings are constructed.  We start with 
a big right triangle resting on its hypotenuse, and divide it into 
two right triangles by an altitude.  Then we quadratically tile each 
of those triangles.  The trick is to choose the dimensions in such a 
way that the same tile can be used throughout.  If that can be 
done then evidently $N$, the total number of tiles, will be the
sum of two squares, $N = n^2 + m^2$,  one square for each of the two 
quadratic tilings.  On the other hand, if we start with an $N$ of 
that form,  and we choose the 
tile to be an $n$ by $m$ right triangle, then we can construct 
such a tiling.  We call these tilings {\em biquadratic}. 
More generally, a {\em biquadratic tiling} of triangle $ABC$ is one in which $ABC$ has a right angle at $C$, and can be divided by an altitude 
from $C$ to $AB$ into two triangles, each similar to $ABC$, which can be tiled respectively by $n^2$ and $m^2$ copies of a triangle similar to $ABC$. 
A larger biquadratic tiling, with  $n=5$ and $m=7$ and hence $N=74$, is shown on the right side of  Figure~\ref{figure:biquadratic}.

If the original triangle $ABC$ is chosen to be isosceles, and is 
then quadratically tiled, 
then each of the $n^2$ triangles can be divided in half by an altitude;  hence any isosceles triangle can be decomposed into $2n^2$ congruent 
triangles.   If the original triangle is equilateral,  then it can be first decomposed into $n^2$ equilateral triangles, and then these 
triangles can be decomposed into 3 or 6 triangles each,  showing that any equilateral triangle can be decomposed into $3n^2$ or $6n^2$
congruent triangles. 
For example, we can 12-tile an equilateral triangle in two different ways,  starting with a 3-tiling and then subdividing each triangle 
into 4 triangles (``subdividing by 4''),  or starting with a 4-tiling and then subdividing by 3.

\begin{figure}[ht]    
\begin{center}
\includegraphics[width=0.35\textwidth]{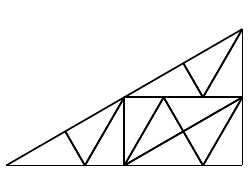}
\hskip 1cm
\includegraphics[width=0.35\textwidth]{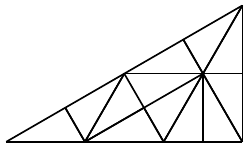}
\end{center}
\caption{Two 12-tilings.}
\label{figure:12-tilings}
\end{figure}

\begin{figure}[ht]    
\begin{center}
\includegraphics[width=0.35\textwidth]{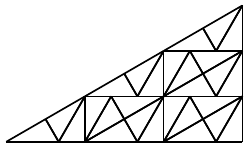}
\hskip 1cm
\includegraphics[width=0.35\textwidth]{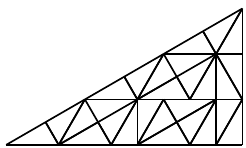}
\end{center}
\caption{Two 27-tilings.}
\label{figure:27-tilings}
\end{figure}

There is another family of $N$-tilings, in which $N$ is of the form $3m^2$, and both the tile and the tiled triangle are 30-60-90 triangles.  We call these the {\em triple-square} tilings. 
The case $m=2$ makes $N=12$.  There are two ways to 12-tile a 30-60-90 triangle with a 30-60-90 triangle.  
One is to first quadratically 4-tile it, and then subtile the four triangles with the 3-tiling of Figure 1.  This produces the first 
12-tiling in Figure~\ref{figure:12-tilings}.  Somewhat surprisingly, there is another way to tile the same triangle with the same 12 tiles, also shown in Figure~\ref{figure:12-tilings}.
The next member of this family is $m=3$, which makes $N=27$.  
Two 27-tilings are shown in Figure~\ref{figure:27-tilings}.
 Similarly, there are two 48-tilings (not shown).

 Whenever there is an $N$-tiling of the right triangle $ABM$, there is a $2N$-tiling of
the isosceles triangle $ABC$.  Using the biquadratic tilings (see 
Figure~\ref{figure:biquadratic}) and  triple-square tilings (see Figure~\ref{figure:12-tilings} and Figure~\ref{figure:27-tilings}),
we can produce $2N$-tilings when $N$ is a sum of squares or three times a sum of squares.  We call these tilings {\em double biquadratic}
and {\em hexquadratic}, respectively. 
  For example,  one has two 26-tilings and two 148-tilings, obtained by reflecting the
two tilings in Figure~\ref{figure:biquadratic}
about either of the 
sides of the triangles shown in those figures;  and one has 24-tilings and 54-tilings obtained from  
Figure~\ref{figure:12-tilings} and Figure~\ref{figure:27-tilings}.
In the latter two cases, $ABC$ will be equilateral if the reflection is through the longer leg.  

In the case when the sides of the tile $T$ form a Pythagorean triple $n^2 + m^2 = k^2$ and  $N = 2k^2$, then we can tile one 
half of  $ABC$ with a quadratic tiling and the other half with a biquadratic tiling.  The smallest example is when the tile
has sides 3, 4, and 5, and $N = 50$.  See Figure~\ref{figure:severalisosceles2}.
 One half is 25-tiled quadratically, and the other half is divided into two smaller
right triangles which are 9-tiled and 16-tiled quadratically.  This shows that the tiling of $ABC$ does not have to be 
symmetric about the altitude.
\FloatBarrier 

As we shall see below, the work of Laczkovich implies that there are 
only four possible shapes of the tile:  right-angled, $\gamma = 2\alpha$,
$\gamma = 2\pi/3$, and $3\alpha + 2\beta = \pi$.  The last case is 
taken up in another paper, since the techniques apply also to tilings of 
non-isosceles triangles $ABC$ with $3\alpha+2\beta = \pi$.  The first
three are studied in this paper.  We obtain, in the second and third case,
necessary conditions on $N$, but not necessary and sufficient conditions.
In the case of a right-angled tile, our conditions are necessary 
and sufficient.   

The oldest 
examples of tilings of isosceles $ABC$ with $\alpha \neq \frac \pi 2$
have $N$ even. We prove in this paper that $N$ must be even when the tile is 
right-angled, unless $ABC$ is a right isosceles triangle. 
 That is not true for all shapes of the tile, as can be seen by the tilings
exhibited and discussed in Section \ref{section:herdt}.

\subsection{Definitions and Notation}
We first note that this paper is about triangles $ABC$ that are 
isosceles and not equilateral.  Let that be understood; then for
the rest of this paper, ``isosceles'' means ``exactly two sides
are equal.''%
\footnote{\,That was Euclid's definition of ``isosceles.''}  

We give a mathematically precise definition of ``tiling'' and fix some terminology and notation.
Given a triangle $T$ and a larger triangle $ABC$,  a {\em tiling} of triangle $ABC$ by triangle $T$ is 
a set of triangles $T_1,\ldots,T_n$ congruent to $T$, whose interiors are disjoint, and the closure of whose union is triangle $ABC$. 

Let $a$, $b$, and $c$ be the sides of the tile $T$, and angles $\alpha$, $\beta$, and $\gamma$ be the angles 
opposite sides $a$, $b$, and $c$.
The letter ``$N$'' will always be used for the number of triangles used in the tiling. 
 An {\em $N$-tiling} of $ABC$ is a tiling that uses
$N$ copies of some triangle $T$.  
The meanings of $N$, $\alpha$, $\beta$, $\gamma$, $a$, $b$, $c$, $A$, $B$, and $C$ will be fixed throughout this paper.  
We do not assume $\alpha \le \beta$ in general;  although that may
sometimes be justified by symmetry, we often will consider some 
equation such as $3\alpha + 2\beta = \pi$, in which case we do not 
want to assume $\alpha \le \beta$.

\section{History}

Above we exhibited quadratic and biquadratic tilings
in which the tile is similar to $ABC$.  There 
are hexagonal tilings, not exhibited in this paper, but see
\cite{beeson-noseven} for pictures.   These
involve $N$ being square, a sum of two squares, or three times a square. 
The biquadratic tilings were known in 1964, when the 
paper \cite{golomb} was published.  This is the earliest paper on 
the subject of which we are aware.%
\footnote{\,The simplest hexagonal tiling
is attributed to Major MacMahon (1921) in the notes accompanying a 
plastic toy purchased at an AMS meeting in 2012.}
 Snover
{\em et.\,al.}\,\cite{snover1991} took up the challenge of showing that 
these are the only possible values of $N$, when the tile is similar to $ABC$.
The following theorem completely answers the question, ``for which $N$ does there exist an $N$-tiling
in which the tile is similar to the tiled triangle?''
\begin{theorem} [Snover {\em et.\,al.}\,\cite{snover1991}] 
\label{theorem:snover} 
 Suppose $ABC$ is $N$-tiled by tile $T$ similar to $ABC$.
If $N$ is not a square, then $T$ and $ABC$ are right triangles.  Then either

(i) $N$ is three times
a square and $T$ is a 30-60-90 triangle, or 

(ii)  $N$ is a sum of squares $e^2 + f^2$,  the right angle of $ABC$ is split 
by the tiling, and the acute angles of $ABC$ have rational tangents $e/f$ and $f/e$,

\noindent
  and these two alternatives are 
mutually exclusive. 
\end{theorem}

Soifer's book \cite{soifer} appeared in 1990, with a second
edition in 2009.  He considered two ``Grand Problems'': for which 
$N$ can {\em every} triangle be $N$-tiled, and for which $N$ can 
{\em every} triangle be dissected into similar, but not necessarily
congruent triangles.  
The 2009 edition has an added chapter in 
which the biquadratic tilings and a theorem of Laczkovich occur.

Miklos Laczkovich published six papers
\cite{laczkovich2010,laczkovich1990,laczkovich1995,laczkovich1998,laczkovich2012,laczkovich-szekeres1995}
on triangle and polygon tilings.  According to Soifer, the 1995 paper
was submitted in 1992.  Laczkovich, like Soifer, 
studied dissecting a triangle into smaller {\em similar} triangles,
not {\em congruent} triangles as we require here.  
If those similar triangles are rational (i.e., if the ratios of their sides
are rational) then if we divide each of them into small enough 
quadratic subtilings, we can achieve an $N$-tiling into {\em congruent}
triangles, but of course $N$ may be large.  Laczkovich focused primarily  
  on the shapes of $ABC$ (or more 
 generally, convex polygons) and 
of the tile.  His theorems 
 give us an exhaustive list of the possible shapes of $ABC$
and the tile, which we will need in our proof that there is no 7-tiling.
This list can be found in Section~\ref{section:laczkovich}.  
However, his
theorem published in the last chapter of \cite{soifer} does mention
$N$.  It states that given an integer $k$,  there exists an $N$-tiling for 
some $N$ whose square-free part is $k$.   

\section{Laczkovich}\label{section:laczkovich}

A basic fact is that, apart from a small number of cases that can be 
explicitly enumerated, if there is an $N$-tiling of $ABC$ by a tile 
with angles $(\alpha, \beta, \gamma)$, then
 the angles $\alpha$ and $\beta$ are not rational 
multiples of $\pi$.  This theorem follows from Theorems~4.1, 5.1, 
and 5.3 of \cite{laczkovich1995}.  Laczkovich is dealing with 
a more general situation,  tiling an arbitrary triangle by tiles
that are only required to be similar, not congruent. We extract
the following theorem from his results by specializing to isosceles
$ABC$ and congruent tiles.

\begin{table}[ht]
\begin{center}
\setlength{\extrarowheight}{.5em}
\begin{tabular}{rr}
$ABC$  &  the tile     \\
\hline
$(\beta,\beta,\pi- 2\beta)$  &  similar to $ABC$ \\
$(\beta,\beta,2\alpha)$ &  $\gamma = \pi/2$ \\
$(\alpha,\alpha,\pi-2\alpha)$  & $\gamma = 2\alpha$ \\
$(\alpha,\alpha,\pi-2\alpha)$ & $\gamma = 2\pi/3$ \\
$(\alpha,\alpha,\alpha+2\beta)$ &$3\alpha + 2\beta = \pi$ \\
$(\beta,\beta,3\alpha)$ &$3\alpha + 2\beta = \pi$ \\
$(\alpha+\beta,\alpha+\beta,\alpha)$ &$3\alpha + 2\beta = \pi$ 
\end{tabular}
\end{center}
\caption{Possible tilings of isosceles triangles, according to Laczkovich.}
\label{table:isosceles}
\end{table}

\begin{theorem}[Laczkovich \cite{laczkovich1995}]
\label{theorem:laczkovich-isosceles}
Let isosceles (and not equilateral) triangle $ABC$ be $N$-tiled by a tile with
angles $(\alpha,\beta,\gamma)$.  Then the possible shapes of 
$ABC$ and the tile are given by Table~\ref{table:isosceles}.
 In the table, 
the triples giving the angles of the tile are $(\alpha, \beta, \gamma)$
after a suitable permutation, i.e., they are unordered triples. 
If the angles of the tile are commensurable, then one of the first two lines applies. 
Otherwise $\alpha$ is not a rational multiple of $\pi$.
\end{theorem}

\begin{remark} 
For example, in the second line of the table, we do not list 
separately $(\alpha,\alpha,2\beta)$, as that is already covered 
by the entry $(\beta,\beta,2\alpha)$, and the fact that we do not 
assume $\alpha < \beta$.
\end{remark}

\begin{proof} [Proof of Theorem~\ref{theorem:laczkovich-isosceles}]
This theorem is proved in \cite{laczkovich1995},
but it is not stated in quite this way; therefore  
 we spell out in detail how this 
statement follows from theorems explicitly stated in \cite{laczkovich1995}.
Let isosceles (and not equilateral) triangle $ABC$ be $N$-tiled by a tile with
angles $(\alpha,\beta,\gamma)$.  Then either all three angles are
rational multiples of $\pi$, or not.  

\case{1} Not all of the angles are rational
multiples of $\pi$.  Then by Theorem~4.1 of \cite{laczkovich1995},
where $T$ in that paper is our $ABC$,  one of six cases holds.
Cases (i), (ii), and (iv) are the first three lines of our table.
Case (iii) says $ABC$ is equilateral, which we have ruled out by 
hypothesis.  Case (v) says $3\alpha + 2\beta = \pi$ and the base 
angles of isosceles $ABC$ must be $\alpha$ or $\beta$ or $\alpha + \beta$, by 
Theorem~2.4 of \cite{laczkovich1995}; so that is lines 5, 6, and 7
of our table. Finally, case~(vi) concerns the shape 
$(\alpha, \alpha,\alpha+3\beta)$  with $\gamma = 2\pi/3$,  which
is another way of writing line~4 of the table, since if 
$\gamma = 2\pi/3$, then $\alpha+3\beta = \pi-2\alpha$.

\case{2} All three angles are rational multiples of 
$\pi$.  Then Theorem~5.1 of \cite{laczkovich1995} applies.  That theorem
is about dissections into similar (rather than congruent) triangles,
and according to the subsequent Theorem~5.3, the last four cases 
(cases (iv), (v), (vi), and (vii))
in Theorem~5.1 cannot hold for dissections into congruent triangles. 
Cases (i) and (ii) are the first lines of our table.  Case (iii) requires
$ABC$ equilateral, which we have ruled out by hypothesis.  Case (iv)
has $ABC$ a right triangle with one angle $\pi/6$, which is not isosceles
and hence irrelevant here.  
\end{proof}

We note in passing the following immediate consequence of 
Laczkovich's theorem:  If an isosceles triangle $ABC$ is tiled
by a right-angled tile $(\alpha,\beta,\frac \pi 2)$, then the 
base angles of $ABC$ are either equal to $\beta$ or to $\alpha$.
That follows, because in Table~\ref{table:isosceles}, there is 
only one entry corresponding to a right-angled tile, namely
the second line.  Readers are invited to try to prove that 
directly, without appeal to Laczkovich,  in order to gain a deeper
appreciation for Laczkovich's work.

\section{Some Number-Theoretic Facts}
The facts in this section may not be 
well-known to all our readers, and their proofs are short. 

The following lemma can be found in any number theory textbook,
for example in \cite{hardy-wright}[Theorem~366, page~299].

\begin{lemma}\label{lemma:sumsofsquares} 
An integer $N$ can be written as a sum of two integer squares if and only if the 
squarefree part of $N$ is not divisible by any prime of the form $4n+3$.
\end{lemma}

\begin{lemma}\label{lemma:doublesquares}
An integer $N$ is a sum of two squares if and only if 
$2N$ is a sum of two squares.
\end{lemma}

\begin{proof} The lemma follows immediately from the identities
\begin{align*}
(p-q)^2 + (p+q)^2 &= 2(p^2+q^2),\\
\intertext{and}
\left( \frac{p-q} 2 \right)^2 + \left(\frac{p+q} 2\right)^2 &= \frac 1 2 \,(p^2+q^2).
\end{align*}
 This lemma is also a corollary of
Lemma~\ref{lemma:sumsofsquares}, of course, but that is not needed.
\end{proof}

The following lemma identifies those relatively few rational multiples of $\pi$ that have rational tangents or 
whose sine and cosine satisfy a polynomial of low degree over $\Q$.

\begin{lemma}  \label{lemma:euler}
Let $\theta = 2m \pi/n$,  where $m$ and $n$ have no common factor.
Suppose $\cos \theta$ is algebraic of degree 1 or 2 over $\Q$.
Then $n$ is one of $1,2,3,4,5,6,8,10,12$.   If both $\cos \theta$ and
$\sin \theta$ have degree 1 or 2 over $\Q$, then $n$ is one of
$1,2,3,4,6,8,12$.
\end{lemma}

\begin{proof}  Let $\varphi$ be the Euler totient function.
If $n \le 2$, then $\cos\theta = \pm 1$ has degree 1.  If $n \ge 3$,
then by \cite{niven}, Theorem~3.9, page~37, $\cos\theta$ has degree
$\varphi(n)/2$; so $\cos\theta$ has degree 1 or 2 if and only if
$\varphi(n) = 2$ or $4$, that is, $n$ is one of $3,4,6$ or $5,8,10,12$.
That proves the first claim.

For the second claim, it suffices to show that $\sin\theta$ has degree
greater than 2 when $n = 5$ or $n = 10$.  In those cases $\theta$ is,
up to sign and up to adding a multiple of $\pi$, one of $\pi/5$ and
$2\pi/5$, and
$$ \sin^2 \frac \pi 5 = \frac{5-\sqrt 5}{8}, \qquad
   \sin^2 \frac{2\pi} 5 = \frac{5+\sqrt 5}{8}. $$
These are the roots of $16y^2 - 20y + 5$, so $\sin\theta$ is a root of
$$ 16x^4 - 20x^2 + 5, $$
which is irreducible over $\Q$ by Eisenstein's criterion with the
prime $5$.  Hence $\sin\theta$ has degree 4.
\end{proof}

\begin{lemma} \label{lemma:rationalsquares}
If the integer $n$ is a sum of two rational squares then 
it is a sum of two integer squares.
\end{lemma}

\begin{proof} Suppose $n = (p/q)^2 + (s/t)^2$. 
Then $(qt)^2 n = (pt)^2 + (qs)^2$.  Then by Lemma~\ref{lemma:sumsofsquares},
the square-free part of $n$ is not divisible  by any prime  $p$ with $p \equiv 3 \pmod{4}$
Then by a second application of Lemma~\ref{lemma:sumsofsquares},
$n$ is a sum of two integer squares. 
\end{proof}

\section{Tilings of an Isosceles $ABC$ by a Right-angled Tile: Examples}
\begin{figure}[ht]
\begin{center}
\includegraphics[width=\textwidth]{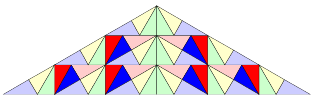}
\end{center}
\caption{A 54-tiling; $N/2$ is three times a square. Tile is 30-60-90.}
\label{figure:54}
\end{figure}

\begin{figure}[ht]
\begin{center}
\includegraphics[width=0.45\textwidth]{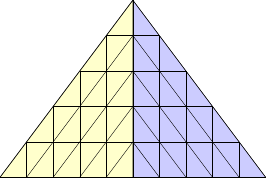}
\includegraphics[width=0.45\textwidth]{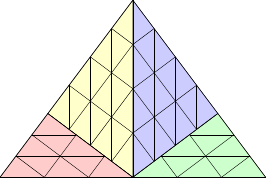}
\end{center}
\caption{$N$ is twice a square or  twice a sum of squares. }
\label{figure:severalisosceles}
\end{figure}

\begin{figure}[ht]
\begin{center}
\includegraphics[width=0.45\textwidth]{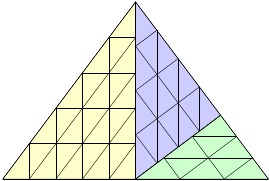}
\end{center}
\caption{50 is both twice a square and twice a sum of squares. }
\label{figure:severalisosceles2}
\end{figure}

\begin{figure}[ht]
\begin{center}
\includegraphics{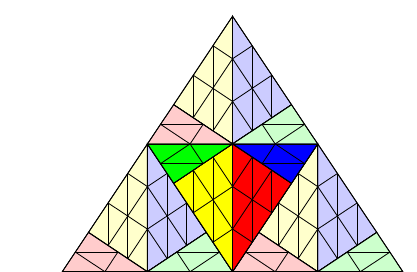}
\end{center}
\caption{$N=104$, eight essential segments, base angles about $56^\circ$.}
\label{figure:8essentials}
\end{figure}

\begin{figure}[ht]
\begin{center}
\includegraphics[width=0.45\textwidth]{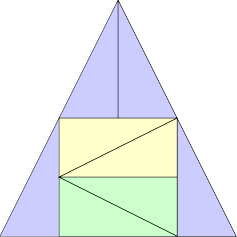}
\end{center}
\caption{The altitude need not be part of the tiling.}
\label{figure:nicelytiled}
\end{figure}

Is it possible
to have more complicated tilings without essential segments? Yes, 
because when two tiles share their hypotenuses, they form a rectangle, 
and we can just draw the diagonal of that rectangle the other way.
In this way we can produce (exponentially) 
many different tilings, but they differ only in 
this trivial way.
And sometimes, as shown in  Figure~\ref{figure:nicelytiled}, 
even those rectangles can be rotated.  That figure
also shows that a tiling need not necessarily include
the altitude of $ABC$.

\begin{figure}[ht]
\includegraphics{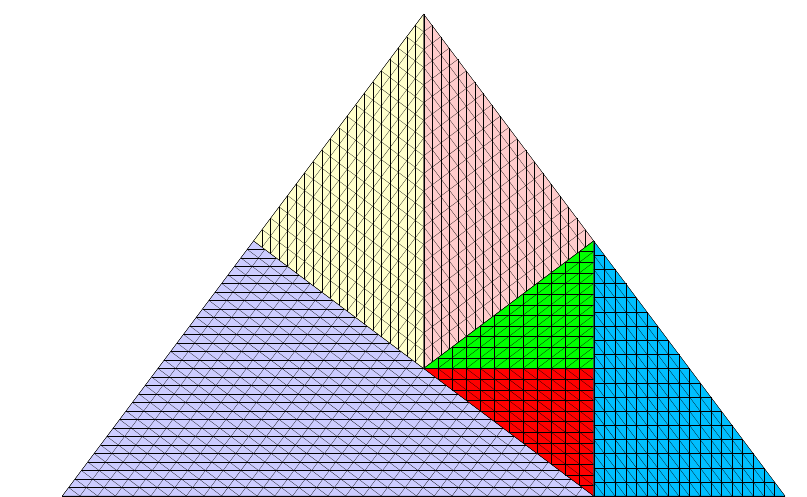}
\caption{$N = 2312$, $N/2 = 34^2$, $(a,b,c) = (3,4,5)$.}
\label{figure:2312}
\end{figure}

In the tilings
based on two biquadratic tilings, there are no $c$
edges on $AB$ and $BC$, while in the tilings based 
on two quadratic tilings, there are only $c$ edges.
There are of course some hybrid tilings when a square
is also a sum of squares, in which
$AB$ falls under one case and $BC$ under the other.
If $N/2$
is not a square (as is the case for the biquadratic tilings)
then there are no $c$ edges on $AB$ and $BC$, as we see 
in the biquadratic tilings (and prove in the next section).  
 
All these tilings, in which $N/2$ is a sum of squares, 
involve essential segments (where tiles
of different lengths occur on the two sides of an internal line).   
One sees such linear relations in two of the tilings
illustrated in Figure~\ref{figure:severalisosceles}.

\section{Laczkovich's Graphs $\Gamma_a$}
In trying to prove the impossibility of certain tilings 
directly, it is easy to become involved in complicated 
arguments with many cases, involving complicated diagrams.
Laczkovich had the brilliant idea to abstract some of these
arguments using graph theory.  The definition will not be 
grasped immediately, but instead will require time and the 
study of examples to understand. However, it leads to very 
elegant proofs of theorems that are much more complicated
or impossible to prove more directly.  To emphasize its
importance, we devote a whole section just to the definition.

 Given a tiling of a triangle $ABC$,
 an {\em internal segment} 
 is a line segment connecting two vertices
of the tiling that is contained in the union of the boundaries of the 
tiles, and lies in the interior of $ABC$ except possibly for its endpoints.
A {\em maximal segment} is an internal segment that is not part of a longer
internal segment.   A segment is {\em terminated}  at a vertex $P$ if it has
tiles on both sides with vertices at $P$.  In that case there may or may not 
be a continuation of that segment past $P$.  
A {\em left-terminated segment} is an internal segment 
$XY$ that is terminated at $X$.  
A {\em left-maximal segment} is an internal segment $PQ$ that 
cannot be extended past $P$ to a longer internal segment $UPQ$.
In these two concepts, we are using directed segments; so $PQ$
is not the same as $QP$ in this context. The ``left'' in ''left maximal''
refers to the fact that $P$ is listed to the left of $Q$ in $PQ$.
A tile is {\em supported by}
$XY$ if one edge of the tile lies on $XY$.   The internal segment $XY$ is 
said to have {\em all $c$'s on the left} if the endpoints $X$ and $Y$ are
vertices of tiles supported by $XY$ and lying on the left side of $XY$,
and all tiles supported by $XY$ lying on the left of $XY$ have their 
$c$ edges on $XY$.  The phrase {\em all $c$'s on the right} is defined similarly.  
Here again $XY$ is a directed segment, so the concept ``left side'' of $XY$
makes sense; but this is a different sense of the English word ``left''
than in ``left-terminated.''

An internal segment $XY$ is said to {\em witness the relation} 
$jc = \ell a + mb$ if one side has $j$ more $c$ edges than the other,
and the other has $\ell$ more $a$ edges and $m$ more $b$ edges than
the first.   The simplest example is when
 $XY$ has all $c$'s on one side, and exactly $j$ of them (that 
is, the length of $XY$ is $jc$), and on the other side $XY$ supports
$\ell$ tiles with their $a$ edges on $XY$ and $m$ tiles with their
$b$ edges on $XY$ (in any order) and no other tiles, and the endpoints
$X$ and $Y$
are vertices of tiles on both sides of $XY$.    
Similarly we use the terminology ``$XY$ witnesses a relation 
$jb = \ell a + mc$.''

An internal segment that witnesses a relation is called 
an {\em essential segment}.  The definition allows that
an essential segment might have different numbers of 
tiles of lengths $a,b,c$ on its two sides, without necessarily having
all the tiles on one side be the same length, but often it 
is the case that all the tiles on one side are the same length.   

To be sure that you understand the concept of ``essential segment'',  
identify the eight essential 
segments in Figure~\ref{figure:8essentials}.  Also identify in 
that figure some internal segments $PQ$ that are not essential segments,
because each side of $PQ$ supports tiles with three $b$ and two $a$ sides 
on $PQ$.  Those $PQ$ connect the midpoints of the sides of $ABC$ in 
that figure.

The following definition is equivalent to the one given 
in \cite[page~346]{laczkovich2012}, except that there condition (iv) is
automatic because of an additional assumption about the tile. 

\begin{definition}
 [The directed graph $\Gamma_a$]
  \label{definition:Gamma} Given a tiling
of some triangle, 
the {\em nodes} of the graph $\Gamma_a$ are the endpoints of links as defined below. 
A {\em link} of $\Gamma_a$ connects
vertices $X$ and $Y$ if each (all) of the following holds:

(i) the segment $XY$ is a left-maximal
internal segment having all $a$ edges on one side (say Side~1) of $XY$, and

(ii) on the other side of $XY$ (say Side~2) the first tile (the one with a 
vertex at $X$) does 
not have its $a$ edge on $XY$, and

(iii) at vertex $Y$, there is
another tile supported by $XY$ extended past $Y$ on Side~1 of $XY$ with a vertex at $Y$,
 that does not have its
$a$ edge on Side~1, and 
\smallskip
 
(iv)  no tile supported by $XY$ on Side~2 of $XY$ has a vertex at $Y$.
\smallskip

   The directed graphs $\Gamma_b$ and $\Gamma_c$
are defined
similarly. 
\end{definition}

\begin{remark}  We use ``link'' instead of ``edge'' with 
these graphs, to avoid confusion with tile ``edges.''
If $XY$ is a link in $\Gamma_a$, then $XY$ does not terminate at $Y$,
because there is another tile past $Y$ whose $a$ side is not on $XY$;
and also because $Y$ lies on the interior of an edge of a tile on 
the other side of $XY$.  In graph theory, sometimes the word ``vertex''
is used.  We use ``node'' here because ``vertex'' is already in use for 
triangle tilings.  Every node of these graphs is a vertex of the tiling,
but not every vertex is a node.  
\end{remark}

The reader is recommended to identify all three graphs $\Gamma_a$, $\Gamma_b$,
and $\Gamma_c$  in 
Figure~\ref{figure:8essentials}.  Hint:  $\Gamma_c$ is empty;
$\Gamma_b$ has four links, 
starting at the midpoints of the sides of $ABC$;
 $\Gamma_a$ is, in this example, the same graph as $\Gamma_b$.

\section{Tilings of an Isosceles $ABC$ by a Right-angled Tile: Theory}
Laczkovich studied the possible shapes of tiles that can tile 
an isosceles triangle, but did not characterize the possible $N$.
We do so in this section for right-angled tiles and $N$ even. 
We have two ways to tile an isosceles triangle by a right 
triangle:  either tile each of its two halves by a 
quadratic tiling,  in which case $N$ is twice a square, 
or tile each of its halves with a biquadratic tiling,
in which case $N$ is twice a sum of squares.  See
Figures~\ref{figure:severalisosceles} and \ref{figure:severalisosceles2}.
The main theorem in this section shows that these are
the only possible values of $N$.  However, Figure~\ref{figure:2312} shows
that, when $N/2$ is a square, there are also more complicated tilings.

\begin{lemma}\label{lemma:Gamma}
 Suppose isosceles (or equilateral) triangle $ABC$ with base angles 
$\beta$ is $N$-tiled by tile $(\alpha,\beta,\pi/2)$ with sides
$(a,b,c)$, and $\alpha$ is not a rational multiple of $\pi$,
or $\alpha$ is an odd multiple of $\beta$.
Let $PQ$ be a link in $\Gamma_c$.  Then there 
are two adjacent tiles with vertices at $Q$ whose common boundary
contains
an $a$ or $b$ edge of one tile, and a $c$ edge of the other tile.
\end{lemma}

\begin{remark} Consider Figure~\ref{figure:54}, in which 
$\beta = \pi/6$ and $\alpha = \pi/3$, so $\alpha$ is not an 
odd multiple of $\beta$.  Observe that the present lemma fails 
in the tiling of Figure~\ref{figure:54}, showing that the hypothesis
that $\alpha$ is an odd multiple of $\beta$ cannot be removed. 
\end{remark}

\begin{proof}[Proof of Lemma~\ref{lemma:Gamma}] Let $\Delta_1,\ldots,\Delta_k$ be 
tiles with vertices at $Q$, numbered so that $\Delta_1$ is 
supported by $PQ$ (and hence has side $c$ on $PQ$),
$\Delta_i$ and $\Delta_{i+1}$ are adjacent, and $\Delta_k$ 
has one edge extending $PQ$ past $Q$.   Since $PQ$ is a link
in $\Gamma_c$, there does exist such a tile $\Delta_k$, and 
$\Delta_k$ has its $a$ or $b$ edge extending $PQ$.
There may or may not be tiles on the other side of $PQ$
with a vertex at $Q$, but if so, we do not list them 
among the $\Delta_i$.

If there are an even number of tiles with vertices at $Q$ and 
an $\alpha$ or $\beta$ angle at $Q$, then each has a $c$ edge
ending at $Q$.  Since $\Delta_1$ has its $c$ edge ending at $Q$
and $\Delta_k$ does not, the remaining odd number of $c$ edges
cannot all be paired with other $c$ edges supported by the same line.
Therefore, there is an $a/c$ or $a/b$ edge, as claimed.
Therefore, we may assume the number of such tiles is odd. 
At most one tile can have its right angle at $Q$, and it cannot be 
$\Delta_1$.

Now, if $\alpha$ is not a rational multiple of $\pi$,  then 
either $k=2$ and both angles
are right angles, or $k=4$ and there are two $\alpha$ and two
$\beta$ angles, or $k=3$ and there are one each of $(\alpha,\beta,\frac \pi 2)$.
In all these cases, the above condition that there are an even 
number of tiles at $Q$ with an $\alpha$ or $\beta$ angle at $Q$
is fulfilled.   

If $\alpha$ is an odd multiple $m$ of $\beta$, then the same argument
works, as the number of tiles with vertices at $Q$ that do not have
their right angle at $Q$ will still be even.
\end{proof}

\begin{lemma}[Propagation]\label{lemma:propagation}
 Suppose isosceles (or equilateral) triangle $ABC$ with base angles 
$\beta$ is $N$-tiled by tile $(\alpha,\beta,\pi/2)$ with sides
$(a,b,c)$, and $\alpha$ is not a rational multiple of $\pi$,
or $\alpha$ is an odd multiple of $\beta$.
 Suppose
there is no relation $jc = ua + vb$ with $j > 0$ and $u,v \ge 0$, 
and $u,v,j$ integers.  Then
\smallskip

(i) Let $PQ$ be a link in $\Gamma_c$.  Then there is a vertex $R$
such that $QR$ is a link in $\Gamma_c$.
\smallskip

(ii) The in-degree and out-degree of each node of $\Gamma_c$ is 
exactly 1.
\end{lemma}

\begin{remark} Note that in Figure~\ref{figure:8essentials}, 
it is not true that the in-degree of every vertex in $\Gamma_b$ is equal to the
out-degree of that vertex.  That is because, in that tiling,
there is a relation $2b = 3a$.  Please take the time to 
verify that in Figure~\ref{figure:54},
the conclusion of this lemma fails,  but so does the hypothesis
that $\alpha$ is an odd multiple of $\beta$; this shows the necessity
of that hypothesis. 
\end{remark}

\begin{proof}[Proof of Lemma~\ref{lemma:propagation}]  According to Lemma~\ref{lemma:Gamma},
there is an outgoing $a/c$ edge from $Q$.  Let $R$ lie on 
the internal segment containing that edge, as far as possible
from $Q$ such that there are only $c$ edges on one side 
of $QR$.  Then $R$ is not a vertex of a tile on the other side of 
$QR$, since that would give rise to relation $jc = ua + vb$,
where $j$ is the excess of the number of $c$ edges on one side
of $QR$ over the other.  Then, by definition of $\Gamma_c$,
$QR$ is a link in $\Gamma_c$.  That completes the proof of (i).

Next we observe that the in-degree of each node $Q$ of $\Gamma_c$
is at most 1.  For, if $PQ$ is a link of $\Gamma_c$, then 
$PQ$ is part of an internal segment of the tiling that extends
past $Q$ and on one side of $PQ$, $Q$ is not a vertex of any tile on that side.
Hence, no other internal segment can pass through $Q$; hence there
is no other link of $\Gamma_c$ ending at $Q$.

By part~(i), the out-degree of each node of $\Gamma_q$ that
has positive in-degree is at least 1.  Hence, the out-degree always
is greater than or equal to the in-degree.  Since every link
has one head and one tail, the total in-degree is equal to the total
out-degree.  Therefore, the in-degree and out-degree are equal
at every node.  Therefore, if the in-degree is positive, both 
the in-degree and out-degree are 1.  If the in-degree is zero,
so must the out-degree be zero, but then that node is not in 
$\Gamma_c$ at all.  
\end{proof}

\begin{corollary}  \label{lemma:allc}
Suppose isosceles (or equilateral) triangle $ABC$ with base $AC$ and base angles 
$\beta$ is $N$-tiled by tile $(\alpha,\beta,\pi/2)$ with sides
$(a,b,c)$, and $\alpha$ is not a rational multiple of $\pi$,
or $\alpha$ is an odd multiple of $\beta$.   Suppose
there is no relation $jc = ua + vb$ with $j > 0$ and $u,v \ge 0$, 
and $u,v,j$ integers.  Then $AC$ is composed only of $c$ edges
and there are no $c$ edges on $AB$ or $BC$,  or $AB$ and $BC$
are composed only of $c$ edges and there is no $c$ edge on $AC$.
\end{corollary}

\begin{proof}      According to Lemma~\ref{lemma:propagation},
each link  $PQ$ in $\Gamma_c$ has a corresponding link $QR$. 

Suppose there is a tile with a $c$ edge on $AB$. 
Unless the entire segment $AB$ supports only 
tiles with their $c$ edges on $AB$, there will be a segment 
$PQ$ lying on $AB$ such that $PQ$ is composed of $c$ edges
but beyond $Q$ there is another tile with an $a$ or $b$ edge
on $PQ$.  That contradicts Lemma~\ref{lemma:propagation},
as there will be an outgoing link of $\Gamma_c$ from $Q$, but
at a boundary point, there can be no incoming link. 
Hence, if there is any $c$ edge on $AB$, then $AB$ is composed
entirely of $c$ edges.  The same argument applies to $BC$ and $AC$.

At vertex $A$, there is an angle $\beta$.  There must be a 
single tile there, with its $\beta$ angle at $A$, since 
either $\beta < \alpha$ or $\alpha$ is not a rational multiple of $\pi$.
The $c$ edge of that tile must lie on $AB$ or $AC$.  If it lies
on $AB$, then $AB$ is composed entirely of $c$ edges.  If it 
lies on $AC$, then $AC$ is composed entirely of $c$ edges.  If $AC$
is composed of $c$ edges, then $AB$ and $BC$ do not contain any 
$c$ edges, since if they contained one, then there would also be 
a $c$ edge at $A$ or $C$, which is impossible since the tiles
at $A$ and $C$ have only one $c$ edge, and it is on $AC$. 
\end{proof}

\begin{lemma} \label{lemma:isosceleshelper2}
Suppose isosceles (or equilateral) triangle $ABC$ with base angles
$\beta$ at $A$ and $C$ is $N$-tiled by a tile with angles
 $(\alpha,\beta,\pi/2)$ and sides $(a,b,1)$.   Suppose 
 $\beta \neq \pi/6$ and 
$\sqrt{N/2}$ is 
 irrational (i.e., $N$ is not twice a square).
  Then $a$ and $b$ belong
 to $\Q(\sqrt{N/2})$.
 \end{lemma}
 
\begin{remark} The tiling in Figure~\ref{figure:54}
shows that the exception for $\beta = \pi/6$ is necessary.
\end{remark}
 
\begin{proof}[Proof of Lemma~\ref{lemma:isosceleshelper2}]
 Let $X = \vert AB\vert$.  
Twice the area of $ABC$ is the cross product of the two equal sides,
which is 
\begin{eqnarray*}
X^2 \sin 2\alpha &=& 2 X^2 \sin \alpha \cos \alpha \\
&=& 2X^2 \sin \alpha \sin \beta \ = \ 2X^2 ab.
\end{eqnarray*}
Twice the area of the tile is $ab$.  Since $N$ tiles cover $ABC$
we have the area equation
$$ 2X^2 = N. $$
Define
$$ \lambda := \sqrt{N/2}. $$
Then $X = \lambda$.  Let $M$ be the midpoint of 
the base $AC$.  Then triangle $ABM$ has a right angle 
at $M$, angle $\alpha$ at $B$, and angle $\beta$ at $A$,
so it is similar to the tile.  Therefore, $AM = X \sin \alpha =
\lambda a$.  Therefore, $\vert AC \vert = 2a \lambda$.
Since there is a tiling of $ABC$, 
there
are non-negative integers $(p,q,r)$  and $(s,t,u)$ 
such that 
$$ \lambda = pa + qb + r \qquad\mbox{and}\qquad 2a\lambda = sa + tb + u. $$

We write this as a system of equations in unknowns $a,b$:
\begin{align*}
pa + qb &= \lambda -r,\\
\intertext{and}
(s-2\lambda)a + tb &= -u.
\end{align*}
The determinant $D = pt-q(s-2\lambda)$.
If $D \neq 0$, then both $a$ and $b$ belong to $\Q(\lambda)$, since
that field contains $D$ and the right-hand sides of the system. 

Now suppose, for proof by contradiction, that  $D=0$.  Then since $\lambda$ is irrational, we have 
$q = 0$ and $pt = qs = 0$. Since $pa + qb = \lambda -r \neq 0$, we have
$p \neq 0$ and hence $t = 0$.  Then the two equations become
$$ pa = \lambda-r \qquad\mbox{and}\qquad (s-2\lambda)a = -u. $$
Multiplying these two equations we have
$$ -2\lambda^2 + \lambda(s+2r) - sr = -pu. $$
Since $2\lambda^2 = N$, we have
$$ -N + \lambda(s+2r) - sr = -pu. $$
  Since $\lambda$ is irrational, and $s$ and $r$ are non-negative,
we have $s=r=0$.  Hence, $AC = u$ is made up only of $c$ edges (each of length 1)
and $X = pa + qb$ , so $AB$ has no $c$ edges.  Since also $q=0$
(as shown above), $X = pa$,  so $AB$ also has no $b$ edges, and is 
therefore composed entirely of $a$ edges.  
A similar argument applies
to show that $CB$ also is composed entirely of $a$ edges. 

Therefore, $a = X/p = \lambda / p$ belongs to $\Q(\lambda)$.  We cannot,
however, immediately conclude that $b$ is in $\Q(\lambda)$. 
We can, however, conclude that there is no relation 
$jc = ua + vb$ with integers $j, u, v$:   Suppose there
was such a relation.
  Since $c = 1$, that would mean $b = (j-ua)/v \in \Q(\lambda)$.
  Therefore, as claimed,
there is no relation $jc = ua + vb$.

Now consider the tiles with a vertex at $B$, where triangle $ABC$
has an angle of $2 \alpha$.  Since sides $AB$ and $CB$ are
composed entirely of $a$ edges, the tiles at $B$ supported
by $AB$ and $CB$  have their $a$ edges on $AB$ or $BC$,
and hence do not have their $\alpha$ angles at $B$. 
In particular, the case when there is only one tile at $B$
is ruled out, since it cannot have an $a$ edge on both $AB$
and $BC$; and the case of just two tiles at $B$, both with 
$\alpha$ angles at $B$, is also ruled out,  since neither would 
have an $a$ edge on $AB$ or $BC$.  
Therefore
there must be some tiles with $\beta$ angles at $B$. Either 
$\alpha$, or $2\alpha$, or $2\alpha-\pi/2$ must be a multiple
of $\beta$.  In the latter case, $2 \alpha - \pi/2 = \alpha -\beta$
is a multiple of $\beta$, so $\alpha$ is a multiple of $\beta$.
In all of these cases, then, $2\alpha$ is a multiple of $\beta$,
say $2\alpha = m \beta$, with $m > 1$.  Then 
$\beta = \frac \pi 2 - \alpha$, so
$$ 2\beta = \pi - 2 \alpha = \pi - m \beta. $$
Hence $(m+2)\beta = \pi$, and
$$ \beta = \frac {\pi} {m+2} = \frac {2\pi} {2(m+2)}. $$
Now $\cos \beta = \sin \alpha = a$ belongs to $\Q(\lambda)$.
By Lemma~\ref{lemma:euler},  
$2(m+2)$ is equal to one of $1,2,3,4,5,6,8,10,12$.
Therefore, $m+2$ is one of $1,2,3,4,5,6$.  Since $m > 1$,
we have $m \in \{2,3,4\}$.

\case{1} $m=2$. Then $\alpha = \beta = \pi/4$, so $a=b$. 
Since $a \in \Q(\lambda)$,  we have $b \in \Q(\lambda)$,
completing this case.

\case{2} $m = 3$.  Then $2\alpha = 3\beta$, so $\alpha = 3\pi/10$ and
$\beta = \pi/5$.  We have shown above that $AB$ is composed entirely of
$a$ edges, so $\lambda = X = pa$, and hence
$$ a^2 = \frac{\lambda^2}{p^2} = \frac{N}{2p^2}, $$
which is rational.  On the other hand, with $c = 1$ we have
$a = \sin \alpha = \cos \beta = \cos(\pi/5) = (1+\sqrt 5)/4$, so
$$ a^2 = \frac{3 + \sqrt 5}{8}, $$
which is irrational.  This contradiction shows that Case~2 is impossible.

\case{3} $m = 4$.  Then $2\alpha = 4 \beta$,
$\alpha = \pi/3$ and $\beta = \pi/6$,  which is ruled
out by hypothesis (and has to be, because of 
Figure~\ref{figure:54}).  Note that Lemma~\ref{lemma:propagation} does not 
apply, since $2\alpha$ is an even multiple of $\beta$. 
\end{proof}

\begin{lemma} \label{lemma:isosceleshelper3}
 Suppose isosceles triangle  $ABC$ with base angles
$\beta$ at $A$ and $C$ is $N$-tiled by a tile with angles
 $(\alpha,\beta,\pi/2)$ and sides $(a,b,c)$.  
Suppose $\beta \not\in \{\tfrac \pi 6, \tfrac \pi 4, \tfrac \pi 3\}$
 and $\sqrt{N/2}$ is irrational.
  Then
 \smallskip
 
 (i) $a/c$ and $b/c$ are 
 rational multiples of $\lambda = \sqrt{N/2}$, and 
 \smallskip
 
 (ii) $AC$ is composed only of $c$ edges, and there are no 
 $c$ edges on $AB$ or $BC$.
 \end{lemma}
 
\begin{proof} Without loss of generality,
we may assume $c=1$.   Since by hypothesis, $\sqrt{N/2}$
is irrational and $\beta \neq \frac \pi 6$, we 
can apply  Lemma~\ref{lemma:isosceleshelper2} to conclude
that  $a$ and $b$
belong to $\Q(\lambda)$.

We will show that
$$\lambda = X = \vert AB \vert.$$
Twice the area of $ABC$ is equal on the one hand to $Nab$,
and on the other to 
$$X^2 \sin 2\alpha =X^2 2 \sin \alpha \cos \alpha = 2X^2 ab.$$
Hence, $N = 2X^2$.  However, $N = 2\lambda^2$ by definition of $\lambda$.
Hence, $X = \lambda$, as claimed.

Since $b/c = \sin\beta$ and $a/c = \cos\beta$ are both rational
multiples of $\lambda$,  they are irrational, but 
 $b/a= \tan \beta$ is  rational.  By Lemma~\ref{lemma:euler},
if $\beta$ is a rational multiple of $\pi$, it is $2\pi m/n$,
where $n = 2,3,4,6,8,$ or $12$.  Since there are two 
base angles of $\beta$, we have $\beta < \pi/2$; and 
by hypothesis $\beta \neq \pi/6$ and $\beta \neq \pi/4$ and $\beta \neq \pi/3$.
Therefore $\alpha$ and $\beta$ are not 
rational multiples of $\pi$. 

The vertex
angle $B$ is equal to $2\alpha$, since the base angles are $\beta$ and $\alpha+\beta = \tfrac \pi 2$.
We will show that the vertex 
angle is filled by two tiles with their $\alpha$ angles at $B$: Since  $\alpha$ and $\beta$
are not rational multiples of $\pi$, if $p\alpha + q\beta = 2\alpha$, and $q \neq 0$,
then $\beta = \frac{2-p} q \alpha$, which is a contradiction.  Therefore, as claimed,
the vertex angle is filled by two tiles with their $\alpha$ angles at $B$.  

 Let $a = x\lambda + y$ and $b= z\lambda + w$,
where $x,y,z,w$ are rational.  Then for some nonnegative integers
$p,q,r$, 
$$ \lambda = pa + qb + r = (px + qz)\lambda + (py + qw + r).$$
Since $\lambda$ is irrational, 
\begin{equation}
py + qw + r = 0. \label{eq:8067}
\end{equation}  
We also have 
$$ 1 = a^2 + b^2 = (x\lambda + y)^2 + (z\lambda + w)^2 = (x^2 + z^2) \frac N 2 + y^2+w^2 + 2(xy + zw)\lambda.$$
Since $\lambda$ is irrational we have 
$$ xy + zw = 0. $$

\case{1}  $xy \neq 0$.  Then $xy$ and $zw$ have different signs. 
 We will show that $x > 0$ and $z > 0$.   We will prove this 
by cases, according to the sign of $xy$.  First suppose $xy > 0$. Then
$x$ and $y$ are of the same sign.
Since $a = \lambda x + y > 0$, it follows that $x,y > 0$.  
The tile at $B$ supported by $AB$ has its 
$\alpha$ angle at $B$, as shown above.  Hence,
either its $b$ or $c$ edge is on $AB$.  That is, not both $q$ and $r$ are zero. 
Hence Equation~(\ref{eq:8067}) 
implies that $qw \le 0$. Hence $w \le 0$. Since $zw < 0$,  we have $z > 0$ and $w < 0$.
Thus, the claim $x > 0$ and $z > 0$ holds if $xy >0$.

Now suppose $xy < 0$.  Then $zw > 0$, and we find $z > 0$ and $w > 0$, 
since $b = z\lambda + w > 0$.    Since not both $q$ and $r$ are zero, 
and $w > 0$,  we have $qw + r > 0$.  Then  Equation~(\ref{eq:8067}) implies $py < 0$;
hence $y < 0$.  Since $xy < 0$, we conclude $x > 0$,  establishing the 
claim $x > 0$ and $z > 0$ also in case $xy < 0$.  Thus, we have 
proved that $xy \neq 0$ implies $x > 0$ and $z > 0$.   

We will show that there is no relation $jc = j = ua + vb$ with non-negative integers $j,u,v$ and
$j > 0$.  Suppose such a relation exists; then
\begin{eqnarray*}
j &=& ua + vb \\
&=& u(x\lambda + y) + v(z \lambda + w) \\
&=& (ux + vz) \lambda + (uy + vw).
\end{eqnarray*}
Since $\lambda$ is irrational, we have $ux + vz = 0$.   Since $x$ and $z$ are positive
and $u,v \ge 0$,  this implies $u = v = 0$,   which is a contradiction, since $j > 0$.

Then by Corollary~\ref{lemma:allc},  $AC$ is composed of all $c$ edges, and 
there are no $c$ edges on $AB$.  Let $s$ be the number of $c$ edges on $AC$. 
 Since the length of $AC$ is $2a \lambda$, 
we conclude that $a = s/(2\lambda) = (s/N) \lambda$, since $2\lambda^2 = N$.
Hence, $a$ is a rational multiple of $\lambda$.   Since $AB$ has no $c$ edges,
$r = 0$ and $X = \lambda = pa + qb$.  Hence, $b =( \lambda - pa)/q $ is also 
a rational multiple of $\lambda$.  That completes the proof in Case~1.

\case{2} $xy = 0$.  Then also $zw = 0$.  We will show this case implies $y = w = 0$.
  We argue by cases on whether $x = 0$ or not.
 If $x = 0$, then $a = y$ is rational and $y > 0$.
If $z = 0$, then $b = w > 0$,  contradicting Equation~(\ref{eq:8067}).  If $w = 0$, then $b = z\lambda$,
and
$$ 2 a \lambda = sa + tb + u = (xs + tz)\lambda + (sy + u). $$
Hence, since $\lambda$ is irrational, $sy + u = 0$; since $s, u \ge 0$
and $y > 0$, we have $s = u = 0$,
and $AC$ is composed only of $b$ edges, which is impossible, since the tile at $A$
cannot have its $b$ edge on $AC$.

On the other hand, if $x \neq 0$ then since $xy = 0$ we have $y = 0$.  Then $qw + r = 0$
by Equation~(\ref{eq:8067}).   If $z = 0$, then $b = w > 0$ and therefore $q=r=0$, which is a 
contradiction.    The only remaining possibility is $y=w=0$, as claimed. 
Then $a = x\lambda$ and 
$b = z\lambda$.   Then $X = \lambda = pa + qb + r = (x+z)\lambda + r$,
so $r = 0$ and 
\begin{eqnarray*}
Y &=& 2a\lambda  \\
&=& 2x\lambda^2  \\
&=& xN  \\
&=& sa + tb + u \\
&=& (sx + tz)\lambda + u,
\end{eqnarray*}
 which implies $s=t=0$.
That completes Case~2.
\end{proof}

 \begin{theorem} \label{theorem:rationalisosceles}
Suppose $ABC$ is isosceles with base angles $\beta$,
or $ABC$ is equilateral, 
 and $ABC$ is tiled by triangle $T$ similar to half of $ABC$. If 
$\alpha$ is a rational multiple of $\pi$, then either
\smallskip

(i) $N$ is even and $N/2$ is a square, or
\smallskip

(ii)  $N$ is a square and $\beta = \pi/4$ and $ABC$ has base
angles $\pi/4$, or 
\smallskip

(iii) $N/2$ is three times a square, and the tile is a 30-60-90 triangle:
either $\beta = \pi/6$ and $\alpha = \pi/3$, so that $ABC$ has angles
$(\frac \pi 6, \frac \pi 6, \frac {2\pi} 3)$, or $\beta = \pi/3$ and
$\alpha = \pi/6$, so that $ABC$ is equilateral.
\end{theorem}

\begin{remark}  One possible tiling under case (iii) of the 
theorem is  illustrated in 
Figure~\ref{figure:54}.  $N$ can be odd in case (ii), since half 
the triangle $ABC$ is similar to $ABC$, so quadratic tilings are allowed.
\end{remark}

\begin{proof} 
We begin by remarking that $N/2$ is a rational square if and only 
if it is an integer square, since if it is a rational square then 
$2N$ is an integer that is a rational square, hence it is an integer 
square. Then it is the square of an even integer $2m$, so $N/2 = m^2$.
Hence, if $\sqrt{N/2}$ is rational, then $N$ is even and $N/2$ is a
square.  Then condition (i) holds. 

Therefore, we may assume that 
$\sqrt{N/2}$ is irrational.  We now divide into cases according 
as $\beta = \frac \pi 6$ or not.
\smallskip

First assume $\beta \neq \frac \pi 6$.  
Then by Lemma \ref{lemma:isosceleshelper2}, $a = \sin \alpha$ and $b = \cos \alpha$ 
belong to $\Q(\sqrt{N/2})$. 
 By hypothesis, $\alpha$ is a rational multiple of $\pi$.
These two facts make Lemma \ref{lemma:euler} applicable, so 
we can drastically limit the possible values of $\alpha$.
Namely, by Lemma \ref{lemma:euler}, $\alpha$ and $\beta$ are odd
multiples of $2\pi/n$, where $n$ is one of $6,8,12$; that is,
they are odd multiples of $\frac \pi 3, \frac \pi 4, or \frac \pi 6$.
Since they are both less than $\frac \pi 2$, $\alpha$ and $\beta$
must be exactly $\frac \pi 3$, $\frac \pi 4$, or $\frac \pi 6$. 
Those are the values of $\alpha$ and $\beta$ allowed in the statement 
of the lemma.  We arrived at that conclusion under the assumption
that $\beta \neq \frac \pi 6$, but since that conclusion includes
$\beta = \frac \pi 6$, it holds without that assumption. That is,
we have proved outright that $\alpha$ and $\beta$ are equal to 
$\frac \pi 3$, $\frac \pi 4$, or $\frac \pi 6$.

It remains to show that $N$ has one of the stated values.
Let $X = pa + rb + q$ be the length of $AB$.
Then the area equation is 
$$  Nab = X^2 \sin 2 \alpha = 2X^2 a b,$$
so $N = 2X^2$.

\case{1} $\alpha = \beta = \frac \pi 4$.  Then with $c=1$ we have
$a = b = 1/\sqrt 2$, so the area of each tile is $\frac 1 2$. 
 By  Lemma~\ref{lemma:isosceleshelper3}, $q=0$, so 
 $$X = pa + rb =  \frac{p+r} {\sqrt 2}.$$
Squaring both sides, we have
$$X^2 = (p+r)^2 / 2.$$
Since $N=2X^2$, 
$$N = (p+r)^2.$$
So $N$ is a square.  That is conclusion~(ii) of the lemma. 
That completes Case~1.

\case{2} $\alpha = \frac \pi 6$.  Then $\beta = \frac \pi 3$, so $ABC$ is
equilateral.  With $c = 1$ we have $a = \sin \alpha = 1/2$ and
$b = \cos \alpha = \sqrt 3/2$, so $X = \vert AB \vert = pa + qb + r$
belongs to $\Q(\sqrt 3)$; say $X = u + v\sqrt 3$ with $u,v$ rational.
Since $N = 2X^2$, we have
$$ N/2 = X^2 = u^2 + 3v^2 + 2uv\sqrt 3, $$
so $uv = 0$.  If $v = 0$, then $X = \sqrt{N/2}$ is rational, contrary
to our assumption that $\sqrt{N/2}$ is irrational.  Hence $u = 0$ and
$N/2 = 3v^2$.  Write $v = s/t$ with $s$ and $t$ relatively prime
integers.  Then $Nt^2 = 6s^2$, so $t^2$ divides $6s^2$; since $s$ and
$t$ are relatively prime, $t^2$ divides $6$, so $t = 1$.  Hence
$N/2 = 3s^2$ is three times a square, and conclusion (iii) holds.

\case{3}  $\alpha = \frac \pi 3$.  Then $\beta = \frac \pi 6$.
Then $a =\sin \alpha = \sqrt{3}/2$ and $b =\sin \beta = 1/2$; hence 
$X = \vert AB \vert = pa + qb + r $ belongs  to $\Q(\sqrt 3)$.  
Then the proof is completed verbatim as in Case~2.
\end{proof}

\begin{lemma}\label{lemma:baseangles}
  Suppose isosceles (and not equilateral) triangle $ABC$ 
is $N$-tiled by a  tile with angles $(\alpha,\beta,\frac \pi 2)$.
 Then the 
base angles of $ABC$ are equal to $\alpha$ or to $\beta$.
\end{lemma}

\begin{proof}[Proof of Lemma~\ref{lemma:baseangles}] This is an 
immediate consequence of Laczkovich's work: the first and second lines
of Table~\ref{table:isosceles} are the only ones allowing a right-angled
tile, and the first line can apply to an isosceles $ABC$ only if 
$ABC$ and the tile are both right isosceles triangles.
 \end{proof}

\begin{theorem} \label{theorem:isosceles2}
  Suppose isosceles triangle $ABC$
is $N$-tiled by a  tile with angles $(\alpha,\beta,\frac \pi 2)$.
Then
 \smallskip

(i) $N$ is a square and $\alpha = \beta = \frac \pi 4$, or
\smallskip

(ii) $N$ is twice a square, or
\smallskip

(iii) $N = 6k^2$ and $\beta = \pi/6$ or $\alpha = \pi/6$, or
\smallskip

(iv) $N$ is an even sum of squares, so $N/2$ is also a sum of squares.
\end{theorem}

\begin{remark}
All four cases occur.  Case~(ii) is possible for every such $N$ and every
right-angled tile; case~(iii) is illustrated by the $54$-tiling in
Figure~\ref{figure:54}; and case~(iv) is possible for every such $N$, with a
suitable tile, by a double biquadratic tiling as in
Figure~\ref{figure:severalisosceles}.
\end{remark}

\begin{proof}[Proof of Theorem~\ref{theorem:isosceles2}]  By Lemma~\ref{lemma:baseangles},
the base angles of $ABC$ are equal to $\alpha$ or to $\beta$.
Since the conclusion of the theorem is insensitive to which angle
is labeled $\alpha$, we may assume the base angles are $\beta$.
By Theorem~\ref{theorem:rationalisosceles}, the
conclusion is correct when $\alpha$ is a rational multiple of $\pi$;
indeed in that case either (i), (ii), or (iii) holds.
Therefore, we may assume that 
$\alpha$ is not a rational multiple of $\pi$.
If $N/2$ is a rational square, then $2N$ is a rational square, and
an integer, hence an integer square. Hence, $N$ is even, and $N/2$
is an integer.  Since $N/2$ is a rational square and an integer,
it is also an integer square, so $N$ is twice a square, and case 
(ii) of the theorem holds.

Therefore, we may assume that 
$$\lambda = \sqrt{N/2} \mbox{ is irrational}.$$
Let $X$ be the length of $AB$.   We will show that 
$$ X = \lambda. $$
Twice the area of $ABC$ is 
$$X^2 \sin 2\alpha = 2X^2 \cos \alpha \sin \alpha = 2X^2 ab.$$
It is also $Nab$, since there are $N$ tiles each of area $ab/2$.
Therefore, $N = 2X^2$.  However, $N = 2\lambda^2$ by definition of $\lambda$,
and both $X$ and $\lambda$ are positive.  Therefore, $X = \lambda$,
as claimed.

Let $(a,b,c)$ be the sides of the tile; we may choose the scale so that $c=1$.
By Lemma~\ref{lemma:isosceleshelper3}, $a$ and 
$b$ are rational multiples of $\lambda$, and $AC$ is 
made up of $c$ edges, and the equal sides $AB$ and $BC$ are
made up entirely of $a$ and $b$ edges.
Then $a/b$ is rational; let $a/b = p/q$ in lowest terms; so $p$ and $q$ are coprime.
Then since $a^2 + b^2 = 1$, 
\begin{equation}
  a = \frac p {\sqrt{p^2+q^2}} \mbox{\ and \ }  b = \frac q {\sqrt{p^2+q^2}}.\label{eq:1434}
\end{equation}
Since $a$ and $b$ are rational multiples of $\lambda = \sqrt{N/2}$, it follows that 
$\sqrt{p^2+q^2}$ is a rational multiple of $\sqrt{N/2}$.   Let $\xi$ be a 
rational number such that 
$$  \xi \sqrt{N/2}= \sqrt{p^2 + q^2}.$$
Squaring both sides and multiplying by 2,
$$ \xi^2 N  = 2(p^2+q^2).$$
Put $\xi = u/v$ with coprime integers $u,v$.  Then 
\begin{equation}
  u^2 N = 2v^2(p^2+q^2). \label{eq:1442}
\end{equation}
Assume $N$ is odd. Then $2 \mid u^2$.  Then $u^2 N \equiv 0 \pmod 4$.  Then $2 \mid u$.
Since $u$ and $v$ are coprime, $2 \nmid v$.  Then $ p^2 + q^2 \equiv 2 \pmod 4$.
Then $p$ and $q$ are both odd.

We now show that $N$ odd is impossible.  Suppose $N$ is odd; then, as
just shown, $p$ and $q$ are both odd.  We derive a contradiction by a
parity argument based on the directions of the edges of the tiling.

\smallskip
\noindent{\em Direction classes.}  Fix the direction of $AB$ as
$\theta_0$.  Every tile edge has direction (taken mod $\pi$) of the form
$$ \theta_0 + i\alpha + j\frac \pi 2 $$
with integers $i,j$.  Indeed, the union of the tile boundaries is
connected, any two edges of one tile differ in direction by an angle of
the tile, and two tiles sharing a segment of positive length have
parallel edges there; all tile angles are integer combinations of
$\alpha$ and $\frac \pi 2$, since $\beta = \frac \pi 2 - \alpha$.
Since $\alpha$ is not a rational multiple of $\pi$, the direction
determines $i$ uniquely and $j$ modulo $2$.  Define the {\em class} of a
direction to be $i+j \bmod 2$.  Rotation by $\beta$ preserves the class,
since it changes $(i,j)$ by $(-1,1)$ or $(1,-1)$; rotation by $\alpha$
or by $\frac \pi 2$ changes it.

In each tile the angle between the $a$ edge and the $c$ edge is $\beta$,
and the $b$ edge is perpendicular to the $a$ edge.  Hence the $a$ and $c$
edges of a tile have the same class and the $b$ edge has the other class.
In particular, {\em each tile has exactly one leg ($a$ or $b$ edge) of
class $0$.}  The side $AB$ has class $0$ by choice of $\theta_0$; $AC$
makes angle $\beta$ with $AB$ and $BC$ makes angle $2\alpha$ with $AB$,
so all three sides of $ABC$ have class $0$.

\smallskip
\noindent{\em Counting.}  Rescale the tiling so that $a = p$, $b = q$, and
$c = w := \sqrt{p^2+q^2}$.  Here $w$ is irrational, because $a/c = p/w$
is a nonzero rational multiple of the irrational number $\lambda$.
Every tile edge lies either on the boundary of $ABC$ or on exactly one
maximal internal segment, and the edges on each side of a maximal
internal segment $\sigma$ exactly cover $\sigma$.  If side~$i$ of
$\sigma$ ($i = 1,2$) carries $n_i$ legs of total length $L_i$
(an integer) and $k_i$ hypotenuses, then
$$ L_1 + k_1 w = L_2 + k_2 w, $$
so $k_1 = k_2$ and $L_1 = L_2$, since $w$ is irrational.  Since $p$ and
$q$ are odd, $n_i \equiv L_i \pmod 2$; hence $n_1 + n_2 \equiv 2L_1
\equiv 0 \pmod 2$.

Now count the legs of class $0$.  By the italicized statement above
there are exactly $N$ of them.  Each lies either on a maximal internal
segment of class $0$, or on $AB$ or $BC$ (there are no legs on $AC$,
which is composed of $c$ edges).  By the previous paragraph, the
internal segments contribute an even number.  By
Lemma~\ref{lemma:isosceleshelper3}, $AB$ and $BC$ carry no $c$ edges, so
if $AB$ carries $m$ legs and $BC$ carries $m'$ legs, then
$m \equiv X \equiv m' \pmod 2$, where $X = |AB| = |BC|$ is an integer
in the present scale.  Therefore
$$ N \equiv m + m' \equiv 2X \equiv 0 \pmod 2, $$
contradicting the assumption that $N$ is odd.  Hence $N$ is even.

\smallskip
Finally, by (\ref{eq:1442}),
$$ N = \left(\frac v u\right)^2 \bigl((p+q)^2 + (p-q)^2\bigr)
\quad\mbox{and}\quad
\frac N 2 = \left(\frac v u\right)^2 (p^2+q^2), $$
so $N$ and the integer $N/2$ are sums of two rational squares, hence
of two integer squares by Lemma~\ref{lemma:rationalsquares}.  Thus
(iv) of the lemma holds.
\end{proof}

\begin{corollary} \label{lemma:notprime}  Suppose isosceles
triangle $ABC$ is $N$-tiled by a right triangle and $N>2$.  Then $N$ is 
not a prime.  Moreover, $N$ is not twice a prime $p$ with $p \equiv 3 \pmod{4}$
except for $N=6$. 
\end{corollary}

\begin{proof} That $N$ is not prime follows immediately
from 
Theorem~\ref{theorem:isosceles2}.
The second part follows from the theorem and 
 Lemma~\ref{lemma:sumsofsquares}. 
 \end{proof}

Theorem~\ref{theorem:isosceles2} 
gives necessary and sufficient conditions on $N$ for the existence of an $N$-tiling
of {\em some} isosceles $ABC$ by a right-angled tile, if $N$ is even.
It remains to specify exactly {\em which} isosceles $ABC$ can be 
$N$-tiled, when $N$ is given.
  The following theorem spells it out.
 
\begin{theorem}\label{theorem:whichABC}
Given a positive integer $N>1$, the possible isosceles
triangles $ABC$ that can be $N$-tiled by a right-angled tile
are as follows. Here the sides of the tile are $(a,b,c)$
and the angles are $(\alpha,\beta,\frac{\pi}{2})$.

\begin{enumerate}
\renewcommand{\labelenumi}{(\roman{enumi})}
\item If $N/2$ is a square, any isosceles triangle can be
$N$-tiled (by a double quadratic tiling).

\item If $N/2$ is a sum of two squares but not a square, then an isosceles
triangle $ABC$ with base angles $\beta$ can be $N$-tiled with tile
$(\alpha,\beta,\frac{\pi}{2})$ if and only if
$$ \tan\beta=\frac{r}{p}, \qquad\text{where}\qquad \frac{N}{2}=r^2+p^2 $$
for some nonnegative integers $r$ and $p$.

\item If $N$ is a square, the right isosceles triangle $ABC$
can be $N$-tiled by a quadratic tiling.

\item If $N$ is six times a square, then the isosceles
triangle with base angles $\frac{\pi}{6}$ can be $N$-tiled
by the tile with $\alpha=\frac{\pi}{3}$ and
$\beta=\frac{\pi}{6}$.

\item If none of the above applies, then no isosceles triangle
can be $N$-tiled by any tile.
\end{enumerate}
\end{theorem}

\begin{remark}  Since $N/2$ may sometimes be 
expressible in more than one way as a sum of two squares,
there can sometimes be more than one possible $ABC$ and 
tile for a given $N$, but only finitely many.  Moreover,
if $N$ is both a square and a sum of squares, there are 
more possibilities, as in Figure~\ref{figure:severalisosceles}.  It will
be very difficult to provide a full characterization of all
$N$-tilings.
\end{remark}

\begin{proof}  [Proof of Theorem~\ref{theorem:whichABC}]
(i).  Just divide $ABC$ by 
its altitude $BD$ and tile each half with a quadratic tiling.
\smallskip

(ii).  First suppose $N/2 = r^2 + p^2$.  Choose the tile with legs
proportional to $p$ and $r$, so that $\tan \beta = r/p$.  Then each
half of $ABC$, cut off by the altitude from $B$, is a right triangle
similar to the tile, and it can be $(r^2+p^2)$-tiled by a biquadratic
tiling.  Combining the tilings of the two halves gives an $N$-tiling
of $ABC$.

Conversely, suppose $ABC$ is $N$-tiled by the tile
$(\alpha,\beta,\frac \pi 2)$, where $N/2$ is a sum of two squares but
not a square.  We must show that $\tan\beta = r/p$ for some
nonnegative integers $r,p$ with $N/2 = r^2 + p^2$.

First suppose $\alpha$ is a rational multiple of $\pi$.  By
Theorem~\ref{theorem:rationalisosceles}, since $N/2$ is not a square,
either $N = n^2$ and $\beta = \pi/4$, or $N/2 = 3k^2$.  The latter is
impossible, since $3k^2$ contains an odd power of $3$ and so is not a
sum of two squares.  In the former case, $N/2 = n^2/2$ is an integer,
so $n$ is even, and $N/2 = (n/2)^2 + (n/2)^2$ with
$\tan \beta = 1 = (n/2)/(n/2)$, as claimed.

Now suppose $\alpha$ is not a rational multiple of $\pi$.  Since $N/2$
is not a square, $\lambda = \sqrt{N/2}$ is irrational, as shown in the
proof of Theorem~\ref{theorem:isosceles2}; also $N$ is even, by that
theorem.  By Lemma~\ref{lemma:isosceleshelper3}, $a/b$ is rational, the
base $AC$ is composed only of $c$ edges, and $AB$ is composed only of
$a$ and $b$ edges.  Write $a/b = e/f$ in lowest terms, and rescale
the tiling so that
$$ a = e, \qquad b = f, \qquad c = w := \sqrt{e^2+f^2}. $$
Since $b$ is opposite $\beta$, we have $\tan\beta = b/a = f/e$ and
$\cos\beta = a/c = e/w$.  Let $X = |AB|$ and $Y = |AC|$.  Then $X$ is
an integer, being a sum of $a$ and $b$ edges, and $Y = sw$, where
$s$ is the number of $c$ edges on $AC$.  Since the base angles of
$ABC$ are $\beta$, we have $Y = 2X\cos\beta$, that is,
$sw = 2Xe/w$, so
$$ X = \frac{s(e^2+f^2)}{2e}. $$
The height of $ABC$ is $X\sin\beta = Xf/w$, so the area of $ABC$ is
$\frac 1 2\, sw \cdot Xf/w = sXf/2$.  Since the area of each tile is
$ef/2$, we have $Nef/2 = sXf/2$, so $Ne = sX$.  Hence
$$ \frac N 2 = \frac{sX}{2e} = t^2(e^2+f^2), \qquad\text{where}\qquad
   t = \frac{s}{2e}. $$
We show that $t$ is an integer.  Write $t = g/h$ in lowest terms.
Since $h$ divides $2e$, and $t(e^2+f^2) = X$ is an integer, $h$ also
divides $e^2+f^2$.  Since $e$ and $f$ are relatively prime,
$\gcd(e, e^2+f^2) = 1$, so $h$ divides $2$.  Suppose $h = 2$.  Then
$e^2+f^2$ is even, so $e$ and $f$ are both odd and
$e^2+f^2 \equiv 2 \pmod 4$; also $g$ is odd.  Then
$N = 2t^2(e^2+f^2) = g^2(e^2+f^2)/2$ is odd, contradicting the fact
that $N$ is even.  Hence $h = 1$ and $t$ is an integer.  Setting
$r = tf$ and $p = te$, we have
$$ \frac N 2 = r^2 + p^2 \qquad\text{and}\qquad \tan\beta = \frac f e = \frac r p, $$
as claimed. 
\smallskip

(iii).  Any triangle can be quadratically tiled.
\smallskip

(iv).  If $N = 6k^2$
then a tiling with $\beta = \frac \pi 6$ is possible; see
 Figure~\ref{figure:54}.  It remains to show that no other tile
 is possible.  Let $N=6k^2$.  Then $N$ is not a square, and 
 $N$ is not twice a square.  Since $N$ contains an odd power of 3,
 $N$ is not a sum of two squares, and the same is true of $N/2$.
 Hence, no other case in the theorem can apply, and $\beta = \frac \pi 6$
 is the only possibiity.
\smallskip

(v).   By Theorem~\ref{theorem:isosceles2}, these
cases are exhaustive.
\end{proof}

\section{Possible Values of $N$ in Tilings with Commensurable Angles}
We wish to add a third column to Laczkovich's Table~\ref{table:isosceles},
giving the possible forms of $N$ if there is an $N$-tiling of an isosceles triangle $ABC$
by the tile in that row.  For example, when $ABC$ is similar to the 
tile, then $N$ must be a square, so we put $n^2$ in the third column.
While we are at it, we add a fourth column with a citation to the 
result about the value of $N$.   The fact that the rows of the table
are exhaustive is due to \cite{laczkovich2012}, as explained in \cite{beeson-noseven}.
Laczkovich did not set forth this exact table, and the ``explanation'' is just
to set out how the original table of possibilities from \cite{laczkovich1995} grew smaller as certain
lines were shown not to correspond to actual tilings.

  The revised and extended table is Table~\ref{table:laczkovich-extended}.
All the entries in this table give necessary and
sufficient conditions on $N$ for the tilings to exist.  
  That the condition in that last row is 
necessary we prove in the following lemma:

\begin{lemma}\label{lemma:equilateral120}
If an equilateral triangle is $N$-tiled by a tile with angles
$(\frac \pi 6, \frac \pi 6, \frac{2\pi}3)$, then $N$ is three times a
square.
\end{lemma}

\begin{proof}
Scale so that the tile has legs $1$ and base $\sqrt 3$; its area is
$\frac 1 2 \sin \frac{2\pi} 3 = \sqrt 3/4$.  Let $X$ be the side of the
equilateral triangle.  Its area is $\sqrt 3 X^2/4$, so $X^2 = N$.  Each
side of the triangle is covered by legs and bases of tiles, so
$X = p + q\sqrt 3$ for some nonnegative integers $p,q$, the numbers of
legs and bases on that side.  Then
$$ N = X^2 = p^2 + 3q^2 + 2pq\sqrt 3, $$
so $pq = 0$.

Suppose $q = 0$, so that a side of the triangle supports only legs.  A
leg joins the $\frac{2\pi}3$ vertex of its tile to a $\frac \pi 6$
vertex, so each of the $p$ tiles supported by that side has its
$\frac{2\pi}3$ angle at a vertex on that side.  The angle $\frac \pi 3$
at a corner of the triangle cannot contain a $\frac{2\pi}3$ angle, so
these $p$ angles lie at the $p-1$ vertices in the interior of the
side.  Hence some vertex has two angles $\frac{2\pi}3$ on the same side
of a line, which is impossible since $\frac{4\pi}3 > \pi$.  Therefore
$q > 0$ and $p = 0$, and $N = 3q^2$.
\end{proof}

\begin{theorem} \label{theorem:laczkovich-extended}
Suppose $(\alpha,\beta,\gamma)$ are all rational multiples of 
$2\pi$, and the isosceles or equilateral triangle $ABC$ is $N$-tiled by a tile with angles 
$(\alpha,\beta,\gamma)$.  Then $ABC$, $(\alpha,\beta,\gamma)$,
and $N$ correspond to one of the lines in Table~\ref{table:laczkovich-extended}.

\begin{table}[ht]
\begin{center}
\setlength{\extrarowheight}{.5em}
\begin{tabular}{rrrl}
$ABC$  &  the tile & form of $N$ & citation    \\
\hline
$(\beta,\beta,2 \alpha)$  &  similar to $ABC$ & $n^2$ & \cite{snover1991} \\
$(\beta,\beta,2\alpha)$ &  $(\alpha, \beta, \frac \pi 2)$ & $2n^2$ & Theorem~\ref{theorem:rationalisosceles}\\
{\large $(\frac \pi 6, \frac \pi 6, \frac {2 \pi} 3)$} &{\large  $(\frac \pi 6,\frac \pi 3, \frac \pi 2)$} & $6n^2$ & Theorem~\ref{theorem:rationalisosceles}\\
equilateral  & {\large $(\frac \pi 6, \frac \pi 3, \frac \pi 2)$} & $2n^2$ or $6n^2$ & Theorem~\ref{theorem:rationalisosceles} \\
equilateral  & {\large $(\frac \pi 6, \frac \pi 6, \frac {2\pi} 3)$ }
                  & $3n^2$ &  \cite{laczkovich2012} and Lemma~\ref{lemma:equilateral120}
\end{tabular}
\end{center}
\caption{$N$-tilings by tiles with commensurable angles, with form of $N$.
The rows are not mutually exclusive: for example, the triangle
$(\frac \pi 6, \frac \pi 6, \frac {2\pi} 3)$ also admits $2n^2$-tilings,
by the second row.}
\label{table:laczkovich-extended}
\end{table}
\end{theorem}

\begin{proof}  As discussed above, Laczkovich characterized
the pairs of tiled triangle and tile, as given in
Table~\ref{table:isosceles}.%
\footnote{\, Again, we remind readers who may check with 
\cite{laczkovich1995} that there are four entries in Laczkovich's 
Theorem~5.1 that are shown in the subsequent Theorem~5.3 not to apply
to tilings by congruent triangles, so they do not appear in our tables.
}
It remains to characterize the 
possible $N$ for each line.  In several cases, lines in 
Table~\ref{table:isosceles} split into two or more lines in 
Table~\ref{table:laczkovich-extended}, which supplies the required
possible forms of values of $N$.  That table
 lists in its last column citations
to the literature or theorems in this paper for each line.
\end{proof}

\section{Tilings with $\gamma = 2\alpha$ and $(a,b,c)$ Commensurable}
In this section and the next,
we take up the row of Laczkovich's second table in 
which $ABC$ is isosceles with base angles $\alpha$ and is tiled
by a tile with $\gamma = 2\alpha$, and $\alpha$ is not a rational
multiple of $\pi$.    The condition $\gamma = 2\alpha$ can also be 
written as $3\alpha + \beta = \pi$.  Unlike the similar-looking
condition $3\alpha + 2\beta = \pi$, this condition does not imply 
$\gamma > \pi/2$.  The vertex angle of $ABC$ is then $\pi-2\alpha = \alpha +\beta$.  
The tile  $(4,5,6)$ satisfies $\gamma = 2\alpha$; this is shown below as 
an example of Lemma~\ref{lemma:luthar}. 

The numbers $(a,b,c)$ are called {\em commensurable} if their ratios
are rational.  In that case we say the {\em tile is rational}.
 If the edges of a triangle are commensurable, then the 
triangle is similar to one with integer edges.  The remarkable fact
is that, if $ABC$ is isosceles with base angles $\alpha$ and vertex
angle $\alpha+\beta$,  then it can be $N$-tiled for some $N$ by 
a tile with angles $(\alpha,\beta,\gamma)$ if and only if the 
edges of the tile are commensurable.  This fact is really two 
different theorems:
\begin{itemize}
\item If the tile is rational, there is an $N$-tiling for some $N$, and 
\item If the tile is not rational, there is no $N$-tiling for any $N$.
\end{itemize}
The first statement is due to Laczkovich \cite{laczkovich1995}.
We will explain his proof in this section.  The second statement
is proved in the next section.

Laczkovich  \cite{laczkovich1995}  proves that an isosceles triangle 
with angles as described, can 
be dissected into triangles {\em similar} to the tile, plus 
one parallelogram; then using the commensurability condition, these
triangles and the parallelogram can all be tiled with the same size 
of tile.  The only problem is that the tile will have to be 
{\em really tiny}.

\begin{figure}[ht]
\begin{center}
\includegraphics[width=0.9\textwidth]{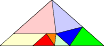}
\end{center}
\caption{Laczkovich's dissection of isosceles $ABC$ into triangles similar to $(4,5,6)$ and a parallelogram.}
\label{figure:bigisosceles}
\end{figure}

\begin{figure}[ht]
\begin{center}
\includegraphics[width=0.9\textwidth]{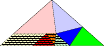}
\end{center}
\caption{We look for a tiling starting like this.}
\label{figure:bigisosceles2}
\end{figure}

We illustrate this with the tile $(4,5,6)$.
 The idea of Laczkovich's construction (Figure~3 in 
\cite{laczkovich1995}) is shown in Figure~\ref{figure:bigisosceles}.
Laczkovich's idea is to quadratically tile each triangle, and then
tile the parallelogram.  As Laczkovich pointed out,
the commensurability of the tile edges mean that with a small enough
tile, this will succeed.
We illustrate the idea in
 Figure~\ref{figure:bigisosceles2}.
Observe that the tiles
shown in that figure will not work, because with that size of tile,
we cannot continue the tiling into the next (blue) triangle, as 
if the tile is $(4,5,6)$, the boundary between light green and 
blue is five 6-edges, for a total of 30, which cannot be made of 4-edges,
as it would have to be to tile the blue triangle.  Clearly we 
should have chosen a smaller tile,  for example half that size.
With a tile half that size, we run into similar trouble at the 
next boundary.  To choose the tile correctly, we introduce a variable
for each triangle to count the number of tiles on each side of that
triangle. Then there is a linear equation at each boundary.  If we 
assume that the parallelogram will be tiled by $n$ tiles on its 
diagonal side and $m$ on its horizontal side, then these variables
will satisfy the following equations.  The equations show that 
everything is determined once the number of tiles on each side of 
the red triangle is chosen.%
\footnote{For readers without colors: red is the triangle to the right
of the parallelogram, and the other colors are in counterclockwise order
from red.}
\begin{eqnarray*}    
	&red &= \ 1, \\
	&orange &= \ 4 \ red/6, \\
	&lightgreen &= \ 5\  orange/4, \\
	&blue &= \ 6 \ lightgreen/4, \\
	&green &= \ 5\  blue/4, \\
	&lightblue &= \ 6\  blue/4, \\
	&pink &= \ 5\  lightblue/4, \\
	&m &=\  (5\  pink - 4\  orange)/6, \\
	&n &= \ (5\  red)/5,
\end{eqnarray*}
and the total number of tiles is
\begin{eqnarray*}
	&total &= \ red\,^2 + orange^2 + blue^2 +\, lightgreen^2  \\
	&&  \quad +\, green^2 
	      + lightblue^2 + pink^2 + 2mn.
\end{eqnarray*}
Solving these equations, starting with $red = 1$, we get these 
answers (in the order of variables listed above):
$$ 1,\quad 
\frac 2 3,\quad  
\frac 5 6,\quad  
\frac 5 4,\quad  
\frac {25} {16},\quad  
\frac {15} 8\quad  
\frac {75} {32},\quad  
\frac {869} {576},\quad
1.
$$
This reveals that Figure~\ref{figure:bigisosceles2} is misleading,
in that the parallelogram of Figure~\ref{figure:bigisosceles} is 
not accurately tiled--the tiled parallelogram is off by less than
a pixel.  To get it accurately tiled requires a {\em much} smaller
tile.  To clear those denominators we have to start with 
$ red = 576$ instead of $ red = 1$.  Then we get
$$576, \quad 384, \quad 480, \quad 720, \quad 900, \quad 1080, \quad 1350, \quad 869, \quad 576. $$
So the number of tiles required (namely the sum of the squares of 
the color numbers plus $2mn$) is $6028020$.

Here matters stood for about eight years.  Then, in 2024, Bryce Herdt found 
$N$-tilings for $N = 1125$ and then $720$.  These tilings are exhibited in 
Section \ref{section:herdt}.

\section{No Tilings with $\gamma = 2\alpha$ and $(a,b,c)$ Incommensurable}
In this section, we rule out tilings of an isosceles triangle
with base angles $\alpha$ and vertex angle $\alpha+\beta$, in case
the tile edges are incommensurable.  To state the point another way,
if there is a tiling of such an isosceles $ABC$, then the tile 
must be rational.

\subsection{Stars and Centers}
Suppose isosceles triangle $ABC$ is tiled by a tile with angles 
$(\alpha,\beta,\gamma)$, not a right triangle,  
and $\alpha$ is not a rational
multiple of $\pi$.  We will consider and analyze the possible 
configurations formed by tiles at a vertex of the tiling.
We begin by ruling out certain possibilities.

\begin{lemma}\label{lemma:twobeta}
Let isosceles triangle $ABC$ with base angles $\alpha$ (at $A$ and $C$) be  \mbox{$N$-tiled}
by a tile with angles $(\alpha,\beta,2\alpha)$, 
and suppose that the tile is
not a right triangle and $\alpha$ is not a rational multiple of $\pi$.
Let $P$ be a vertex on the boundary of $ABC$. Then there are not 
two $\beta$ angles of tiles at $P$, and there are not two 
$\gamma$ angles of tiles at $P$.
\end{lemma}

\begin{proof}  Suppose, for proof by contradiction, that
two tiles have their $\beta$ angles at the same boundary vertex.
Then for some nonnegative integers $u,v,w$ we have
$$
\pi = u \alpha + (v+2) \beta + w \gamma.
$$
Since $\gamma = 2\alpha$, 
$$
\pi = (u+2w)\alpha + (v+2) \beta.
$$
Since $3\alpha + \beta = \pi$,
$$
0 = (u+2w-3) \alpha + (v+1)\beta.
$$
Solving for $\beta$, 
$$
\beta =\left( \frac{3-u-2w}{v+1}\right) \alpha.
$$
Putting that value for $\beta$ into $3\alpha+\beta = \pi$, we can
solve for $\alpha/\pi$:
$$ \alpha/\pi = 3 + \left( \frac{3-u-2w}{v+1}\right). $$
That is rational, contradicting the hypothesis that $\alpha$
is not a rational multiple of $\pi$.  That completes
the proof that two $\beta$ angles do not occur at the same boundary vertex.
\smallskip

The proof for two $\gamma$ angles is similar.  First, if there
are no $\beta$ angles, then 
\begin{eqnarray*}
\pi &=& u \alpha + (w+2) \gamma \\
&=&  (u+ 2w+4) \alpha,
\end{eqnarray*}
which is a contradiction, since $\alpha$ is not a rational multiple of $\pi$.
And if there is one $\beta$ angle, then 
$$ \pi = u \alpha + (w+2) \gamma + \beta. $$
Hence, subtracting $\alpha + \beta + \gamma = \pi$, we have
$$ 0 = (u-1) \alpha + (w+1) \gamma.$$
Since $\gamma = 2\alpha$,
$$ 0 =  (u-1 + 2(w+1)) \alpha =  (u+2w+1) \alpha,$$
which is a contradiction, since $u+2w+1 > 0$.
Since we already proved there cannot be more than one $\beta$,
we are finished. 
\end{proof}
\smallskip

  We define the {\em angle sum} of a vertex
to be the sum of the angles of the tiles sharing that vertex.
Except for the vertices $A$, $B$, and $C$, that angle sum will 
always be either $\pi$ or $2\pi$.  It will be $\pi$ if and only 
if the vertex lies on the interior of the boundary of a tile or of 
$ABC$.  For short, we refer to a vertex with angle sum $\pi$ as a
{\em boundary vertex}, though it need not be on the boundary of $ABC$.

Consider a boundary vertex.
A {\em normal boundary vertex} has three tiles, with angles 
$\alpha$, $\beta$, and $\gamma$.   A {\em star} has 
three $\alpha$ angles and a $\beta$.   These are the only possibilities
for a boundary vertex, since $\alpha$ is not a rational multiple of $\pi$,
as spelled out in Lemma~\ref{lemma:twobeta}.
Next, consider an interior vertex.   A {\em normal interior vertex} has 
two each of $(\alpha, \beta, \gamma)$ angles.  A {\em center} has
three $\gamma$ and two $\beta$ angles.  For example, 
there is a center
in Figure~\ref{figure:bigisosceles}, more or less in the center of the figure.
There may also be interior vertices other than centers that are not 
normal; these will have either $4\alpha + 2\beta + \gamma$ or $6\alpha + 2\beta$.
These vertices we call {\em interior stars}.  The case of angles
$4\alpha +2 \beta + \gamma$ we call a {\em single interior star} and 
the other case is a {\em double interior star}.

\begin{lemma} \label{lemma:starsandcenters}
Suppose isosceles triangle $ABC$ is tiled by a tile with angles 
$(\alpha,\beta,\gamma)$ with $\alpha$ not a rational multiple of $\pi$.
Let $\C$ be the number of centers and ${\mathcal S}$ the number of stars,
counting double interior stars twice. 
Then ${\mathcal S} + 1  = \C$.  In particular there is at least
one center.
\end{lemma}

\begin{example}  When the tiling begun 
in Figure~\ref{figure:bigisosceles} 
is completed, there will be one center and zero stars. All 
the vertices introduced by quadratic tilings will be normal vertices.
If we then combine four copies of this tiling to create a quadratic
tiling of a triangle twice the size, there will be four centers,
balanced by three interior stars in the midpoints of the sides,
where the three of the four copies have common vertices, and the 
vertices of $ABC$, $\C = 4 =  {\mathcal S} + 1$.
\end{example}

\begin{proof} [Proof of Lemma~\ref{lemma:starsandcenters}]
 Each tile has one each of $\alpha$, $\beta$, and 
$\gamma$ angles.  At the vertices of $ABC$ we have three $\alpha$ angles 
and one $\beta$ angle (just as we have at a star).  Counting the 
vertex angles we have equal numbers of $\alpha$, $\beta$, and $\gamma$ 
angles at each normal vertex.  At each vertex we define the 
{\em excess} or {\em deficit} of each of $(\alpha,\beta,\gamma)$ to 
be the difference between the number of those angles at the vertex
and the number at a normal vertex.   At a star we have two excess $\alpha$ 
angles and a deficit of one $\gamma$ angles.  At a single interior
star the same applies; at a double interior star we have double that
contribution.   At a center we 
have an excess of one $\gamma$ and a deficit of two $\alpha$.
 Adding up the excesses and deficits from the 
vertices of the tiling, including $A$, $B$, and $C$, we must get zero.
The vertices of $ABC$ count the same as a star.  One center will 
``balance'' one star, in the sense that the deficits and excesses
add to zero.  For example, in Figure~\ref{figure:bigisosceles}, we 
have one center, and no stars; so the center balances the vertices
of $ABC$, which count as a star. 

If there are 
no interior stars, then the number of stars, plus one for $ABC$,
will equal the number of centers.   If, however, there are interior
stars, those will require additional centers to balance them, one 
center for each single interior star and two for each double interior 
star.   Since we defined ${\mathcal S}$ by this double-counting of double interior
stars, we still have $\C = {\mathcal S}+1$. 
\end{proof}

\subsection{The Tile is Rational}

If $Q$ is a vertex of a tiling,  and $QR$ is an internal segment 
of the tiling supporting a tile on one side with its $a$ edge on $QR$
and a vertex at $Q$, and supporting a tile on the other side with 
its $b$ edge on $QR$ and a vertex at $Q$, then we say $QR$ is an 
{\em $a/b$ edge}, or an {\em $a/b$ edge at $Q$}.   We define {\em $a/c$ edge} and {\em $b/c$ edge} similarly.    Note that an $a/b$ edge
is also a $b/a$ edge.

\begin{lemma} \label{lemma:centeroutgoing}
Let the isosceles triangle $ABC$ with base angles $\alpha$  be $N$-tiled by
a tile with angles $(\alpha,\beta, 2\alpha)$,
with $\alpha$ not a rational multiple of $\pi$.
Let $Q$ be a center in the tiling.  Then 
either  there is  an $a/b$ edge
at $Q$ or there are both $a/c$ and $b/c$ edges at $Q$.
\end{lemma}

\begin{proof}  Assume there is no $a/b$ edge.
We must prove there is an $a/c$ and a $b/c$ edge at $Q$.
 Each tile with a 
vertex at the center  $Q$ has an $a$ edge ending at $Q$, since the angles at $Q$
are all $\beta$ or $\gamma$.  At a center, five tiles meet, 
so that is a total of five $a$ edges.
Since five is an odd number, these edges cannot all be paired with 
other $a$ edges.  We have assumed there is no $a/b$ edge;
hence there is an $a/c$ edge.
Similarly, there are three $b$ edges ending at $Q$, belonging to the 
tiles with their $\gamma$ angles at $Q$.  Since three is odd, these
cannot each be paired with another $b$ edge. Since there are
no $a/b$ edges, there is a $b/c$ edge. 
\end{proof} 

\begin{lemma} \label{lemma:essentialsegments}
Let the isosceles triangle $ABC$ with base angles $\alpha$  be $N$-tiled by
a tile with angles $(\alpha,\beta, 2\alpha)$ with $\alpha$ not a rational multiple of $\pi$.
Then the tiling contains essential segments with associated relations
$jb = ua + vc$ and  $Ja = Ub + Vc$.
\end{lemma} 

\begin{proof}  Suppose, for proof by contradiction,
that there are no such essential segments. 
We consider the directed graph $\Gamma_b$ defined in 
Definition~\ref{definition:Gamma}.  We wish to identify the terminal 
links in this graph.  To that end we must consider the possible 
configurations that can arise when an internal segment $UQ$  of the tiling
supports on the same side a series of (one or more) tiles with their
$b$ edges, followed by a tile with an $a$ or $c$ edge. 
Let $PQV$ be the three successive vertices, with $PQ$ of length $b$ 
and $QV$ of length $a$ or $c$.  
 We consider all possible configurations in which $Q$ is a vertex where three 
tiles on one side of a line $PQ$ each have a vertex at $Q$, contributing
one angle each, so the angle sum at $Q$ is 
 $\alpha +\beta + \gamma$ on one side of $PQ$.  All these configurations 
 are shown in Figure~\ref{figure:GammaB}. 
  The figure shows that in each of those cases, 
there is a unique
 outgoing segment $QT$ that is a $b/c$ edge or a $b/a$ 
 edge.  If this segment is extended far enough,
 we will come to the last
$b$ edge on the side that has a $b$ edge at $Q$.  Since there are
no essential segments, that point cannot be a vertex of a tile
on the other side of $QT$, so $QT$ is an outgoing link in $\Gamma_b$. 

\begin{figure}[ht]
\includegraphics[width=0.3\textwidth]{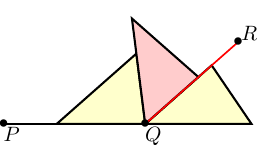}
\includegraphics[width=0.3\textwidth]{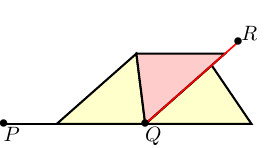} 
\includegraphics[width=0.3\textwidth]{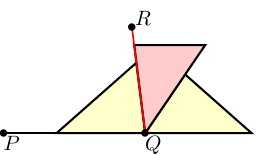} \\
\includegraphics[width=0.3\textwidth]{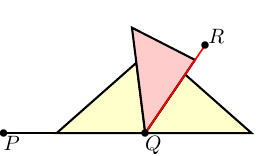} 
\includegraphics[width=0.3\textwidth]{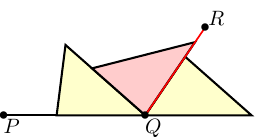} 
\includegraphics[width=0.3\textwidth]{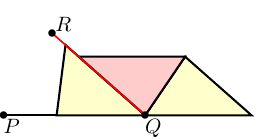} \\
\includegraphics[width=0.3\textwidth]{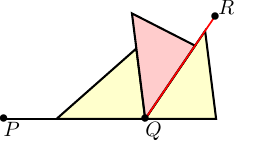} 
\includegraphics[width=0.3\textwidth]{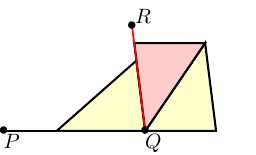} 
\includegraphics[width=0.3\textwidth]{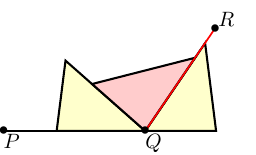} \\
\includegraphics[width=0.3\textwidth]{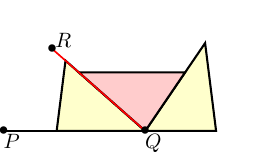} 
\includegraphics[width=0.3\textwidth]{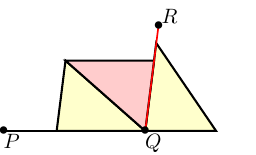} 
\includegraphics[width=0.3\textwidth]{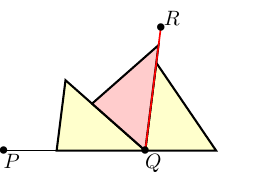} 
\caption{A link $PQ$ in $\Gamma_b$ gives rise to another link through $QR$.}
\label{figure:GammaB}
\end{figure} 

On the other hand, if $Q$ is a star, the possible configurations
are more complicated, and some of them have zero outgoing $b/c$
or $b/a$ edges, 
while others have two.  It turns out that we do not need to make use of 
that fact, so we do not give a diagram of these configurations. 

If $PQ$ is a link, then the line $PQ$ extends past $Q$ as an interior
segment of the tiling, and $Q$ lies on the interior of an edge of a
tile on Side~2 of $PQ$.  Consequently:
\smallskip

(a) no link ends at a point of the boundary of $ABC$; in particular,
the in-degree of $A$, $B$, $C$, and of any star on the boundary of
$ABC$, is zero;
\smallskip

(b) no link ends at a vertex with angle sum $2\pi$, since such a vertex
does not lie on the interior of any tile edge.  In particular, the
in-degree of every center and of every interior star is zero;
\smallskip

(c) at any vertex $Q$, at most one link ends at $Q$, since $Q$ lies on
the interior of the edge of exactly one tile.
\smallskip

\noindent
We now compare out-degree and in-degree vertex by vertex.  At a normal
vertex with angle sum $\pi$, if there is an incoming link $PQ$, there
is an outgoing link $QR$, as shown above; so the out-degree is greater
than or equal to the in-degree.  At a normal vertex with angle sum
$2\pi$, the in-degree is zero by (b).  At $A$, $B$, and $C$, and at
every interior star, the in-degree is zero, so again the out-degree is
greater than or equal to the in-degree.  At a star (angle sum $\pi$),
the in-degree is at most $1$ by (c), so the out-degree minus the
in-degree is at least $-1$.  A star may have an incoming link, and
some configurations at a star have no outgoing link; we do not need
to analyze them.  At a center $Q$, the in-degree is zero by (b).  By
Lemma~\ref{lemma:centeroutgoing}, there is an $a/b$ edge or a $b/c$
edge at $Q$.  Extending it to the last $b$ edge on the side that has a
$b$ edge at $Q$, we reach a point that is not a vertex of a tile on the
other side, since there are no essential segments; hence there is an
outgoing link at $Q$, and the out-degree minus the in-degree is at
least $1$.

Let $\C$ be the number of centers, and let ${\mathcal S}$ be as in
Lemma~\ref{lemma:starsandcenters}.  The number of stars is at most
${\mathcal S}$.  Adding up over all vertices, the total out-degree
minus the total in-degree is at least
$$ \C - {\mathcal S} \ =\ 1, $$
by Lemma~\ref{lemma:starsandcenters}.  Since each link has one head
and one tail, the total out-degree equals the total in-degree, which is
a contradiction.

We have reached a contradiction from the assumption that there is
no essential segment with associated relation of the form $jb = ua+vc$.
Hence, there is such a segment.

\begin{figure}[ht]
\includegraphics[width=0.3\textwidth]{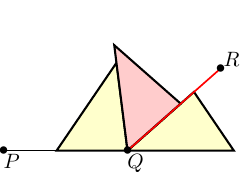}
\includegraphics[width=0.3\textwidth]{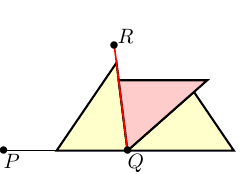} 
\includegraphics[width=0.3\textwidth]{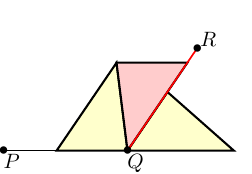}\\
\includegraphics[width=0.3\textwidth]{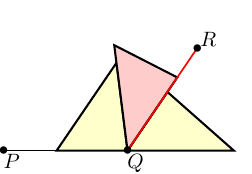} 
\includegraphics[width=0.3\textwidth]{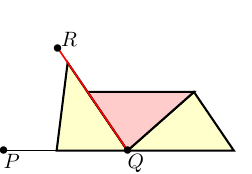} 
\includegraphics[width=0.3\textwidth]{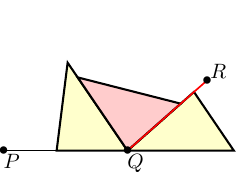} \\
\includegraphics[width=0.3\textwidth]{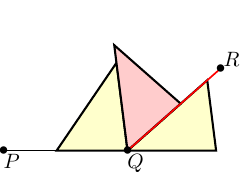} 
\includegraphics[width=0.3\textwidth]{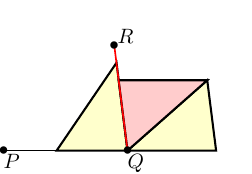} 
\includegraphics[width=0.3\textwidth]{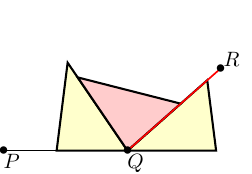} \\
\includegraphics[width=0.3\textwidth]{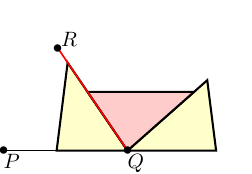} 
\includegraphics[width=0.3\textwidth]{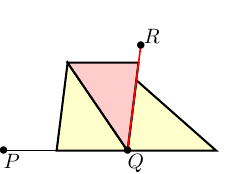} 
\includegraphics[width=0.3\textwidth]{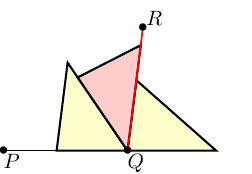}
\caption{A link $PQ$ in $\Gamma_a$ gives rise to another link through $QR$.}
\label{figure:GammaA}
\end{figure}

Next we prove the existence of essential segments with relations
of the forms $Ja = Ub + Vc$.   This is proved 
in the same way, using the graph $\Gamma_a$ instead of 
$\Gamma_b$.  Again there
is an outgoing link from each center, since 
by Lemma~\ref{lemma:centeroutgoing}, there is either an $a/b$ or an 
$a/c$ edge at $Q$, and that edge is part of an outgoing link since 
there are no essential segments.   Again we have to prove that
at a normal boundary vertex there is a unique outgoing link. 
See Figure~\ref{figure:GammaA}.
\end{proof}

\begin{lemma} \label{lemma:extendlink4} 
Let $ABC$ be an isosceles triangle with base angles $\alpha$,
tiled by tile $(\alpha,\beta,\gamma)$ with $\gamma = 2\alpha$.
Let $PQR$ be an (internal or boundary) segment of the tiling with only $c$ 
edges on one side of $PQ$, such that the tile on that side supported
by $QR$ with a vertex at $Q$ has an $a$ or $b$ edge on $QR$.
Then there is a $c/a$ or $c/b$ edge emanating from $Q$.
\end{lemma}

\begin{proof} We only consider tiles on the one side 
of $PQR$ mentioned in the lemma.  There are two cases:
either the vertex angles of those tiles at $Q$ are $\{\alpha,\beta,\gamma\}$,
or there are three $\alpha$ angles and one $\beta$ angle.

\case{1} There are three tiles at $Q$, and one of them has its $\gamma$ angle
at $Q$. Then that tile has no $c$ edge at $Q$. The other two have a 
$c$ edge ending at $Q$.  One of those lies on $PQ$.  The
other does {\em not} lie on $QR$, by hypothesis.  Since there
is no other $c$ edge ending at $Q$, that second $c$ edge forms
either a $c/b$ edge or a $c/a$ edge.

\case{2} There are four tiles at $Q$, with one $\beta$ and three $\alpha$ angles at $Q$.
 Each of the four tiles has a $c$ edge ending at $Q$.
One lies on $PQ$ and three lie on interior segments.  Since three is odd,
one of those $c$ edges is not paired with another $c$ edge, and hence
constitutes a $c/a$ segment or a $c/b$ segment.  
\end{proof}

\begin{lemma} \label{lemma:crelation4}
Let the isosceles triangle $ABC$ with base angles $\alpha$ (at $A$ and $C$) be $N$-tiled by
a tile with angles $(\alpha,\beta, 2\alpha)$ with $\alpha$ not a rational multiple of $\pi$.
Then the tiling contains an essential segment with associated relation
$jc = pa + qb$.
\end{lemma} 

\begin{proof}
Suppose, for proof by contradiction,
that there is no such relation.
Then, by Lemma~\ref{lemma:extendlink4}, if $PQ$ is a link
in the graph $\Gamma_c$,  then there is a $c/b$ or $c/a$
segment emanating from $Q$.  Extend that segment to the 
maximal segment $QR$ supporting only tiles with $c$ edges
on $QR$.  Since there is no relation $jc = pa + qb$,
$R$ cannot be the vertex of a tile on the other side 
of $QR$.  Therefore, $QR$ is a link in $\Gamma_c$.

Therefore, the out-degree of every node $Q$ in $\Gamma_c$ is 
at least one.  The in-degree of $\Gamma_c$ is always
at most one.  Since the total out-degree is equal to the 
total in-degree, it follows that every node of $\Gamma_c$
has both in-degree and out-degree equal to 1.  Since
no link of $\Gamma_c$ can terminate on the boundary of $ABC$,
there can be no links of $\Gamma_c$ emanating from a vertex
on the boundary of $ABC$.

 We will show there cannot be two $\gamma$ angles at any vertex on 
the boundary of $ABC$.  Since $\gamma = 2\alpha$,
we have $3\alpha + \beta = \pi$.   If there were two $\gamma$
angles at vertex $V$ then we would have $2\gamma + p \alpha + q\beta = \pi$
for some integers $p,q$, so $(4+p)\alpha +q \beta = \pi$.
We have $q > 0$ since $\alpha$ is not a rational multiple of $\pi$.
Hence, $q \ge 1$ since $q$ is an integer. 
Subtracting
$3\alpha + \beta = \pi$ we have $(p+1)\alpha + (q-1)\beta = 0$. 
This is impossible, since $\alpha$, $\beta$, and $p+1$ are all positive,
and $q-1 \ge 0$.  That completes the proof that there cannot be two 
$\gamma$ angles at any boundary vertex.

  We will show there is at least one 
$c$ edge on $AC$.  For if not, every tile supported 
by $AC$ has its $\gamma$ angle at a vertex on $AC$.
As just proved, there cannot be two $\gamma$
angles at any one vertex. However, no $\gamma$ angles 
occur at $A$ or $C$, where there are only $\alpha$ 
angles. Then there is one more $\gamma$
angle on $AC$ than possible vertices to receive them,
so by the pigeonhole principle, some vertex on $AC$ has
two $\gamma$ angles, which is a contradiction.  Therefore, as claimed,
there is at least one $c$ edge on $AC$. 

The vertex angle at $B$ is $\pi-2\alpha = \alpha + \beta$,
so there are two tiles at $B$, none of which has a $\gamma$ angle.
Then, by the same argument as for $AC$, there is at least one $c$ edge on $AB$,
and at least one $c$ edge on $BC$.

Since there is a single tile at $A$
with its $\alpha$ angle at $A$, the $b$ edge of that 
tile lies on $AC$ or on $AB$.  Then there exists a
segment $PQ$ lying on $AB$ or on $AC$ supporting only 
tiles with $c$ edges on $PQ$, and with an $a$ or a  $b$ edge beyond $Q$.
Then by Lemma~\ref{lemma:extendlink4}, there is a $c/a$ segment
or a $c/b$ segment emanating from the boundary point $Q$, say
$QR$. Choose $R$ as far as possible from $Q$ such that $QR$
bounds only $c$ tiles on one side.  Then $R$ is not a vertex
of a tile on the other side of $QR$, since that would give
rise to a relation $jc = pa + qb$ with $j > 0$.  Hence
$QR$ is a link in $\Gamma_c$.  That is a contradiction,
since $Q$ is on the boundary of $ABC$. 
\end{proof}

\begin{theorem} \label{theorem:rationaltile}
Let the isosceles triangle $ABC$ with base angles $\alpha$  be $N$-tiled by
a tile with angles $(\alpha,\beta, 2\alpha)$ with $\alpha$ not a rational multiple of $\pi$.  Then the tile is rational; that is, the ratios of its sides are rational.
\end{theorem}

\begin{remark}   If the tile is rational, then after scaling
we can assume its sides are integers with no common factor.
\end{remark}

\begin{proof}[Proof of Theorem~\ref{theorem:rationaltile}] 
By Lemma~\ref{lemma:essentialsegments},
there is an essential segment witnessing a relation $ja = pb + qc$,
and another essential segment witnessing $Jb = Pa + Qc$. 
 If either $q = 0$ or $Q = 0$, then $a/b \in \mathbb{Q}$. 
 So we may assume $q > 0$ and $Q > 0$. Therefore, we have (solving for $c$ in both):
$$ c = \frac{ja - pb}{q} = \frac{Jb-Pa}{Q}. $$
Hence $(jQ+qP)a = (Jq+pQ)b$, and therefore
$$ \frac a b = \frac {Jq+pQ}{jQ+qP}. $$
Therefore, $a/b$ is rational.

By Lemma~\ref{lemma:crelation4}, there is an essential segment
witnessing a relation $jc = pa + qb$.  Then $c/a = p/j + qb/a$ is rational.
\end{proof}

\section{On the Number of Tiles Required when $\gamma = 2\alpha$}
We continue to consider tilings of isosceles $ABC$ with base
angles $\alpha$ and vertex angle $\alpha + \beta$.
In the previous two sections, we showed that the tile has to 
be rational, and that in that case, an $N$-tiling always exists,
for some $N$.
Next we will
 try to show that some values of $N$ are impossible.
We have three theorems along that line.  First, $N$ cannot be squarefree;
in particular, $N$ cannot be a prime number.  Second, given $N$, the
possible tiles and triangles $ABC$ form an explicit finite set.  Third,
$N$ has to be ``at least so big'', i.e., we have a lower bound on $N$,
namely $N \ge 45$.

\subsection{Characterization of the Tile}
\begin{lemma}\label{lemma:twoalphatiles}
Suppose $(a,b,c)$ are the integer sides of a triangle with
angles $(\alpha,\beta,2\alpha)$.  Then 
$$ c^2 = a^2 + ab.$$
\end{lemma}

\begin{remark} Rational triangles with $\gamma = 2\alpha$ correspond to 
solutions of this equation with $c < a+b$ and $b < a+c$ and $a < b+c$.
  Examples include $(4,5,6)$ and $(9,7,12)$. 
\end{remark}

\begin{proof}[Proof of Lemma~\ref{lemma:twoalphatiles}]
By the law of cosines, and since $\gamma = 2\alpha$,
\begin{eqnarray*}
c^2 &=& a^2 + b^2 - 2ab\cos \gamma \\
&=& a^2 + b^2 - 2ab \cos 2\alpha \\
&=& a^2 + b^2 - 2ab (2\cos^2 \alpha-1)\\
&=& a^2 + b^2 + 2ab - 4ab\cos^2 \alpha.
\end{eqnarray*}
By the law of sines, 
$$\frac {\sin \alpha} a = \frac {\sin \gamma} c = \frac{\sin 2 \alpha} c
= \frac {2\sin \alpha \cos \alpha} c,$$ so $ \cos \alpha = c/(2a)$.
Hence
\begin{eqnarray*}
c^2 &=& a^2 + b^2  + 2ab - bc^2/a \\
&=& (a+b)^2 -bc^2/a.
\end{eqnarray*}
Hence $c^2(1+b/a) = (a+b)^2$.  Multiplying both sides by $a/(a+b)$,
we have $c^2 = a(a+b) = a^2 + ab$.
\end{proof}

The following lemma gives a more nuanced characterization of 
$(a,b,c)$.  It was published in \cite{luthar}, but we give a
short proof here.%
\footnote{\,We are indebted to Gerry Myerson for pointing out this 
representation of $(a,b,c)$ to us on MathOverflow.} 

\begin{lemma}\label{lemma:luthar} 
Let $(a,b,c)$ be integers with no common factor, and suppose the 
triangle with sides $(a,b,c)$ has angles $(\alpha, \beta, 2\alpha)$.
Then $(a,b,c) =(k^2, m^2-k^2,mk)$ for some relatively prime integers $(k,m)$,
with $2k > m > k$.  

Conversely, let $(a,b,c)$ be a triple of integers 
$(a,b,c) =(k^2, m^2-k^2,mk)$ with $2k > m > k$ and $m$ and $k$ relatively
prime.  Then $(a,b,c)$ form a triangle, and 
it has angles $(\alpha, \beta, 2\alpha)$.
\end{lemma}

\begin{example} $(4,5,6)$ satisfies this equation 
with $k = 2$ and $m=3$.  Therefore, it is an example of a tile
satisfying $\gamma = 2\alpha$.  Although $(1,3,2)$ satisfies this 
equation with $k = 1$ and $m=2$, it does not correspond to a triangle. 
\end{example}

\begin{remark} Thus, $b$ and $c$ are relatively prime,
but $a$ and $c$ have a common factor $k$ (if $k \neq 1$).
Also, $c > a$ but not necessarily $c > b$, and $\gamma$ can be 
more or less than a right angle.
\end{remark}

\begin{proof}[Proof of Lemma~\ref{lemma:luthar}]
By Lemma~\ref{lemma:twoalphatiles}, we have
$$ c^2 = a^2 + ab = a(a+b). $$
We claim that $a$ and $a+b$ are relatively prime.  Suppose a prime $p$
divides both.  Then $p$ divides $a$ and $b$, so $p^2$ divides
$a^2 + ab = c^2$, so $p$ divides $c$.  That contradicts the hypothesis
that $a$, $b$, and $c$ have no common factor.  Hence $a$ and $a+b$ are
relatively prime, as claimed.  Since their product $c^2$ is a square,
each of them is a square: there are positive integers $k$ and $m$ with
$$ a = k^2, \qquad a + b = m^2, \qquad c = mk. $$
Then $b = m^2 - k^2$.  Since $k^2$ and $m^2$ are relatively prime, so
are $k$ and $m$.  Since $b > 0$, we have $m > k$.  Since the sides form
a triangle, $b < a + c$, that is,
$$ (m-k)(m+k) = m^2 - k^2 < k^2 + mk = k(k+m), $$
so $m - k < k$, that is, $2k > m$.  That completes the proof of the
first claim of the lemma.

Conversely, suppose $(a,b,c) =(k^2, m^2-k^2,mk)$,
and $(a,b,c)$ form
a triangle.
Then one can check that   
\begin{eqnarray*}
c^2 &=& a^2 + ab \\
&=& a (a+b). 
\end{eqnarray*}
We must show that this equation implies $\gamma = 2\alpha$.
Since $c^2-a^2 = ab$, we have by the law of sines,
$$\sin^2 \gamma - \sin^2 \alpha = \sin \alpha \sin \beta.$$
Factoring, we have 
$$(\sin \gamma - \sin \alpha)(\sin \gamma + \sin \alpha) = \sin \alpha \sin \beta.$$
Since $\gamma + \alpha = \pi-\beta$, this gives $\sin(\gamma - \alpha) = \sin \alpha$,
so $\gamma-\alpha = \alpha$, since $\alpha < \pi/2$.  Hence $\gamma = 2\alpha$, as desired.

Finally, we will show that $(a,b,c) =(k^2, m^2-k^2,mk)$ can be the 
lengths of the sides of a triangle.  That is, the sum of any two is greater
than the third.  The calculations follow.

The sum of the second and third sides is greater than the first:
$b + c > a$ becomes $m^2-k^2 + mk > k^2$,
or $m^2 + mk > 2k^2$, which follows from $m > k$.

The sum of the first and third sides is greater than the second:
$a+c > b$ is $k^2 + mk > m^2-k^2$, or $k(k+m) > (m+k)(m-k)$,
or $k > m-k$, which follows from $2k > m$.

The sum of the first and second sides is greater than the third:
$a+b > c$ is $k^2 + (m^2-k^2) > mk$, which follows from $m > k$.
\end{proof}

\begin{lemma}\label{lemma:edgerelations}
Let $(a,b,c)$ be integers with no common factor that are the sides of a
triangle with angles $(\alpha,\beta,2\alpha)$.  Then there are no
nonnegative integers $p,q$ such that any of the following holds:
\smallskip

(i) $c = pa + qb$,
\smallskip

(ii) $b = pa + qc$,
\smallskip

(iii) $a = pb + qc$.
\smallskip

\noindent
Indeed, in case (ii), there are no such integers $p,q$ even if negative
values are allowed.
\end{lemma}

\begin{proof}
By Lemma~\ref{lemma:luthar}, there are relatively prime integers $k,m$
with $m/2 < k < m$ such that
$$ a = k^2, \qquad b = m^2-k^2, \qquad c = mk. $$
We note that $k \ge 2$: if $k = 1$, then $m/2 < 1 < m$ gives
$1 < m < 2$, which is impossible for an integer $m$.
Note also that $m - k \ge 1$.
\smallskip

(i) Suppose $c = pa + qb$ with $p,q \ge 0$.  Then
$$ mk = pk^2 + q(m^2-k^2). $$
Reducing modulo $k$, we have $qm^2 \equiv 0 \pmod k$.  Since $k$ and $m$
are relatively prime, $k$ divides $q$.  If $q = 0$, then $mk = pk^2$, so
$m = pk$, so $k$ divides $m$; since $k$ and $m$ are relatively prime,
$k = 1$, contradicting $k \ge 2$.  Hence $q \neq 0$, and since $k$
divides $q$, we have $q \ge k$.  Then, since $p \ge 0$ and $m-k \ge 1$,
$$ mk \ \ge\ q(m-k)(m+k) \ \ge\ k(m-k)(m+k) \ \ge\ k(m+k) \ >\ km, $$
which is a contradiction.
\smallskip

(ii) Suppose $b = pa + qc$ with integers $p,q$ (of any sign).  Then
$$ m^2 - k^2 = pk^2 + qmk. $$
Reducing modulo $k$, we have $m^2 \equiv 0 \pmod k$.  Since $k$ and $m$
are relatively prime, $k = 1$, contradicting $k \ge 2$.
\smallskip

(iii) Suppose $a = pb + qc$ with $p,q \ge 0$.  Then
$$ k^2 = p(m^2-k^2) + qmk. $$
Reducing modulo $m$, we have $k^2 \equiv -pk^2 \pmod m$, that is,
$(p+1)k^2 \equiv 0 \pmod m$.  Since $k$ and $m$ are relatively prime,
$m$ divides $p+1$.  Since $p+1 > 0$, we have $p+1 \ge m$, so
$p \ge m-1$.  Then, since $q \ge 0$, $m - k \ge 1$, and $m \ge k+1$,
$$ k^2 \ \ge\ p(m-k)(m+k) \ \ge\ (m-1)(m+k) \ \ge\ k(2k+1) \ >\ k^2, $$
which is a contradiction.
\end{proof}

\begin{remark}
The hypothesis that the sides are commensurable cannot be dropped.
The golden triangle $(a,b,c) = (1,\varphi,\varphi)$, with
$\varphi = (1+\sqrt 5)/2$, has angles $(\pi/5, 2\pi/5, 2\pi/5)$, so
$\gamma = 2\alpha$, and it satisfies $c = 0\cdot a + 1\cdot b$ and
$b = 0\cdot a + 1\cdot c$.  The right isosceles triangle
$(1,1,\sqrt 2)$ has $\gamma = 2\alpha = \pi/2$ and satisfies
$a = 1\cdot b + 0\cdot c$.  In both examples $\alpha$ is a rational
multiple of $\pi$.  Also, nonnegativity is needed in (iii):
for the tile $(4,5,6)$ we have $a = 2b - c$.
\end{remark}

\subsection{Possible Shapes of the Tile}
In this section, we consider the possible shapes of 
a tile $(a,b,c)$ with $\gamma = 2\alpha$.  We begin
by observing that $a<c$ is the only obvious restriction
on the ordering of the edges.  As well as $a < b$,
we can have $b < a$, as in
 $(9,7,12)$, or  $a < c < b$, as in $(9,16,15)$.

 One way
of describing the shape of a triangle is by the ratios of 
its sides.  Here we give lower bounds on some of those ratios.
Actually, we use only the bound on $c/b$, and that only once,
 but it does seem necessary. We present all the bounds 
anyway, as they improve the reader's mental picture
of these (possible) tilings.

Empirically, these bounds are tight,
though we have not proved that they are best-possible.  There
are apparently no positive lower bounds on $b/a$ or $b/c$, although again
we have not proved that.  When $b/a$ is very small, an isosceles triangle
$ABC$ with $\gamma = 2\alpha$ will be very close to equilateral, and
by Laczkovich's method (see Figure~\ref{figure:bigisosceles2}), it can be
tiled with trillions of tiny needle-shaped tiles.  For example,
$(a,b,c) = (61504, 497, 61752)$ has $a/b > 123$, and we get
$N = 7227976999088825426993$, more than a billion trillion,
and the side and base of $ABC$ are $470043721566328$ and
$471939059153289$.  In the notation of Theorem~\ref{theorem:finitetiles} below,
$(k,m) = (248, 249)$ and $s = 1895337586961$.

This example also shows that the numbers of tiles on different sides of $ABC$
need not be related.  Here $\beta$ is tiny, $\alpha$ is close to $\pi/3$, and
$\gamma$ is close to $2\pi/3$.  This is not at all the situation in which $ABC$
itself is long and narrow.  $AB$ might be tiled by millions of tiles with
their tiny $b$ edges on $AB$, while $BC$ is tiled by relatively few tiles with
their long $a$ or $c$ edges on $BC$.  There is no obvious relation between the
number of tiles supported by one side and the number supported by another.
The proof of Theorem~\ref{theorem:finitetiles} below needs no such relation: it
uses only the tiles on one side, and only the crude bound that there are at
most $N$ of them.

\subsection{The Area Equation}
\begin{lemma} [Area equation] \label{lemma:areaequation}
Suppose isosceles triangle $ABC$, with base angles $\alpha$,
 is $N$-tiled by a tile 
with sides $(a,b,c)$ and angles $(\alpha,\beta,2\alpha)$.
  Let $X$ be the length of the equal 
sides $AB$ and $BC$.  Then $X^2 = Nab$.
\end{lemma}

\begin{proof}   Let $\gamma = 2\alpha$. 
The base angles of $ABC$ are $\alpha$, 
so $\pi = 2\alpha + \angle B = \gamma + \angle B$. Also, $\pi = \alpha + \beta + \gamma$,
so $\angle B = \alpha + \beta$.  
Twice the area of $ABC$ is given
by the magnitude of the cross product of $BA$ and $BC$, 
namely $X^2 \sin (\alpha + \beta)$.
Twice the area of the tile is given by $ab \sin \gamma$.
Since $\gamma = \pi- (\alpha + \beta)$, twice the area of 
the tile is also $ab \sin(\alpha+\beta)$.
The area of $ABC$ is $N$ times the area of the tile. Hence, 
$$ X^2 \sin(\alpha + \beta) = N ab \sin(\alpha + \beta). $$
Dividing both sides by  $\sin(\alpha+\beta)$, we have the area equation of the lemma.
\end{proof}

\begin{lemma} [Area equation, second form] \label{lemma:areaequation2}
Suppose isosceles triangle $ABC$, with base angles $\alpha$,
 is $N$-tiled by a tile 
with sides $(a,b,c)$ and angles $(\alpha,\beta,2\alpha)$.
  Let $X$ be the length of the equal 
sides $AB$ and $BC$, and let $Y$ be the length of $AC$.
  Then $XY = Nbc$, and $Y = (c/a) X$.
\end{lemma}

\begin{proof} Twice the area of $ABC$ is given by the 
magnitude of the cross product of $AB$ and $AC$, namely 
$XY \sin \alpha$.  Twice the area of the tile is $bc \sin \alpha$.
The area of $ABC$ is $N$ times the area of the tile. Hence, 
$XY \sin \alpha = N bc \sin\alpha$.  Dividing by $\sin \alpha$, we have
$$ XY = Nbc, $$
which proves the first claim of the lemma.  By 
Lemma~\ref{lemma:areaequation}, we have 
$X^2 = Nab$.  Dividing $XY = Nbc$ by $X^2 = Nab$ we have
$$ \frac Y X = \frac {N bc}{Nab} = \frac c a.$$
\end{proof}

\subsection{The Non-primality of $N$}  

In this section, we will show that $N$ cannot be a prime number.
What is more, $N$ cannot even be squarefree.  

\begin{lemma} \label{lemma:notprimehelper}
Let $ABC$ be isosceles with base angles $\alpha$, and 
$\alpha$ not a rational multiple of $\pi$.
Suppose $ABC$ is $N$-tiled by a tile with $\gamma = 2\alpha$
and integer sides $(a,b,c)$ with no common factor.  
 Then the squarefree part
of $N$ divides $b$ and the squarefree part of $b$ divides $N$.
If $N$ is squarefree, then 
$N$ divides $b$, and $b/N$ is a square, i.e., $b = N\ell^2$ 
for some integer $\ell$. 
\end{lemma}

\begin{example} In Figure~\ref{figure:bigisosceles2}, we
indicated a tiling with tile $(4,5,6)$ and $N = 6028020$.  Factoring that number,
we find $N = 2^2 \cdot 3^4 \cdot 5 \cdot 61^2$.  So the squarefree 
part of $N$ is $b = 5$, in accordance with this lemma.
   This provides a check on the
computation of the value of $N$, since it is not at all apparent 
from the construction of the tiling that $5$ has to be the 
squarefree part of $N$.
\end{example}

\begin{proof}[Proof of Lemma~\ref{lemma:notprimehelper}]
By Lemma~\ref{lemma:luthar}, there exist
relatively prime integers $m,k$ with $ 0 < m/2 < k < m$ such that
$$a = k^2, \quad b = m^2-k^2, \quad c = km.$$
Let $X$ be the length of the equal sides $AB$ and $BC$,
and $Y$ the length of $AC$.  Then by Lemmas~\ref{lemma:areaequation} and~\ref{lemma:areaequation2},
we have
\begin{align}
 X^2 &= Nab, \label{eq:area1}\\
\intertext{and}
 XY &= Nbc. \label{eq:area2}
\end{align}
Squaring Equation~(\ref{eq:area2}) and dividing by
Equation~(\ref{eq:area1}), and using $c^2/a = m^2$, we have
$$ Y^2 = Nb\,\frac{c^2}{a} = Nbm^2. $$
Hence $Nb = (Y/m)^2$ is a rational square, and since it is an integer,
it is an integer square.  Therefore $N$ and $b$ have the same squarefree
part, which proves the first two claims.  If $N$ is squarefree, then
$N$ is the squarefree part of $b$, so $b = N\ell^2$ for some
integer $\ell$.
\end{proof}

\begin{theorem} \label{theorem:notsquarefree}
Let $ABC$ be isosceles with base angles $\alpha$, and 
$\alpha$ not a rational multiple of $\pi$.
Suppose $ABC$ is $N$-tiled by a tile with angles $(\alpha, \beta,2\alpha)$
and sides $(a,b,c)$.
Then $N$ is not squarefree.   In particular, $N$ is not prime.
\end{theorem}

\begin{proof}
Suppose, for proof by contradiction, that $N$ is squarefree.
By Theorem~\ref{theorem:rationaltile}, the tile is rational,
so we can assume without loss of generality that $(a,b,c)$ are integers
with no common factor.  By Lemma~\ref{lemma:luthar}, there exist
relatively prime integers $m,k$ with $ 0 < m/2 < k < m$ such that 
$$a = k^2, b = m^2-k^2, \text{\ and \ } c = km.$$
Let $X$ be the length of the equal sides $AB$ and $BC$,
and $Y$ the length of $AC$. Then by Lemma~\ref{lemma:areaequation},
$X^2 = Nab$.   Since $b = N\ell^2$, by Lemma~\ref{lemma:notprimehelper},
$X^2 = N^2 a \ell^2$.   Since $a=k^2$ we have 
$$ X^2 = N^2k^2 \ell^2.$$
Taking the square root of both sides, we have
\begin{equation}
X = Nk\ell. \label{eq:9356}
\end{equation}
The tiling gives rise to a relation 
\begin{equation}
X = pa + qb + rc \label{eq:9360}
\end{equation}
Applying (\ref{eq:9356}), we have
\begin{equation}
Nk\ell = pk^2 + q(m-k)(m+k) + rkm. \label{eq:9361}
\end{equation}
 Taking 
that equation modulo $k$ we find
$$ 0 \equiv q m^2 \pmod k. $$
Since $k$ and $m$ are relatively prime, we can divide by $m^2$:
\begin{equation}
q \equiv 0 \pmod k. \label{eq:9370}
\end{equation}
Putting $N\ell^2 =  b  = (m-k)(m+k)$ into Equation~(\ref{eq:9361}), we have 
$$ k (m-k)(m+k)/\ell = pk^2 + q(m-k)(m+k) + rkm. $$
Hence
$$ 0 = pk^2 + (q-k/\ell)(m-k)(m+k) + rkm. $$
Since $m-k > 0$, we have $q-k/\ell \le 0$, and 
either $q-k/\ell < 0$ or $p=r=0$. 
We argue by cases:

\case{1} $q-k/\ell < 0$. Then 
$q < k/\ell \le k$. Then by Equation~(\ref{eq:9370}), we have $q=0$.
Therefore, no tile supported by $AB$ or $BC$ has its $b$
edge on $AB$ or $BC$, since a relation $X = pa + qb + rc$
would arise from each of $AB$ or $BC$ (although perhaps with 
different coefficients $(p,q,r)$).  However, at $B$ there
are two tiles, one with an $\alpha$ angle and one with a $\beta$
angle at $B$. Renaming $A$ and $C$ if necessary, we may assume
the tile with $\alpha$ at $B$ is supported by $AB$. 
Since each tile supported by $AB$ has its $a$ or $c$ edge
on $AB$, each of those tiles has a $\beta$ angle at one of 
its vertices on $AB$.  There are $\alpha$ angles at $A$ 
and at $B$.  Then by the pigeonhole principle, one of the vertices
on $AB$ is a vertex of two tiles with $\beta$ angles there.
That contradicts Lemma~\ref{lemma:twobeta}.
That disposes of Case~1.

\case{2} $p=r=0$.  Then every tile supported by the side $Z$ = 
$AB$ or $BC$ that gave rise to Equation~(\ref{eq:9360}) has its $b$ 
edge on $Z$.  Hence, every tile supported by $Z$
has a $\gamma$ angle on $Z$.  There are no $\gamma$ angles at 
$A$, $B$, or $C$, so by the pigeonhole principle, there must 
be a boundary vertex with two $\gamma$ angles. 
There cannot be two $\gamma$ 
angles at the same boundary vertex, since the only integer relations
between the angles are $\alpha + \beta + \gamma = \pi$ and 
$3\alpha + \beta = \pi$.  This contradiction shows that 
Case~2 is impossible.  
\end{proof}


\subsection{$N$ determines a finite set of tiles and $ABC$}
We wish to show that, given $N$, we can calculate a finite
set of triangles and a finite set of possible tiles
$(a,b,c)$,  such that if there is an $N$-tiling of some
isosceles $ABC$ with base angles $\alpha$ by some tile with
$\gamma = 2\alpha$ and $\alpha$ not a rational multiple of $\pi$,
then $ABC$ and the tile are in those finite sets.
For each $N$, let $T_N$ be the set of pairs of integers $(k,m)$ such that
\begin{eqnarray*}
&& 2 \le k < m < 2k, \qquad \gcd(k,m) = 1, \qquad m \le N + N^2, \\
&& N(m^2-k^2) \mbox{ is a perfect square.}
\end{eqnarray*}
$T_N$ is finite, and it can be computed by checking each pair $(k,m)$ with
$m \le N+N^2$.

\begin{theorem} \label{theorem:finitetiles}
Let isosceles triangle $ABC$ with base angles $\alpha$
be $N$-tiled by a
triangle with angles $(\alpha, \beta, 2\alpha)$,
and $\alpha$ not a rational multiple of $\pi$.
Then there is a pair $(k,m) \in T_N$ such that, after scaling,
the tile has sides
$$ a = k^2, \quad b = m^2-k^2, \quad c = km $$
(opposite $\alpha$, $\beta$, and $2\alpha$ respectively), and
$ABC$ has sides
$$ \vert AB \vert = \vert BC \vert = ks, \qquad \vert AC \vert = ms, $$
where $s$ is the positive integer with $s^2 = N(m^2-k^2)$.
\end{theorem}

\begin{remark}
Thus, given $N$, the pair $(k,m)$ determines both the tile and $ABC$, and there are
only finitely many pairs.  If $T_N$ is empty, there is no such $N$-tiling.
The bound $m \le N + N^2$ is crude, but we only need {\em some} explicit bound.
For example, $T_5 = \{(2,3), (11,21), (19,21), (22,23), (22,27)\}$, and $T_{45}$ has
330 elements.  Of course, there are no actual 5-tilings, as will be shown below;
these are necessary conditions only.

The point of the theorem is that, given $N$, the possible tiles and triangles
form an explicit finite list that can be computed from $N$.  Whether one of
them can actually be $N$-tiled is a separate question.
\end{remark}

\begin{proof}[Proof of Theorem~\ref{theorem:finitetiles}]
By Theorem~\ref{theorem:rationaltile}, the tile is rational,
so we can scale it so that its sides $(a,b,c)$ are integers
with no common factor.  By Lemma~\ref{lemma:luthar}, there exist
relatively prime integers $k, m$ with $k < m < 2k$ such that
$$a = k^2, \quad b = m^2-k^2, \quad c = km.$$
Since $k < m < 2k$, we have $k \ge 2$.

Let $X = \vert AB \vert$.  The tiles supported by $AB$ give a representation
\begin{equation}
X = pa + qb + rc \ =\ pk^2 + q(m^2-k^2) + rkm, \label{eq:finiterep}
\end{equation}
where $p$, $q$, and $r$ are the numbers of tiles with their $a$, $b$, and $c$ edges on $AB$.
A triangle cannot have two edges on the same line, so each tile has at most one edge on $AB$.
Hence $p + q + r \le N$.  Since $a$, $b$, and $c$ are integers, so is $X$.

By Lemma~\ref{lemma:areaequation}, $X^2 = Nab = Nk^2(m^2-k^2)$.  Hence
$(X/k)^2 = N(m^2-k^2)$ is an integer, so the rational number $s = X/k$ is an
integer, since its square is an integer.  Then $s > 0$, $s^2 = N(m^2-k^2)$,
and $\vert AB \vert = \vert BC \vert = X = ks$.
By Lemma~\ref{lemma:areaequation2},
$$ \vert AC \vert = \frac c a X = \frac {km}{k^2}\, ks = ms. $$

It remains to show $m \le N + N^2$.  Reducing (\ref{eq:finiterep}) mod $m$,
$$ X \equiv pk^2 - qk^2 = (p-q)k^2 \pmod m, $$
so $X^2 \equiv (p-q)^2 k^4 \pmod m$.  On the other hand,
$X^2 = Nk^2(m^2-k^2) \equiv -Nk^4 \pmod m$.  Hence
$$ (N + (p-q)^2)\, k^4 \equiv 0 \pmod m. $$
Since $k$ and $m$ are relatively prime, $m$ divides $N + (p-q)^2$.
Since $N + (p-q)^2 > 0$, we have $m \le N + (p-q)^2$.  Finally
$\vert p - q \vert \le p + q \le N$, so $m \le N + N^2$.
\end{proof}

\subsection{What is the Least $N$ Permitting a Tiling?}
The smallest explicitly-known such tiling has $N = 720$,
with tile $(4,5,6)$,
as shown in Figure~\ref{figure:bryce5}. 
We will show below  that $N \ge 45$.  There
is a  swath of ignorance between 45 and 720.

\begin{lemma} \label{lemma:qandrnotzero}
Let isosceles triangle $ABC$ with base $AC$ and vertex $B$
have base angles $\alpha$,
not a rational multiple of $\pi$. Let 
$(a,b,c)$ be integers with no common factor forming
a triangle with angles $(\alpha,\beta,2\alpha)$.  Let
$X$  be the length of side $AB$, where the tile on $AB$ at $B$
has its $\alpha$ angle at $B$.
Suppose there is an $N$-tiling of $ABC$ by $(a,b,c)$
whose boundary tiling on $AB$ corresponds to $X = pa + qb + rc$.
Then $q\neq 0$ and $r \neq 0$.
\end{lemma}

\begin{proof} 
First suppose $q = 0$.  That means
there are no $b$ edges on $AB$, so every tile supported by $AB$ has a 
$\beta$ angle on $AB$.  Since the top and bottom tiles have their 
$\alpha$ angles at $B$ and $A$,  and no vertex has two $\beta$ angles,
 this violates the pigeonhole principle.  Now suppose $r = 0$.
That 
means there are no $c$ edges on $AB$, so every tile supported by $AB$
has a $\gamma$ angle on $AB$.  Since the top and bottom tiles have their 
$\alpha$ angles at $B$ and $A$,  and no vertex has two $\gamma$ angles,
 this violates the pigeonhole principle. 
 \end{proof}

\begin{lemma} \label{lemma:no28}
Let isosceles triangle $ABC$ with base $AC$ and vertex $B$
have base angles $\alpha$, not a rational multiple of $\pi$.  Let
$(a,b,c)$ be integers with no common factor forming a triangle with
angles $(\alpha,\beta,2\alpha)$, and suppose $ABC$ is $N$-tiled by
$(a,b,c)$, with the tile at $B$ supported by $AB$ having its $\alpha$
angle at $B$.  Then $AB$ does not support exactly one tile with its
$c$ edge on $AB$.
\end{lemma}

\begin{proof}
By Lemma~\ref{lemma:luthar}, $a$ and $b$ are relatively prime, $c > a$,
and $b \neq 1$; hence $a$ is not an integer multiple of $b$.  Suppose,
for proof by contradiction, that exactly one tile supported by $AB$,
say $T$, has its $c$ edge on $AB$.  Every other tile supported by $AB$
has its $a$ or $b$ edge on $AB$, and hence has its $\gamma$ angle at
an endpoint of that edge.
\smallskip

{\em Step 1: the tile $T_1$ at $A$ has its $b$ edge on $AB$ and its
$c$ edge on $AC$.}  $T_1$ has its $\alpha$ angle at $A$, so its edges
at $A$ are its $b$ and $c$ edges.  Suppose $T_1$ has its $c$ edge on
$AB$, so $T_1 = T$.  Let $P$ be the other end of that edge.  The other
$n-1$ tiles supported by $AB$ have their $\gamma$ angles at the $n-1$
vertices in the interior of $AB$, since there are no $\gamma$ angles
at $A$ or $B$, and by Lemma~\ref{lemma:twobeta} no vertex on $AB$ has
two $\gamma$ angles.  Hence the tile supported by $AB$ next to $T_1$
has its $\gamma$ angle at $P$.  $T_1$ has its $\beta$ angle at $P$, so
the vertex $P$ has angles $(\alpha,\beta,\gamma)$, and a third tile
$T_2$, lying between these two, has its $\alpha$ angle at $P$.  The
edges of $T_2$ at $P$ are its $b$ and $c$ edges, and one of them lies
along the edge of $T_1$ at $P$ that is not on $AB$, which is the $a$
edge of $T_1$.  That $a$ edge has both endpoints on the boundary of
$ABC$, so it is covered on the side of $T_2$ by tile edges of total
length $a$.  Since $c > a$, $T_2$ has its $b$ edge there, and the rest
of that $a$ edge, of length $a-b < a$, is covered by edges shorter
than $a$, which are $b$ edges.  Then $a$ is a multiple of $b$, a
contradiction.  This proves Step~1.  In particular $T \neq T_1$.
\smallskip

{\em Step 2: each vertex in the interior of $AB$ has exactly three
tiles, with angles $\alpha$, $\beta$, $\gamma$, and exactly one of the
two tiles supported by $AB$ there has its $\gamma$ angle there.}
Let $n$ be the number of tiles supported by $AB$.  The $n-1$ tiles
other than $T$ have their $\gamma$ angles at the $n-1$ vertices in the
interior of $AB$, at most one at each vertex by
Lemma~\ref{lemma:twobeta}; so each such vertex has exactly one of
them.  A vertex on a side of $ABC$ with a $\gamma$ angle has angles
$(\alpha,\beta,\gamma)$, since $\pi = 3\alpha + \beta$ is the only
other way to write $\pi$ as a sum of tile angles.
\smallskip

{\em Step 3.}  Let $PQ$ be the $c$ edge of $T$ on $AB$, where $T$ has
its $\alpha$ angle at $P$ and its $\beta$ angle at $Q$, and let $R$ be
the third vertex of $T$, where $T$ has its $\gamma$ angle.  Then
$\vert PR \vert = b$ and $\vert QR \vert = a$.  Since $T \neq T_1$, neither
$P$ nor $Q$ is $A$.  Since the tile at $B$ supported by $AB$ has its
$\alpha$ angle at $B$, $Q \neq B$.  So $Q$ lies in the interior of $AB$.
By Step~2, and since $T$ has no $\gamma$ angle at $Q$, the other tile
supported by $AB$ at $Q$ has its $\gamma$ angle there, and a third tile
$G$, lying between the two, has its $\alpha$ angle at $Q$.  The edges
of $G$ at $Q$ are its $b$ and $c$ edges, and one of them lies along
$QR$.
\smallskip

{\em Step 4: $R$ lies in the interior of an edge of a tile on the side
of the line $QR$ opposite to $T$.}  The segment $QR$ is covered, on the
side opposite to $T$, by tile edges along the line $QR$, the first of
which is the edge of $G$ at $Q$.  If some edge ended exactly at $R$,
then $a = \vert QR \vert$ would be the sum of the lengths of the edges
from $Q$ to $R$.  The first of these is the $b$ or $c$ edge of $G$;
since $c > a$, it is $b$, and $b < a$, and the remaining edges have
total length $a - b < a$, so they are $b$ edges; then $a$ is a multiple
of $b$, a contradiction.  Hence some edge along the line $QR$ contains
$R$ in its interior.
\smallskip

{\em Step 5.}  Let $T^*$ be the tile of Step~4.  Near $R$, the side of
the line $QR$ opposite to $T$ is covered by $T^*$.  The line $PR$
crosses the line $QR$ at $R$, so beyond $R$ it enters the interior of
$T^*$.  Hence no tile edge along the line $PR$ extends beyond $R$, and
the tile edges along $PR$ on the side opposite to $T$ have total length
exactly $\vert PR \vert = b$.  Let $T_2$ be the first of these tiles, the
one with a vertex at $P$.  If $P$ lies in the interior of $AB$, then by
Step~2 the vertex $P$ has three tiles, namely $T$ with its $\alpha$
angle, the other tile supported by $AB$ with its $\gamma$ angle, and
$T_2$ between them, which therefore has its $\beta$ angle at $P$.  If
$P = B$, the two tiles at $B$ have angles $\alpha$ and $\beta$; $T$
has the $\alpha$, so again $T_2$ has its $\beta$ angle at $P$.  In
either case, the edges of $T_2$ at $P$ are its $a$ and $c$ edges, so
the first edge along $PR$ is not a $b$ edge.  Then none of the edges
along $PR$ is a $b$ edge, since their total length is $b$.  Hence $b$
is a sum of $a$ edges and $c$ edges, contradicting
Lemma~\ref{lemma:edgerelations}.
\end{proof}

\begin{theorem} \label{theorem:N45}
Let isosceles triangle $ABC$ with base angles $\alpha$
 be $N$-tiled by a
triangle with angles $(\alpha, \beta, 2\alpha)$,
and $\alpha$ not a rational multiple of $\pi$.
Then  $N \ge 45$.
\end{theorem}

\begin{proof}[Proof of Theorem~\ref{theorem:N45}]
By Theorem~\ref{theorem:rationaltile}, the tile is rational,
so we can assume without loss of generality that its sides $(a,b,c)$ are integers
with no common factor.  By Lemma~\ref{lemma:luthar}, there exist
relatively prime integers $k, m$ with $k < m < 2k$ such that
$$a = k^2, \quad b = m^2-k^2, \quad c = km.$$
Renaming $A$ and $C$ if necessary, we may assume that the tile at $B$
supported by $AB$ has its $\alpha$ angle at $B$.  Let $X$ be the length
of $AB$.  The tiles supported by $AB$ give a representation
\begin{equation}
X = pa + qb + rc \ =\ pk^2 + qb + rkm, \label{eq:N45rep}
\end{equation}
and by Lemma~\ref{lemma:qandrnotzero}, $q \ge 1$ and $r \ge 1$.

By Lemma~\ref{lemma:areaequation}, $X^2 = Nab = Nk^2b$.  Hence
$(X/k)^2 = Nb$ is an integer, so the rational number $t = X/k$ is an
integer, and
\begin{equation}
X = kt, \qquad t^2 = Nb. \label{eq:N45t}
\end{equation}
Since $k$ divides $X$, Equation~(\ref{eq:N45rep}) shows that $k$ divides $qb$.
Now $b = m^2 - k^2 \equiv m^2 \pmod k$, and $k$ and $m$ are relatively
prime, so $k$ and $b$ are relatively prime.  Hence $k$ divides $q$; write
$q = kq'$, with $q' \ge 1$ since $q \ge 1$.  Dividing
Equation~(\ref{eq:N45rep}) by $k$, we find
\begin{equation}
t = pk + q'b + rm. \label{eq:N45tsum}
\end{equation}

Now suppose, for proof by contradiction, that $N \le 44$.
We claim that $r = 1$.  Suppose instead that $r \ge 2$.  Then by
Equation~(\ref{eq:N45tsum}), $t \ge b + 2m$, so by Equation~(\ref{eq:N45t}),
\begin{equation}
(b+2m)^2 \ \le\ t^2 \ =\ Nb \ \le\ 44b. \label{eq:N45bound}
\end{equation}
Since $(b+2m)^2 \ge b^2 + 4bm$, this gives $b + 4m \le 44$.  Since $k \le m-1$,
we have $b = m^2 - k^2 \ge 2m - 1$, so $6m - 1 \le 44$ and $m \le 7$.  The pairs
of relatively prime integers with $k < m < 2k$ and $m \le 7$ are
$$(k,m) = (2,3),\ (3,4),\ (3,5),\ (4,5),\ (4,7),\ (5,6),\ (5,7),\ (6,7),$$
with $b = 5, 7, 16, 9, 33, 11, 24, 13$ respectively.  For the last four,
(\ref{eq:N45bound}) fails: $(b+2m)^2 = 2209, 529, 1444, 729$, while
$44b = 1452, 484, 1056, 572$.  For the first four, $t$ satisfies
$b + 2m \le t \le \sqrt{44b}$ and $b$ divides $t^2$.  We will show that
in each case, there is no such $t$.
\smallskip

$(k,m) = (2,3)$, $b = 5$: $11 \le t \le 14$, but $5$ divides $t^2$ only if $5$ divides $t$.

$(k,m) = (3,4)$, $b = 7$: $15 \le t \le 17$, but $7$ divides $t^2$ only if $7$ divides $t$.

$(k,m) = (3,5)$, $b = 16$: $t = 26$, but $16$ does not divide $26^2 = 676$.

$(k,m) = (4,5)$, $b = 9$: $t = 19$, but $9$ does not divide $19^2 = 361$.
\smallskip

\noindent
As claimed, in every case there is no such $t$, which proves the claim $r = 1$.
Then $AB$ supports exactly one tile with its $c$ edge on $AB$, which is
impossible by Lemma~\ref{lemma:no28}.  Hence $N \ge 45$.
\end{proof}

\begin{remark}
The technique does not extend to $N=45$: for the tile $(4,5,6)$ ($k=2$, $m=3$)
we have $t^2 = 45 \cdot 5$, so $t = 15 \ge b + 2m = 11$, and indeed
$X = 30 = 2a + 2b + 2c$ has two $c$ edges.
\end{remark}

\section{Tilings of an Isosceles Triangle by a Tile $(\alpha,\beta,2\pi/3)$}

In this section we take up the tilings of 
an isosceles triangle (and not equilateral) $ABC$ with base angles $\alpha$ or $\beta$, by a tiling
with angles $(\alpha,\beta,2\pi/3)$, where $\alpha$ is 
not a rational multiple of $\pi$.   Let $\gamma = 2\pi/3$. 
We could insist that the base angles are called $\alpha$, but 
then we may have to speak of tilings by $(5,3,7)$ instead of $(3,5,7)$;
so it is convenient to allow the base angles to be $\beta$ sometimes.
If the base angles are $\alpha$, then the vertex angle is $\pi-2\alpha = \alpha + 3 \beta$.

\begin{figure}[ht]
\begin{center}
\includegraphics[width=\textwidth]{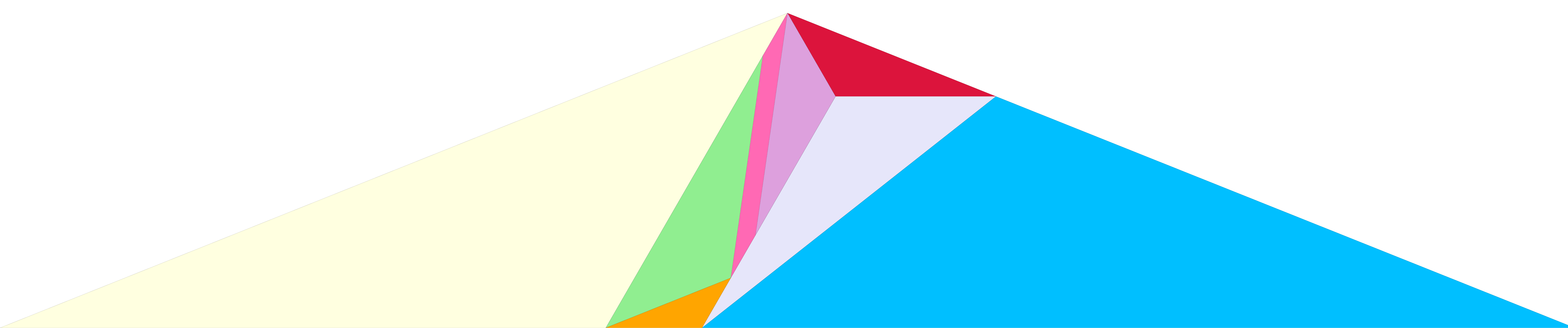}
\end{center}
\caption{With tile $(3,5,7)$, $N$ would be $1878500$, which is too large to draw.}
\label{figure:isosceles2pi}
\end{figure}

 Laczkovich proved \cite[Theorem~2.5]{laczkovich1995}
that there exist  tiles that
 can be used to tile {\em some} such $ABC$, but $N$ constructed
 by his method can be large.  He proceeds by
 first constructing a dissection of $ABC$ into {\em similar} 
rational triangles and parallelograms. 
 Figure~\ref{figure:isosceles2pi} shows such a preliminary dissection. Then to get a tiling by congruent triangles, we have to 
choose a very small tile such that if each of the visible triangles is 
tiled quadratically, then every shared edge is an integer multiple of the 
tile edges.  For example if the red triangle will get $p^2$ tiles and 
the light blue triangle will get $q^2$ tiles, then we must satisfy $pb = qa$,
in this case $5p = 3q$.  There will be another such equation on every 
shared boundary.  To clear all the denominators we will have to use a 
large number of tiles.

\begin{figure}[ht]
\begin{center}
\includegraphics[width=\textwidth]{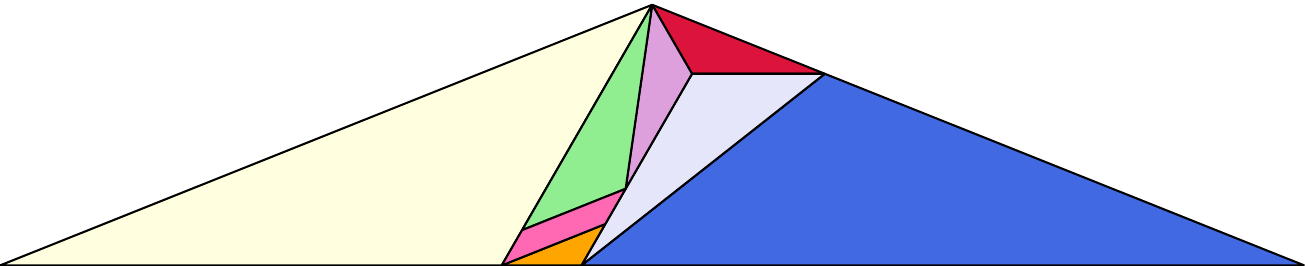}
\end{center}
\caption{With this configuration, $N$ is $3681860$, even larger than in Figure~\ref{figure:isosceles2pi}.}
\label{figure:isosceles2pi2}
\end{figure} 

Figure~\ref{figure:isosceles2pi2} shows that 
a slightly different arrangement of the similar triangles and parallelograms
can make a difference in the resulting number of tiles.  In Figure~\ref{figure:isosceles2pi},
the parallelogram has an  $\alpha$ angle and is tiled using $b$ and $c$ edges on the boundary.
In Figure~\ref{figure:isosceles2pi2}, the parallelogram has a $\beta$ angle and is 
tiled using $a$ and $c$ edges on the boundary.
 That makes the equations at 
the boundary different, even though the other boundary conditions 
are the same and the areas of the two parallelograms are equal.

Nevertheless, these tilings are too large to draw.
In 2024, Bryce Herdt discovered a 2673-tiling, which is 
exhibited in Section \ref{section:herdt}. This is presently the smallest known tiling of an isosceles triangle by 
a tile with $\gamma = 2\pi/3$.

\subsection{The Tile is Rational}
Suppose given an $N$-tiling of some triangle $ABC$ by a tile
with angles $(\alpha,\beta,\gamma)$ and sides $(a,b,c)$.  In this 
section we will prove that the tile has to be rational, i.e., the 
ratios of the sides are all rational, so after a suitable scaling,
they will be integers.  The proof uses the graphs $\Gamma_c$ 
introduced by Laczkovich and described above; several preparatory 
lemmas will be developed first.  

\begin{definition} \label{definition:c/a}
A {\em $c/a$  segment} is a left-terminated
 interior segment $PQ$ of the tiling
supporting two tiles on opposite sides of $PQ$,
each with a vertex at $P$, one with its $c$ edge on $PQ$ and 
one with its $a$ edge on $PQ$.  The segment is said 
to {\em emanate from $P$}. 

A {\em $c/b$  segment} is a left-terminated
 interior segment $PQ$ of the tiling
supporting two tiles on opposite sides of $PQ$,
each with a vertex at $P$, one with its $c$ edge on $PQ$ and 
one with its $b$ edge on $PQ$.    

An {\em $a/b$  segment} is a left-terminated
 interior segment $PQ$ of the tiling
supporting two tiles on opposite sides of $PQ$,
each with a vertex at $P$, one with its $a$ edge on $PQ$ and 
one with its $b$ edge on $PQ$.   
\end{definition}

\begin{remark} The point $Q$ serves only to indicate
the direction of the segment; it can be any point on that ray.
\end{remark}

\begin{lemma} \label{lemma:extendlink} 
Let $ABC$ be an isosceles triangle with base angles $\alpha$,
tiled by tile $(\alpha,\beta,\gamma)$ with $\gamma = 2\pi/3$.
Let $PQR$ be an (internal or boundary) segment of the tiling with only $c$ 
edges on one side of $PQ$, such that the tile on that side supported
by $QR$ with a vertex at $Q$ has an $a$ or $b$ edge on $PQ$.
Then there is a $c/a$ or $c/b$ segment emanating from $Q$.
\end{lemma}

\begin{proof} We only consider tiles on the one side 
of $PQR$ mentioned in the lemma.  There are two cases:
either there are three tiles with a vertex at $Q$, or there
are six.

\case{1} There are three tiles at $Q$. Then one of them has its $\gamma$ angle
at $Q$, and hence no $c$ edge at $Q$. The other two have a 
$c$ edge ending at $Q$.  One of those lies on $PQ$.  The
other does {\em not} lie on $QR$, by hypothesis.  Since there
is no other $c$ edge ending at $Q$, that third $c$ edge forms
either a $c/b$ segment or a $c/a$ segment.

\case{2} There are six tiles at $Q$.  Then all six angles are $\alpha$ or
$\beta$.  Each of the six tiles has a $c$ edge ending at $Q$.
One lies on $PQ$ and five lie on interior segments.  Since five is odd,
one of those $c$ edges is not paired with another $c$ edge, and hence
constitutes a $c/a$ segment or a $c/b$ segment.  
\end{proof}

\begin{lemma} \label{lemma:crelation}
Suppose isosceles triangle $ABC$ with base angles $\alpha$ 
is $N$-tiled by $(\alpha,\beta,2\pi/3)$, with $\alpha$ not a 
rational multiple of $\pi$.  Then there is a relation
$$ jc = pa + qb$$
with nonnegative integers $p,q,j$ and $j > 0$. 
\end{lemma}

\begin{proof}  Suppose, for proof by contradiction,
that there is no such relation.
Then, by Lemma~\ref{lemma:extendlink}, if $PQ$ is a link
in the graph $\Gamma_c$,  then there is a $c/b$ or $c/a$
segment emanating from $Q$.  Extend that segment to the 
maximal segment $QR$ supporting only tiles with $c$ edges
on $QR$.  Since there is no relation $jc = pa + qb$,
$R$ cannot be the vertex of a tile on the other side 
of $QR$.  Therefore, $QR$ is a link in $\Gamma_c$.

Therefore, the out-degree of every node $Q$ in $\Gamma_c$ is 
at least one.  The in-degree of $\Gamma_c$ is always
at most one.  Since the total out-degree is equal to the 
total in-degree, it follows that every node of $\Gamma_c$
has both in-degree and out-degree equal to 1.  Since
no link of $\Gamma_c$ can terminate on the boundary of $ABC$,
there can be no links of $\Gamma_c$ emanating from a vertex
on the boundary of $ABC$.

 We will show there is at least one 
$c$ edge on $AC$.  For if not, every tile supported 
by $AC$ has its $\gamma$ angle at a vertex on $AC$.
Since $\gamma > \pi/2$, there cannot be two $\gamma$
angles at any one vertex. However, no $\gamma$ angles 
occur at $A$ or $C$, where there are only $\alpha$ 
angles. Then there is one more $\gamma$
angle on $AC$ than possible vertices to receive them,
so by the pigeonhole principle, some vertex on $AC$ has
two $\gamma$ angles, which is a contradiction.  Therefore, as claimed,
there is at least one $c$ edge on $AC$. 

The vertex angle at $B$ is $\pi-2\alpha = \alpha + 3 \beta$,
so there are four tiles at $B$, none of which has a $\gamma$ angle.
Then, by the same argument as for $AC$, there is at least one $c$ edge on $AB$,
and at least one $c$ edge on $AC$.

Since there is a single tile at $A$
with its $\alpha$ angle at $A$, the $b$ edge of that 
tile lies on $AC$ or on $AB$.  Then there exists a
segment $PQ$ lying on $AB$ or on $AC$ supporting only 
tiles with $c$ edges on $PQ$, and with an $a$ or a  $b$ edge beyond $Q$.
Then by Lemma~\ref{lemma:extendlink}, there is a $c/a$ segment
or a $c/b$ segment emanating from the boundary point $Q$, say
$QR$. Choose $R$ as far as possible from $Q$ such that $QR$
bounds only $c$ tiles on one side.  Then $R$ is not a vertex
of a tile on the other side of $QR$, since that would give
rise to a relation $jc = pa + qb$ with $j > 0$.  Hence
$QR$ is a link in $\Gamma_c$.  That is a contradiction,
since $Q$ is on the boundary of $ABC$. 
\end{proof}

\begin{lemma} \label{lemma:abrelation}
Suppose isosceles triangle $ABC$ with base angles $\alpha$
is $N$-tiled by $(\alpha,\beta,2\pi/3)$, with $\alpha$ not a
rational multiple of $\pi$.  Then there are relations
$$ ja = pb + qc \qquad\mbox{and}\qquad Jb = Pa + Qc $$
with nonnegative integers $j,p,q,J,P,Q$ and $j > 0$, $J > 0$.
\end{lemma}

\begin{proof}   For tilings by a triangle with $\gamma = 2\pi/3$,
we classify the ``types'' of vertices that can occur in a tiling
according to the number of $\alpha,\beta$, and $\gamma$ angles occurring
at the vertex.
 A {\em center} is a vertex
of the tiling where three tiles meet, each having its
$\gamma$ angle at that vertex.  A {\em star} is a vertex
$P$ where six tiles lying on one side of a line through $P$
have three $\alpha$ and three $\beta$ angles at $P$. 
A star can occur on the boundary of $ABC$ or in the interior.
A {\em double star} is a vertex where twelve tiles meet,
six with $\alpha$ angles and six with $\beta$ angles.  
The three vertices of $ABC$ together have six angles, three
$\beta$ and three $\alpha$, the same count as a star.
An {\em anomaly} is a vertex where four $\alpha$, four $\beta$, 
and one $\gamma$ angle meet.  A {\em simple vertex} has three 
tiles, one each of $\alpha, \beta, \gamma$.  See Table~\ref{table:vertextypes}.

\begin{table}[htp]
\begin{center}
\begin{tabular}{ c | c }
Name  & $(\alpha,\beta,\gamma)$  \\
\hline
simple    & $(1,1,1)$ \\
center    & $(0,0,3)$  \\
star      & $(3,3,0)$ \\
double star  & $(6,6,0)$ \\
anomaly      & $(4,4,1)$  \\
double simple & $(2,2,2)$
\end{tabular}
\end{center}
\caption{Vertex Types.}
\label{table:vertextypes}
\end{table}%

 We will show that every vertex is one of the types shown in Table~\ref{table:vertextypes}.
 To prove that, let $X$ be a vertex. The sum of the
angles of tiles at $X$ is either $2\pi$ or $\pi$.  

\case{1} The angle sum is $\pi$.  If there is a $\gamma$ angle at $X$ then the other two are
an $\alpha$ and a $\beta$, and $X$ is a simple vertex.  If there is no 
$\gamma$ angle, then since $\pi = 3 \alpha + 3 \beta$, $X$ is a star.

\case{2} The angle sum is $2 \pi$.  Let $w$ be the number of $\gamma$
angles at $X$.  Since $\gamma = 2\pi/3$, the remaining angles have sum
$(3-w)\,2\pi/3 = 2(3-w)(\alpha+\beta)$.  Since $\alpha$ is not a rational
multiple of $\pi$ and $\beta = \pi/3 - \alpha$, a sum of $u$ angles
$\alpha$ and $v$ angles $\beta$ equals $r(\alpha+\beta)$ only if
$u = v = r$.  Hence the remaining angles are $2(3-w)$ angles $\alpha$
and $2(3-w)$ angles $\beta$.  For $w = 3, 2, 1, 0$ this gives a center,
a double simple vertex $(2,2,2)$, an anomaly $(4,4,1)$, and a double
star $(6,6,0)$, respectively.  That completes the proof that every
vertex is one of the types in Table~\ref{table:vertextypes}.

Suppose, for proof by contradiction, that there is no relation as
alleged in the statement of the lemma.  Let $\C$ be the number of
centers, and let ${\mathcal S}$ be the number of stars plus the number
of anomalies plus twice the number of double stars.  For each vertex,
compute the number of $\alpha$ angles plus the number of $\beta$ angles
minus twice the number of $\gamma$ angles at that vertex.  By
Table~\ref{table:vertextypes}, this is $0$ at simple and double simple
vertices, $-6$ at a center, $6$ at a star and at an anomaly, and $12$
at a double star.  The vertices $A$, $B$, $C$ together carry three
$\alpha$ and three $\beta$ angles, contributing $6$.  Summing over all
vertices we get $6\,{\mathcal S} - 6\,\C + 6$.  On the other hand, the
total numbers of $\alpha$, $\beta$, and $\gamma$ angles are each $N$,
so the sum is $N + N - 2N = 0$.  Hence
$$ {\mathcal S} = \C - 1. $$
For example, in Figure~\ref{figure:isosceles2pi2}, we see one
center and no stars.

Now we consider the graph $\Gamma_a$.  Every center has an outgoing
link.  Indeed, at a center $P$ there are three tiles, each with its
$\gamma$ angle at $P$, so each has its $a$ edge and its $b$ edge ending
at $P$, and no other edges end at $P$.  The three rays from $P$ each
carry two of these six edges.  If no ray carried both an $a$ edge and a
$b$ edge, the number of $a$ edges would be even; since it is $3$, some
ray carries an $a$ edge and a $b$ edge.  Let $Q$ be the farthest point
from $P$ along that ray such that $PQ$ supports only $a$ edges on one
side.  If $Q$ were a vertex of a tile on the other side, we would have
a relation $ja = pb + qc$ with $p > 0$, since there is a $b$ edge on
the other side of $PQ$; by hypothesis there is no such relation.
Hence $PQ$ is a link in $\Gamma_a$.

If $PQ$ is a link, then $Q$ lies on the interior of an edge of a tile
on the other side of $PQ$, and $Q$ is not on the boundary of $ABC$.
Hence: at most one link ends at any vertex; no link ends at $A$, $B$,
or $C$; and no link ends at a vertex with angle sum $2\pi$, that is,
at a center, a double simple vertex, an anomaly, or a double star.

Let $PX$ be a link ending at a simple vertex $X$.  We claim there is
an outgoing link at $X$.  At $X$ there are three tiles, with angles
$\alpha$, $\beta$, $\gamma$ at $X$.  The tiles with their $\beta$ and
$\gamma$ angles at $X$ each have an $a$ edge ending at $X$; the tile
with its $\alpha$ angle does not.  One of these two $a$ edges lies on
$PX$.  The other lies on a ray from $X$ whose other edge at $X$
belongs to an adjacent tile and is not an $a$ edge, since no third
$a$ edge ends at $X$.  Let $XQ$ be the longest segment on that ray
supporting only $a$ edges on the side of that $a$ edge.  If $Q$ were a
vertex of a tile on the other side, we would have a relation
$ja = pb + qc$ with $j > 0$, contrary to hypothesis.  Hence $XQ$ is a
link of $\Gamma_a$.

We now compare out-degree and in-degree vertex by vertex.  At a center
the in-degree is $0$ and the out-degree is at least $1$.  At a vertex
with angle sum $2\pi$ that is not a center, and at $A$, $B$, and $C$,
the in-degree is $0$, so the out-degree minus the in-degree is
nonnegative.  At a simple vertex with an incoming link, there is an
outgoing link, as just shown; so again the out-degree minus the
in-degree is nonnegative.  At a star, the in-degree is at most $1$, so
the out-degree minus the in-degree is at least $-1$.  Since the number
of stars is at most ${\mathcal S}$, the total out-degree minus the
total in-degree is at least
$$ \C - {\mathcal S} \ =\ 1. $$
On the other hand, it is zero, since every link has a head and a tail.
This contradiction shows that there is a relation $ja = pb + qc$ with
nonnegative integers $j,p,q$ and $j > 0$.

The relation $Jb = Pa + Qc$ is proved in the same way, using the graph
$\Gamma_b$ in place of $\Gamma_a$, under the hypothesis that there is
no relation $Jb = Pa + Qc$ with $J > 0$.  At a center, the three tiles
have their $b$ edges ending there, an odd number, so some ray from the
center carries a $b$ edge and an $a$ edge, and an outgoing link of
$\Gamma_b$ starts there.  At a simple vertex, the tiles with their
$\alpha$ and $\gamma$ angles there have their $b$ edges ending there,
exactly two, so an incoming link of $\Gamma_b$ is followed by an
outgoing one.  The rest of the argument is unchanged.
\end{proof}

\begin{theorem} \label{theorem:rationaltile120}
 Let $ABC$ be an isosceles (and not equilateral)
 triangle with base angles $\alpha$.
Suppose $ABC$ is  tiled by 
a tile $(\alpha,\beta, 2\pi/3)$ with $\alpha$ not a rational multiple 
of $\pi$.  Then the tile is rational.
\end{theorem} 

\begin{proof}
Suppose $ABC$ is tiled as in the theorem.  By Lemma~\ref{lemma:crelation}
and Lemma~\ref{lemma:abrelation}, there are nonnegative integers
$j,p,q,J,P,Q,K,L,M$, with $j$, $J$, and $K$ positive, such that
\begin{align*}
jc &= pa + qb, \\
Ja &= Pb + Qc, \\
\intertext{and}
Kb &= La + Mc.
\end{align*}
Dividing the first two by $c$ we have
\begin{align*}
j &= p(a/c) + q(b/c),\\
\intertext{and}
-Q &= -J(a/c) + P(b/c).
\end{align*}

\case{1} $q \neq 0$.  These equations can be solved for $a/c$ and
$b/c$, with rational solutions, provided the determinant $pP + Jq$ is
not zero.  Since $p, P \ge 0$, $J > 0$, and $q > 0$, we have
$pP + Jq \ge Jq > 0$.

\case{2} $q = 0$.  Then $jc = pa$.  Since $j > 0$ and $c > 0$, we have
$p > 0$, so $a/c = j/p$ is rational.  Dividing the third relation by
$c$, we get $K(b/c) = L(a/c) + M$; since $K > 0$, $b/c$ is rational.
\end{proof}

\subsection{The Diophantine Equation $c^2 = a^2 + b^2 + ab$}

Let $(a,b,c)$ be the sides of a triangle with angles $(\alpha,\beta,2\pi/3)$.
By the law of cosines, and since $\cos(2\pi/3) = -1/2$, we have
\begin{eqnarray*}
c^2 &=& a^2 + b^2 - 2ab \cos(2\pi/3) \\
&=& a^2 + b^2 + ab.
\end{eqnarray*}
Therefore, this Diophantine equation determines the possible rational 
triangles with a $2\pi/3$ angle. 

\begin{lemma}\label{lemma:relprime} Suppose $c^2 = a^2 + b^2 + ab$,
and $a,b$, and $c$ are integers with no common factor. Then $a,b$, and $c$ are
pairwise relatively prime.
\end{lemma}

\begin{proof}  If prime $p$ divides any two of $a,b$, and $c$, then
it also divides the third one.
\end{proof}

\begin{lemma} \label{lemma:2b+a}
Suppose $a,b$, and $c$  are integers with no common factor that are the sides of a triangle with angles $(\alpha,\beta,2\pi/3)$.  Then $2b+a$ is 
relatively prime to each of $a$, $b$, and $c$, except that if $a$
is even, 2 divides both $a$ and $2b+a$.
\end{lemma}

\begin{proof}   By Lemma~\ref{lemma:relprime}, $2b+a$
is relatively prime to $b$ and $a$, with the exception 
mentioned in the statement.  It remains to prove $2b+a$
is relatively prime to $c$.  Suppose, for proof by contradiction, 
that $p$ is a prime that divides both $c$ and $2b+a$.  Then $p$
is not 2, since then $c$ and $a$ would both be even, contradicting
Lemma~\ref{lemma:relprime}. 

Suppose, for proof by contradiction, that $p=3$.
Since $3$ divides $2b+a$, we have $a \equiv b \pmod{3}$.
Thus $a=b+3t$ for some integer $t$, and
\[
c^2 = a^2+b^2+ab = 3b^2+9bt+9t^2.
\]
Since $3$ divides $c$, we have $9 \mid c^2$.
The displayed equation therefore implies $3 \mid b^2$,
and hence $3 \mid b$. Since $a \equiv b \pmod{3}$,
we also have $3 \mid a$. This contradicts the fact
that $a$ and $b$ are relatively prime.
Hence, $p \neq 3$.

Then we have
\begin{align*}
c &\equiv 0 \pmod{p}, \\
c^2 &\equiv a^2 + b^2 + ab \pmod{p}, \\
\intertext{and}
2b+a &\equiv 0 \pmod{p}.
\end{align*}
Using $c \equiv 0 \pmod{p}$, we obtain
\begin{align}
0 &\equiv a^2 + b^2 + ab \pmod{p}, \label{eq:8806}\\
\intertext{and}
2b+a &\equiv 0 \pmod{p}. \notag
\end{align}
From the second equation we have $a \equiv -2b \pmod{p}$.
Since $a$ and $b$ are relatively prime, and $p \neq 2$,
this implies that neither $a$ nor $b$ is divisible by $p$.
Using $a \equiv -2b \pmod{p}$ in Equation~(\ref{eq:8806}), we have
$$
0 \equiv 4b^2 + b^2 - 2b^2 \equiv 3b^2 \pmod{p}.
$$
Since $p \neq 3$, it follows that $b^2 \equiv 0 \pmod{p}$, and hence
$b \equiv 0 \pmod{p}$, since $p$ is prime.
This contradicts the fact that $p$ does not divide $b$.
\end{proof}

{\em Remark}.  Using the techniques of \cite{cohen},
Corollary 6.3.15, page~353, we are able to parametrize the solutions
of $c^2 = a^2 + b^2 + ab$ by two integer parameters $(s,t)$ or 
one rational parameter $s/t$.  Having worked this out and used it
in preliminary versions,  we eventually found simpler proofs without 
it.  

\subsection{The Area Equation for an Isosceles Tiling with $\gamma = 2\pi/3$}

\begin{lemma} \label{lemma:2pi/3area}
Let isosceles triangle $ABC$ with base angles $\alpha$ be $N$-tiled
by a tile with angles $(\alpha,\beta,2\pi/3)$.
Suppose $\alpha$ is not a rational multiple of $\pi$.  
 Let $X$ be the length of the equal sides $AB$ and $BC$,
 and $Y$ the length of the base $AC$.  
Then the area equation is
$$ X^2 (2b+a) = Nbc^2. $$
and another form of the area equation is 
$$ XY = Nbc. $$
\end{lemma}

\begin{proof}  
 By the law of cosines,
$$ a^2 = b^2 + c^2 - 2bc \cos \alpha. $$
Hence, solving for $\cos \alpha$ and using $c^2 = a^2 + b^2 + ab$,
\begin{eqnarray}
\cos \alpha &=& \frac {b^2 + c^2 -a^2} {2 bc}\nonumber \\
&=& \frac {b^2 + (a^2 + b^2 + ab) - a^2}{2bc}\nonumber \\
&=& \frac {2b^2 + ab} {2bc} \nonumber\\
&=& \frac {2b+a}{2c}. \label{eq:8767}
\end{eqnarray}

Twice the area of $ABC$ is $X^2 \sin(\pi-2\alpha) = X^2 \sin 2\alpha$.
Twice the area of the tile is $bc\sin \alpha$.  Equating the area of $ABC$
to $N$ times the area of the tile, we have
$$ X^2 \sin 2\alpha = N bc \sin \alpha. $$
Hence, since $\sin 2\alpha = 2 \sin \alpha \cos \alpha$, dividing by $\sin \alpha$ we have
$$ 2X^2 \cos \alpha = Nbc. $$
Substituting for $\cos \alpha$ the value from Equation~(\ref{eq:8767}),
$$ 2X^2 \left(\frac {2b+a}{2c}\right) = Nbc. $$
Hence, multiplying by $c$,
$$ X^2 (2b+a) = Nbc^2. $$
That completes the proof of the first formula of the lemma.

To prove the second form:  twice the area of $ABC$ is 
$XY \sin \alpha$. Twice the area of the tile is $bc \sin \alpha$.
Therefore, $XY = Nbc$.  
\end{proof}

\subsection{A Necessary Condition}

\begin{lemma} \label{lemma:2b+adividesN}
Let isosceles triangle $ABC$ with base angles $\alpha$ be $N$-tiled
by a tile with angles $(\alpha,\beta,2\pi/3)$ and sides $(a,b,c)$.
Suppose $\alpha$ is not a rational multiple of $\pi$. 
Then 
\smallskip

(i) $2b+a$ divides $N$, and $Nb/(2b+a)$ is a square, say $m^2$, and
\smallskip

(ii) the side and base of $ABC$ are given by 
$$ X = mc \qquad\mbox{and}\qquad Y = m(2b+a). $$
\end{lemma}

\begin{remark}
This lemma gives us an {\em a priori} bound on $(a,b,c)$, namely $2N$, 
since $c^2 = a^2 + b^2 + ab \le (a+b)^2 \le (2b+a)^2 \le (2N)^2$.
Also, if $N$ is prime, $N = 2b+a$, and $b = m^2$.  
This actually can happen.  For example 
$(a,b,c)=(143,25,157)$ satisfies $2b+a = 193$, which is prime, and 
$b$ is a square.  It can even happen with a prime congruent to 
3 modulo 4: When $(a,b,c) = (39,16,49)$ we have $2b+a = 71$.
\end{remark}

\begin{proof}[Proof of Lemma~\ref{lemma:2b+adividesN}] Let $X$ be the length of the equal sides
$AB$ and $BC$.  According to Lemma~\ref{lemma:2pi/3area},
$$ X^2 (2b+a) = Nbc^2.$$
By Lemma~\ref{lemma:relprime}, $a$, $b$, and $c$ are pairwise
relatively prime.  By Lemma~\ref{lemma:2b+a},  if $a$ is odd, then
$2b+a$ is relatively prime to each of $a$, $b$, and $c$.
On the other hand, if $a$ is even, then $b$ and $c$ are odd,
so $2b+a$ is relatively prime to $c$ and $b$.  
Thus, whatever the parity of $a$, $2b+a$ is relatively 
prime to $b$ and $c$.
Then by the area equation, $2b+a$ divides $N$.

According to the area equation, 
$$ \left(\frac X c\right)^2 = \left( \frac {Nb}{2b+a} \right). $$
Therefore, $Nb/(2b+a)$ is a rational square, and since it is an integer,
it is an integer square, say $m^2$.   That completes
the proof of part~(i).

Proof of (ii).  Since $X/c$ and $m$ are
positive and have equal squares, they are equal, so $X = cm$ as
claimed. We compute $Y$.  
By the law of cosines, and since $c^2 = a^2+b^2+ab$,
$$
\cos \alpha = \frac{b^2+c^2-a^2}{2bc} = \frac{2b^2 + ab}{2bc}
            = \frac{2b + a}{2c}.
$$
Since $Y = 2X\cos\alpha$, where $X = \vert AB \vert$, we have
$$
Y = \frac{X}{c}\,(2b+a).
$$
By (i), $c$ divides $X$, so the right-hand side is an integer.
\end{proof}

\begin{example}  In the tiling whose construction begins with
Figure~\ref{figure:isosceles2pi2}, we have $N = 3681860$,
  $(a,b,c) = (3,5,7)$, so $2b + a = 13$, and $Nb/(2b+a) = 170^2 \cdot 7^2$,
so $m = \sqrt{Nb/(2b+a)} = 1190$,  $X =  mc =  8330$, and $Y = m(2b + a) = 15470$. 
\end{example}

\begin{example}
With $N=33$, $(a,b,c) = (5,3,7)$, $2b+a = 11$,
$Nb/(2b+a) = 9$, $m = 3$, $X = mc = 21$, and $Y = m(2b+a) = 33$.  
We think that no such tiling exists, although the 
present lemma does not rule it out.   In principle one can 
``just check all the possibilities'', but that is easier
said than done. 

By contrast, with the tile $(3,5,7)$, $Nb/(2b+a) = 165/13$ is not an 
integer, so no $33$-tiling with that tile is possible. \end{example}

\begin{example} With $N=37$, $(a,b,c) = (5,16,19)$, $2b+a=37$,
$Nb/(2b+a) = 16$, $m=4$, $X = mc = 76$, and $Y = m(2b+a) = 148$.
Since 37 is prime, no such tiling exists, by \cite{bonfioli2026};
the present lemma alone does not rule it out.
\end{example}

\subsection{Ruling Out Small Values of $N$}

\begin{theorem} \label{theorem:smallN120}
If there is an $N$-tiling of some isosceles triangle $ABC$ with
base angles $\alpha$ by a tile with angles $(\alpha,\beta,2\pi/3)$,
and $\alpha$ is not a rational multiple of $\pi$,
then $N$ is at least 33.  If $N \le 200$ then $N$ is one of the 
values shown in Table~\ref{table:isosceles120}, and the side
and base of $ABC$ must be as given in the table.
\end{theorem}

\begin{table}[ht]
\centering
\begin{tabular}{rrr}
 $N$ & $(a,b,c)$ & $(X,Y)$\\
 \hline
33 &(5, 3, 7) & (21, 33)\\
37 &(5, 16, 19)& (76, 148)\\
46 &(7, 8, 13)&(52, 92) \\
65 &(3, 5, 7)&(35, 65)\\
71 &(39, 16, 49)&(196, 284)\\
74 &(56, 9, 61)&(183, 222)\\
130& (16, 5, 19)&(95, 130)\\
132 &(5, 3, 7)&(42, 66)\\
148 &(5, 16, 19)&(152, 296)\\
154 &(8, 7, 13)&(91, 154)\\
184 &(7, 8, 13)&(104, 184)\\
193 &(143, 25, 157)&(785, 965)
\end{tabular}
\caption{$N < 200$ and $(a,b,c)$ not ruled out by Lemma~\ref{lemma:2b+adividesN}.}
\label{table:isosceles120}
\end{table}

\begin{remark}  We do not suggest that tilings for 
$N$ in the table do, or do not, exist, only that they are not ruled out
by the simple considerations of area and boundary tiling. 
\end{remark}

The prime numbers 37, 71, and 193 are not ruled out immediately, and 
$71 \equiv 3 \pmod 4$.
 Hence, the possibility of $N$ prime
for this kind of tiling is not ruled out by the area equation and 
boundary-tiling conditions; but at least the cases 7, 11, 19 
are eliminated, which is required for a proof that there are no 
$N$-tilings of any triangle for those values of $N$.

\begin{proof}[Proof of Theorem~\ref{theorem:smallN120}]
Let the positive integer $N$ be given, 
and suppose there is an $N$-tiling of some isosceles $ABC$
by a tile $(\alpha,\beta,2\pi/3)$.  By 
Theorem~\ref{theorem:rationaltile120}, the tile is 
rational, so we may suppose its sides are integers $a,b,c$
with no common divisor. 
According to Lemma~\ref{lemma:2b+adividesN},
$2b+a$ divides $N$ (so $a$ and $b$ are at most $N$),
 and $Nb/(2b+a)$ is a square, say $m^2$. 
  Then, since the 
 tile has a $2\pi/3$ angle, $c$ is determined by the equation
 $c^2 = a^2+b^2 + ab$.  If $c$ is not an integer, then 
 we do not consider $(a,b,c)$ further.  Also if $(a,b,c)$
 is not a triangle, because the sum of two of its sides is less
 than the third, we do not consider it further. 
 Table~\ref{table:isosceles120}
was computed by running this algorithm for $N \le 200$.
There are no entries for $N < 33$.
\end{proof}

We note that it would be a waste of time to compute the length 
of the base $Y$ and reject $(a,b,c)$ in case $Y$ is not an integer,
because $Y$ always {\em has} to be an integer $m(a+2b)$,
by Lemma~\ref{lemma:2b+adividesN}.

\subsection{Given $N$, Find the Possible Tiles and $ABC$}

\begin{theorem}\label{theorem:possibletiles120}
Given $N$, we can efficiently compute a finite set 
$\Delta$ of  $(a,b,c,X)$, such that if there is 
an $N$-tiling of some isosceles triangle $ABC$ with
base angles $\alpha$ by a tile 
with $\gamma = 2\pi/3$,  then the tile is $(a,b,c)$
and the side of $ABC$ is $X$, for some $(a,b,c)$ in $\Delta$.
\end{theorem}

\begin{remark} Then by backtracking search, 
applied to each tile $(a,b,c)$ and isosceles triangle
$ABC$ with base angles $\alpha$ and side $X$ 
with $(a,b,c,X)$ in $\Delta$, 
we can determine (in principle) if any $N$-tiling
of any isosceles $ABC$ exists.  We do not undertake that in this paper, 
because of the difficulty of
establishing the correctness of such computer searches beyond a shadow of doubt.
\end{remark}

\begin{proof} [Proof of Theorem~\ref{theorem:possibletiles120}]
 Let $N$ be given.  The
algorithm given in the proof of Theorem~\ref{theorem:smallN120}
determines the possible tiles $(a,b,c)$, in such a way that 
$Nb/(2b+a)$ is a square, say $m^2$.  Then $X = mc$ must be 
the side of triangle $ABC$, if there is any $N$-tiling 
of isosceles $ABC$ by $(a,b,c)$, and $Y = m(2b+a)$ is the 
base, by Lemma~\ref{lemma:2b+adividesN}.
\end{proof}

\section{Open Problems}
The methods and results of this paper
leave us still unable to answer some interesting questions.  Here 
we list several.   In the following, as elsewhere in this paper, 
``isosceles'' means ``isosceles and not equilateral.''

(i) What is the smallest $N$ such that some isosceles triangle with 
base angles $\alpha$ 
can be $N$-tiled by a tile of the form $\gamma = 2\alpha$?  The smallest
such tiling so far explicitly constructed has 720 tiles,
but for all we know there is a 45-tiling.%
\footnote{In this paper we make no judgement about claims by others
to have proved non-existence of certain tilings by exhaustive search.
}
In fact, we do not even know the smallest $N$ such that some isosceles
triangle can be tiled by the tile with sides $(4,5,6)$ (which is the 
tile used in the 720-tiling).

(ii) Find  easily checkable 
necessary and sufficient conditions on $N$ for the existence
of $N$-tilings of some isosceles $ABC$ with $\gamma = 2\alpha$
or $\gamma = 2\pi/3$.   Or, determine the existence or 
non-existence of  such tilings,  one $N$ at a time,
using exhaustive computer search.  You can use Table~\ref{table:isosceles120}
for your initial test data.

\section{Conclusions}

We have studied the possible tilings of an isosceles (and not equilateral)
triangle $ABC$ by 
a tile that is a right triangle,  or by a tile of the form 
$(\alpha,\beta,2\alpha)$ where the base angles of $ABC$ are equal to $\alpha$.
In the case of a tile $(\alpha,\beta, 2\alpha)$, we derived a necessary
condition from the area equation, and
we made use of 
directed graphs inspired by Laczkovich to prove that the tile is necessarily
rational.

We analyzed the case of a right-angled tile thoroughly enough to give a complete
characterization of the possible  values of $N$ for which 
some isosceles $ABC$ can 
be $N$-tiled.  Namely, Theorem~\ref{theorem:isosceles2} says $N$ is one of the following:
\begin{itemize}
\setlength{\itemsep}{0pt}
\setlength{\parskip}{0pt}
\item twice a square  
\item an even sum of squares
\item twice an even sum of squares
\item six times a square, with base angles $\tfrac \pi 6 =\beta$
\item any square,  with $ABC$ a right isosceles triangle and $\beta = \tfrac \pi 4$
\end{itemize}

In the case of a tile with $\gamma = 2\alpha$, we gave a necessary condition,
using the area equation and the law of cosines for the tile.  That this 
necessary condition is not trivial is shown by our proof that 
$N$ cannot be prime.  $N=45$ is the least number for which 
we do not know whether a tiling exists, and $720$ is
the smallest $N$ for which we are certain
that there does exist a tiling. 

Finally, in the case of a tile with $\gamma = 2\pi/3$, we gave a 
necessary condition and an algorithm to check it.  There
are no such tilings for $N < 33$.  There is one for $N = 2673$.
Between those two values of $N$, there are many values of $N$
satisfying our necessary conditions, for which we do not know whether
tilings exist.

In all the possible cases of Laczkovich's tables,
  we have been able to show 
(either in this paper or in unpublished work) that given $N$,
there is a finite set $\Delta$ of tiles $(a,b,c)$ and triangles $ABC$
such that either there is no $N$-tiling falling under that line
of the table, or one of the finite set permits an $N$-tiling.  Hence,
there is (in principle) an algorithm,  albeit inefficient, for determining
if there is an $N$-tiling.   The inefficiency arises from the exponentially
large number of ways of trying to place $N$ tiles of a specific shape into 
a specific $ABC$.

\section{Tilings Found by Bryce Herdt} \label{section:herdt}
In 2024, Bryce Herdt found several new tilings with tiles $(4,5,6)$ (for which $\gamma = 2 \alpha$)
and $(3,5,7)$ (for which $\gamma = 2\pi/3$).  These tilings dramatically lowered the $N$
for the ``smallest known tiling'', both of isosceles triangles and of equilateral triangles.
Here we exhibit Herdt's tilings.

\begin{figure}
\begin{center}
\includegraphics[width=0.4\textwidth]{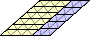}
\end{center}
\caption{Decomposing a parallelogram with top $pc + qb$ and side $bc$.}
\label{figure:bryce1}
\end{figure}

Consider Figure~\ref{figure:bigisosceles}.   The calculations after that figure give the 
number of tiles required on each edge in that figure.  All those numbers
are divisible by 6, except $m = 869$.  
Therefore, if we shrink the diagram by a (linear) factor of 6,  we can still tile all the triangles
in the figure with $(4,5,6)$, and the yellow parallelogram will be 869 by 480 (since it was
$576 \cdot 5$, after shrinking it is $576\cdot 5/6 = 480$. 
That 869 by 480 parallelogram cannot be 
tiled with tiles all in the same orientation.  As Herdt pointed out  
(in 2024), often a parallelogram can be divided into two smaller parallelograms,
which can be tiled with tiles in different orientations.   Figure~\ref{figure:bryce1}
illustrates the technique.  In the case at hand, $869 = 425+444 = 85 b + 74c$, 
while $480$ is divisible by $b$ and $c$ (5 and 6).  
The final value of $N$ will be $6028020/6^2 = 167445$.  That is, unfortunately, still 
too large to draw in a space smaller than several meters.   

\begin{figure}[ht]
\begin{center}
\includegraphics[width=\textwidth]{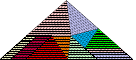}
\end{center}
\caption{Herdt's 1125-tiling by $(4,5,6)$.}
\label{figure:bryce3}
\end{figure}

Starting from a different dissection of the triangle into triangles and parallelograms,
Herdt was able to construct the $1125$-tiling shown in Figure~\ref{figure:bryce3}. Observe
how the parallelogram-dissection technique has been applied twice, in both the lower left
and the lower right corners of the tiling. 

\begin{figure}[ht]
\begin{center}
\includegraphics[width=\textwidth]{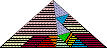}
\end{center}
\caption{Herdt's 720-tiling by $(4,5,6)$}
\label{figure:bryce5}
\end{figure}

Dissecting $ABC$ into similar triangles and parallelograms in a different way,
Herdt was able to reduce $N$ still further, to $720$, as shown in Figure~\ref{figure:bryce5}.

\begin{figure}[ht]
\begin{center}
\includegraphics[width=\textwidth]{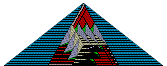}
\end{center}
\caption{Herdt's 2673-tiling by $(5,3,7)$. Here $\gamma = 2\pi/3$.}
\label{figure:bryce4}
\end{figure}

Herdt also found new tilings with the tile $(3,5,7)$, which has $\gamma = 2\pi/3$.
He began by finding a 1215-tiling of an equilateral triangle of side 135; this
was also done by his technique of decomposing parallelograms, starting with a known 
10935-tiling.  Then he flanked that equilateral triangle by two triangles 
similar to $(3,5,7)$, producing the tiling shown in 
 Figure~\ref{figure:bryce4}.  Here $N = 1215 + 2\cdot 27^2 = 2673$.
This is presently the smallest known tiling of an isosceles triangle by 
a tile with $\gamma = 2\pi/3$.
\FloatBarrier

\appendix
\section{Formal Verification in Lean}\label{section:lean}

The results of this paper have been formalized in the Lean~4 proof assistant
(version~4.30.0), using the Mathlib library.  The Lean code is at
\url{https://github.com/mjbeeson/IsoscelesTilingsLean} (tag {\tt v1.0-isosceles}).
Each theorem and lemma listed in Table~\ref{table:lean} has a Lean counterpart,
usually named after its label, whose proof is checked by the Lean kernel using
only the standard axioms ({\tt propext}, {\tt Classical.choice},
{\tt Quot.sound}), with no {\tt sorry} and no {\tt native\_decide}.
Tilings are formalized as in V.~Bonfioli's Lean library \cite{bonfioli-lean}:
$N$ triangles congruent to a given tile (reflections allowed), with disjoint
interiors, whose union is $ABC$.

\smallskip
\noindent{\em What is not formalized.}
\begin{itemize}
\item The existence of tilings: the constructions in the introduction and in
Section~\ref{section:herdt}, the ``possible for any such $N$'' clauses of
Theorem~\ref{theorem:isosceles2}, the ``if'' directions of
Theorem~\ref{theorem:whichABC}, and the sufficiency of the conditions in
Table~\ref{table:laczkovich-extended}.  Only the necessity of the stated
conditions is formalized.
\item Theorem~\ref{theorem:snover} as stated.  The part used in
Table~\ref{table:laczkovich-extended}, that if $N$ is not a square then the
tile is right-angled, is formalized.
\end{itemize}

\noindent{\em Results taken from Bonfioli's library.}
Laczkovich's classification (Theorem~\ref{theorem:laczkovich-isosceles}, from
Theorems~4.1, 5.1 and 5.3 of \cite{laczkovich1995} and Theorem~3.3 of
\cite{laczkovich2012}, for congruent tiles) and the part of
Theorem~\ref{theorem:snover} mentioned above are taken from
\cite{bonfioli-lean}.  The results that depend on them
(Theorems~\ref{theorem:laczkovich-isosceles}, \ref{theorem:isosceles2},
\ref{theorem:whichABC}, \ref{theorem:laczkovich-extended},
Lemma~\ref{lemma:baseangles} and Corollary~\ref{lemma:notprime}) are in a
separate part of the Lean code, which takes several hours to build; the
README of the repository explains how.

\smallskip
\noindent{\em Differences between the Lean proofs and the paper.}
A few Lean proofs take a different route.
Lemma~\ref{lemma:essentialsegments} is formalized only for the relation
$Ja = Ub + Vc$; the proof of Theorem~\ref{theorem:rationaltile} in Lean uses
that relation and Lemma~\ref{lemma:crelation4}.
Lemmas~\ref{lemma:centeroutgoing}, \ref{lemma:extendlink4} and
\ref{lemma:extendlink} are steps inside proofs that are formalized by other
arguments, and have no Lean counterpart of their own.

\begin{table}[ht]
\centering
\footnotesize
\begin{tabular}{ll}
Result & Lean name \\
\hline
Theorem~\ref{theorem:laczkovich-isosceles} & {\tt Laczkovich.laczkovich\_isosceles} \\
Lemma~\ref{lemma:sumsofsquares} & {\tt sumsofsquares} \\
Lemma~\ref{lemma:doublesquares} & {\tt doublesquares} \\
Lemma~\ref{lemma:euler} & {\tt euler} \\
Lemma~\ref{lemma:rationalsquares} & {\tt rationalsquares} \\
Lemma~\ref{lemma:Gamma} & {\tt Gamma} \\
Lemma~\ref{lemma:propagation} & {\tt propagation} \\
Corollary~\ref{lemma:allc} & {\tt allc} \\
Lemma~\ref{lemma:isosceleshelper2} & {\tt isosceleshelper2} \\
Lemma~\ref{lemma:isosceleshelper3} & {\tt isosceleshelper3} \\
Theorem~\ref{theorem:rationalisosceles} & {\tt rationalisosceles} \\
Lemma~\ref{lemma:baseangles} & {\tt baseangles} \\
Theorem~\ref{theorem:isosceles2} & {\tt isosceles2} \\
Corollary~\ref{lemma:notprime} & {\tt notprime} \\
Theorem~\ref{theorem:whichABC} & {\tt whichABC\_ii}, {\tt whichABC\_iv}, {\tt whichABC\_v} \\
Lemma~\ref{lemma:equilateral120} & {\tt equilateral120} \\
Theorem~\ref{theorem:laczkovich-extended} & {\tt laczkovich\_extended} \\
Lemma~\ref{lemma:twobeta} & {\tt twobeta} \\
Lemma~\ref{lemma:starsandcenters} & {\tt starsandcenters} \\
Lemma~\ref{lemma:essentialsegments} & {\tt arelation} (the relation for $a$) \\
Lemma~\ref{lemma:crelation4} & {\tt crelation4} \\
Theorem~\ref{theorem:rationaltile} & {\tt rationaltile} \\
Lemma~\ref{lemma:twoalphatiles} & {\tt twoalphatiles} \\
Lemma~\ref{lemma:luthar} & {\tt luthar} \\
Lemma~\ref{lemma:edgerelations} & {\tt edgerelations} \\
Lemma~\ref{lemma:areaequation} & {\tt areaequation} \\
Lemma~\ref{lemma:areaequation2} & {\tt areaequation2} \\
Lemma~\ref{lemma:notprimehelper} & {\tt notprimehelper} \\
Theorem~\ref{theorem:notsquarefree} & {\tt notsquarefree} \\
Theorem~\ref{theorem:finitetiles} & {\tt finitetiles} \\
Lemma~\ref{lemma:qandrnotzero} & {\tt qandrnotzero} \\
Lemma~\ref{lemma:no28} & {\tt no28} \\
Theorem~\ref{theorem:N45} & {\tt N45} \\
Lemma~\ref{lemma:crelation} & {\tt crelation} \\
Lemma~\ref{lemma:abrelation} & {\tt abrelation} \\
Theorem~\ref{theorem:rationaltile120} & {\tt rationaltile120} \\
Lemma~\ref{lemma:relprime} & {\tt relprime} \\
Lemma~\ref{lemma:2b+a} & {\tt two\_b\_add\_a} \\
Lemma~\ref{lemma:2pi/3area} & {\tt two\_pi\_div\_three\_area} \\
Lemma~\ref{lemma:2b+adividesN} & {\tt two\_b\_add\_a\_dvd\_N} \\
Theorem~\ref{theorem:smallN120} & {\tt smallN120} \\
Theorem~\ref{theorem:possibletiles120} & {\tt possibletiles120} \\
\end{tabular}
\caption{Results of this paper and their Lean names (all in the namespace
{\tt TriangleTiling}).}
\label{table:lean}
\end{table}

\vskip20pt\noindent {\bf Acknowledgement.} 
We are grateful to Miklos Laczkovich for his valuable comments on our work and 
especially for his simplification of the proofs of 
Lemmas~\ref{lemma:Gamma},  
\ref{lemma:isosceleshelper2}, and \ref{lemma:isosceleshelper3}, 
and of course for his many pioneering papers in this subject, on which
this paper rests.

\vskip 20pt
\noindent {\bf Acknowledgment.} The original version of this paper was
written without AI assistance.  In preparing this revision, the author
used Claude (Anthropic), an AI system, in three ways.  First, to check
the proofs: it found a number of errors and gaps, and proposed corrected
statements and proofs of several lemmas.  Second, it proposed new and
shorter proofs of some results, notably Theorem~\ref{theorem:N45}, and
Theorem~\ref{theorem:finitetiles}, which replace a longer argument in an
earlier version.  Third, it wrote the Lean formalization described in
Appendix~\ref{section:lean}, whose proofs are checked by the Lean kernel.  The author
checked all proposed changes to the mathematics, rejected or modified
some, checked that the Lean statements say what the theorems of the paper
say, and is responsible for the final text.  The main results are unchanged.
 

\begin{thebibliography}{10}\footnotesize

\bibitem{beeson-noseven} {\raggedright M. Beeson, No triangle can be cut into seven congruent triangles, preprint, {\tt arXiv:1811.09723}.\par}  


\bibitem{bonfioli2026}{\raggedright V. Bonfioli,
\newblock Congruent dissections of a triangle: prime counts, the counts up to 200,
  and a family of isosceles targets, Theorem~1.4 ({\tt thm:iso}).
\newblock Lean~4 formalization, GitHub repository
  \texttt{ElVec1o/erdos634\_primes\_and\_N200}, release v4.0, 22 October 2026,
  {\tt doi:10.5281/zenodo.23400433}.\par}
  
\bibitem{bonfioli-lean}{\raggedright V. Bonfioli,
 \newblock Formalizing triangle tiling in Lean, version 1.0, 22 October 2026.
\newblock Lean~4 formalization, GitHub repository
  \texttt{ElVec1o/formalizing-triangle-tiling-in-lean}.\par}

\bibitem{laczkovich2010} S. Butler, F. R. K. Chung, R. L. Graham, and M. Laczkovich, Tiling polygons with lattice triangles, {\it Discrete Comput. Geom.} {\bf 44} (2010), 896--903.

\bibitem{cohen} H. Cohen, {\it Number Theory, Volume 1: Tools and Diophantine Equations}, Springer, New York, 2007.

\bibitem{golomb} S. W. Golomb, Replicating figures in the plane, {\it Math. Gaz.} {\bf 48} (1964), 403--412.

\bibitem{hardy-wright} G. H. Hardy and E. M. Wright, {\it An Introduction to the Theory of Numbers}, fourth edition, Clarendon Press, Oxford, 1960.

\bibitem{laczkovich1990} M. Laczkovich, Tilings of polygons with similar triangles, {\it Combinatorica} {\bf 10} (1990), 281--306.

\bibitem{laczkovich1995} M. Laczkovich, Tilings of triangles, {\it Discrete Math.} {\bf 140} (1995), 79--94.

\bibitem{laczkovich1998} M. Laczkovich, Tilings of polygons with similar triangles, II, {\it Discrete Comput. Geom.} {\bf 19} (1998), 411--425.

\bibitem{laczkovich2012} M. Laczkovich, Tilings of convex polygons with congruent triangles, {\it Discrete Comput. Geom.} {\bf 48} (2012), 330--372.

\bibitem{laczkovich-szekeres1995} M. Laczkovich and G. Szekeres, Tilings of the square with similar rectangles, {\it Discrete Comput. Geom.} {\bf 13} (1995), 569--572.

\bibitem{luthar} R. S. Luthar, Integer-sided triangles with one angle twice another, {\it College Math. J.} {\bf 15} (1) (1984), 55--56.

\bibitem{niven} I. Niven, {\it Irrational Numbers}, Carus Mathematical Monographs 11, Mathematical Association of America, Washington, DC, 1967.

\bibitem{snover1991} S. L. Snover, C. Waiveris, and J. K. Williams, Rep-tiling for triangles, {\it Discrete Math.} {\bf 91} (1991), 193--200.

\bibitem{soifer} A. Soifer, {\it How Does One Cut a Triangle?}, Springer, New York, 2009.

\end{thebibliography}

\end{document}